\documentclass[10pt]{template}
\usepackage[utf8]{inputenc}
\usepackage{amsmath,amsthm}
\usepackage{rotating}
\usepackage{amssymb}
\usepackage{multirow}
\usepackage{bbm}
\usepackage{verbatim}
\usepackage{lineno}
\usepackage{calligra}
\usepackage{subfig}
\usepackage{setspace}
\usepackage{cleveref}
\usepackage{graphicx}
\usepackage{booktabs}
\usepackage[table,xcdraw,dvipsnames]{xcolor}
\usepackage[backend=biber, sorting=none, maxbibnames=99, style=vancouver ]{biblatex}
\usepackage{geometry}
\usepackage{float}

\onehalfspace
\bibliography{main.bbl}
\theoremstyle{remark}
\newtheorem*{remark}{Remark}

\newtheorem*{ocprob}{Optimal Control Problem}

\newcommand{\beginsupplement}{%
        \setcounter{section}{0}
        \renewcommand{\thesection}{S\arabic{section}}
        \setcounter{table}{0}
        \renewcommand{\thetable}{ST\arabic{table}}%
        \setcounter{figure}{0}
        \renewcommand{\thefigure}{SF\arabic{figure}}%
     }

\title{ Optimized numerical solutions of SIRDVW multiage model controlling SARS-CoV-2 vaccine roll out: an application to the Italian scenario.}
\author{Giovanni Ziarelli\textsuperscript{1}, Luca Dede'\textsuperscript{1}, Nicola Parolini\textsuperscript{1}, Marco Verani\textsuperscript{1}, Alfio Quarteroni\textsuperscript{1,2}}
\address{$^{1}$MOX, Department of Mathematics, Politecnico di Milano, Italy\\$^{2}$Institute of Mathematics, Ecole Polytechnique Fédérale de Lausanne (EPFL), Switzerland (Professor Emeritus)}
\date{October 2022}

\begin{document}
\maketitle
\section*{Abstract}
In the context of SARS-CoV-2 pandemic, mathematical modelling has played a fundamental role for making forecasts, simulating scenarios and evaluating the impact of preventive political, social and pharmaceutical measures. Optimal control theory can be a useful tool based on solid mathematical bases to plan the vaccination campaign in the direction of eradicating the pandemic as fast as possible. The aim of this work is to explore the optimal prioritisation order for planning vaccination campaigns able to achieve specific goals, as the reduction of the amount of infected, deceased and hospitalized in a fixed time frame, among age classes. For this purpose, we introduce an age stratified $SIR$-like epidemic compartmental model settled in an abstract framework for modelling two-doses vaccination campaigns and conceived with the description of COVID19 disease. Overall, we formalize an optimal control framework adopting the model as state problem by acting on the administrations of vaccine-doses.
An extensive campaign of numerical tests, featured in the Italian scenario and calibrated on available data from Dipartimento di Protezione Civile Italiana, shows that the presented framework can be a valuable tool to support the planning of vaccination campaigns minimizing specific goals.

\section{Introduction}
In January 2020 the tremendous SARS-CoV-2 virus (the agent pathogen of the COVID19 disease) outbroke in the Chinese province of Hubei, with epicenter the city of Wuhan. The first detected infections in Italy date back to the 21st February, when two distinct cases have been detected in Veneto and Lombardy regions. From this point onward up to September 2022 more than 602 millions cases and 6.4 millions deaths have been recorded around the world according to the weekly reports available from the World Health Organization (WHO) \cite{who2022}, which has declared the pandemic alert since March 2020. Since the severe symptomatology and highly--transmissible nature of the disease, Public Health Authorities of many countries have reacted tempestively in order to minimize the infectious risk introducing Non-Pharmaceutical Interventions (so-called NPIs) such as the compulsory adoption of face-masks, hygienic precautions, several measures for minimizing contacts, imposing smart-working whenever possible and sometimes with different kinds of lockdown. This is the Italian scenario, where starting from March 2020 to June 2020 a strict lockdown has been imposed at all levels. Non-pharmaceutical Interventions have been a fundamental tool in order to minimize the spread of the virus and the generation of severer variants. However, the virus has overtaken such preventive measures, letting many different variants arising (especially some of those particularly dangerous to be catalogued as Variants of Concern, VOC, by WHO). We recall the main names of the strains which have been prevalent in the Italian country along with the wild-type: B.1.617.2 (Delta), B.1.117 (Alpha), P.1 (Gamma) and B.1.1.529 (Omicron). The different restrictions have been the sole tool for facing the pandemic until December 2021, when the first vaccine has been approved by the American Food and Drugs administration (FDA) and by the European Medical Agency (EMA). Starting from January 2021 until September 2022 six different vaccines developed with different pharmaceutical techniques have been approved and employed in the Italian vaccination campaign: Pfizer mRNABNT162b2 (Comirnaty), COVID-19 Vaccine Moderna mRNA--1273 (Spikevax), Vaxzevria, Jcovden, Nuvaxovid (Novavax) and Valneva. Some of these vaccines require two administrations in order to acquire effectiveness against transmission and severe symptoms. One of the main difficulties that Italian Public Health authorities had to face during the planning stage of the vaccination campaign was the prioritization order across ages (at least for those for which the vaccine administrations had been safely tested) and working categories, alongside with the choice of the suitable elapsing time between subsequent administrations.

In this complex social and medical scenario, we set this work aiming at contributing to the knowledge acquired in the mathematical epidemiology field. In particular, the mathematical community has contributed in this scenario by employing methods and models in order to study the current pandemic situation in different respects, \textit{e.g.} \cite{bertuzzo2020geography, gatto2020spread, chinazzi2020effect, kuhl2021computational, parolini2021dashboard}. Typical approaches for modelling and forecasting scenarios consider Compartmental \cite{bertozzi2020challenges, bertaglia2021hyperbolic, rozhnova2021model, viana2021controlling, parolini2021suihter, capistran2022filtering} and Agent-Based models \cite{wolfram2020agent, gharakhanlou2020spatio, kerr2021covasim, shamil2021agent}. The former approach adopts dynamical systems for modelling all the necessary features of the disease, and can be easily adapted to available data for the calibration of the possibly many parameters involved. The latter, instead, is particularly focused on capturing the behavior of and interactions among individuals in specific exposure contacts (\textit{e.g.} schools). In this work, we recast the main questions regarding the planning of the vaccination campaign in the framework of the Optimal Control Theory \cite{kirk2004optimal}. Some recent works have already investigated the use of optimal control techniques based on compartmental models, to act on the levels of NPIs to be implemented for minimizing infected \cite{lemecha2020optimal, araz2021analysis, tsay2020modeling, ZAMIR2020105642} or deceased individuals \cite{richard2021age, perkins2020optimal}, sometimes coupling the evolution of the states of infectiousness model with other opinion models as in \cite{silva2021optimal}. Other works have dealt with the optimal allocation of vaccines, based on $SIR$-like models as in the case of \cite{libotte2020determination, ziarelli2021numerical} or taking into account the spatial heterogeneity \cite{lemaitre2022optimal}, as well as age-dependancies \cite{shim2021optimal}.
However, a detailed classification and comparison of vaccination policies reducing infections, deaths or hospitalisations dependently on age is missing. 
Moreover, up to our knowledge, this is the first work providing an accurate analysis of optimal vaccination campaigns for SARS-CoV-2 vaccines dependently on age and considering a state model that can be straight-forwardly calibrated with data that are commonly available in many countries.
In particular, the present work displays optimally controlled vaccination policies under state constraints of an age-stratified compartmental model in an abstract framework which is then applied to the specific Italian case study during the first half of 2021, using real data available daily from the Dipartimento di Protezione Civile Italiana \cite{dpcdata}. An extensive campaign of numerical tests shows that the presented optimal control framework can represent a valuable tool to support, starting from real epidemiological data, the planning of vaccination campaigns aiming at achieving specific goals (e.g. reduction of deceased, infected, hospitalized).

The outline of the paper is as follows. In Section \ref{Methods} we illustrate the rationales behind the chosen compartmental model, detailing the specific features that allow to describe the COVID19 context. Moreover, we introduce the mathematical formulation of the optimal control problem and the main numerical technicalities. In Section \ref{Results} we apply the abstract framework to the Italian scenario highlighting the versatility and the potentialities of the implemented computational architecture. Lastly, in Section \ref{Conclusions} the main conclusions deriving from the analysis are drawn.
\section{Methods}
\label{Methods}
{In this section, we present the mathematical description of an optimal control problem aiming at the planning of a vaccination campaign which is able to achieve specific goals. More precisely, in Section \ref{MathematicalModel} we introduce the  epidemiological differential model governing the control problem, while in Section \ref{FormulationOC} we introduce the optimal control problem. Finally, we provide details about the implemented numerical scheme in Section \ref{NumericalProcedure}.}
\subsection{Mathematical Model}
\label{MathematicalModel}

\begin{figure}
    \centering
    \includegraphics[scale=0.4]{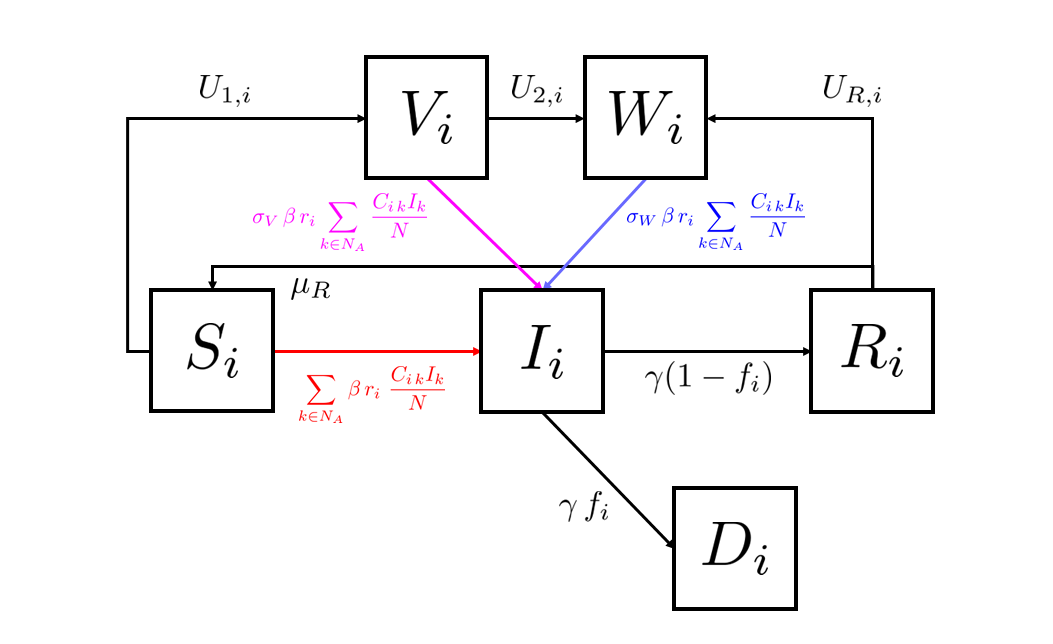}
    \caption{Schematic flowchart of each age-stratification of the adopted $SIRDVW$ multiage model. The colored fluxes account for new infections and they explicitly embed interactions among age-classes.}
    \label{flowchart}
\end{figure}
{We introduce a deterministic age-stratified compartmental model embodying the main features of SARS-CoV-2 transmission and vaccination campaign, cf. \cite{ram2021modified, parolini2022modelling, ivorra2020mathematical}.} 
{More precisely, we consider (see also the flowchart in Figure \ref{flowchart}) the following system of ODEs on the time interval $I=(0,T_f]$

\begin{equation}
    \begin{cases}
    \begin{aligned}
    &\dot{S_i}=- \beta\, r_i \displaystyle \sum_{k \in N_A}   \dfrac{S_i \, C_{ik} \, I_k}{N_i}  - U_{1,i} \dfrac{S_i}{S_i + I_{u,i}} + \mu_R R_i,\\[2pt]
    &\dot{I_i} = \beta\, r_i \displaystyle \sum_{k \in N_A} \dfrac{(S_i + \sigma_V V_i + \sigma_W W_i) \, C_{ik} \, I_k}{N_i} - \gamma I_i,\\[2pt]
    &\dot{R_i} = (1 - f_i(S_i,V_i,W_i))\, \gamma I_i - U_{R,i} - \mu_R R_i, \\[2pt]
    &\dot{D_i} = f_i(S_i,V_i,W_i)\, \gamma I_i, \\[2pt]
    &\dot{V_i} = - \beta\, r_i \displaystyle \sum_{k \in N_A} \dfrac{\sigma_V V_i \, C_{ik} \, I_k}{N_i} + U_{1,i} \dfrac{S_i}{S_i + I_{u,i}} - U_{2,i}, \\[2pt]
    &\dot{W_i} = - \beta\, r_i \displaystyle \sum_{k \in N_A} \dfrac{\sigma_W W_i \, C_{ik} \, I_k}{N_i}  + U_{2,i} + U_{R,i},\\
    &S_i(0) = S_{i,0},\, I_i(0) = I_{i,0},\, R_i(0) = R_{i,0}, \, D_i(0) = D_{i,0}, \, V_i(0) = V_{i,0},\, W_i(0) = W_{i,0},
    \end{aligned}
    \end{cases} \hspace{-1.5cm} \forall t\, \in \mathrm{I},\, i \in \, N_A, 
    \label{ODEmodel}
\end{equation}
where, as usual, each time dependent variable accounts for the number of individuals in different conditions with respect to the disease:   Susceptible ($S$), Infectious ($I$), Recovered ($R$), Deceased ($D$) due to complications related to SARS-CoV-2, Vaccinated with a first dose administered ($V$) and Vaccinated who have completed the vaccination cycle ($W$). Actually, each state is split according to the index $i \in N_A$, where $N_A$ is the set of indexes identifying the different age-classes.}
The parameters involved in the model are described in the following:
\begin{itemize}
    \item $T_f$ (days): final time of the simulation frame;  
    \item $\beta \in (0,1)$: transmission rate, depending on the implemented Non-Pharmaceutical Interventions (NPIs) and virus transmissibility. It is assumed to be constant across all ages as in \cite{marziano2021retrospective};
    \item $\sigma_V, \sigma_W \in (0,1)$: vaccine effectiveness on transmissibility after administration of first dose (the former) or completing the cycle (the latter). It can be interpreted as the ratio of transmissibility between vaccinated individuals and unvaccinated ones. The value 0 means that the vaccine is fully effective, 1 totally ineffective;
    
    \item $\theta_V, \theta_W \in (0,1)$: vaccine effect on mortality after administration of first dose (the former) or completing the cycle (the latter). It can be interpreted as the ratio of probability of getting severe symptoms between vaccinated individuals and unvaccinated ones. The value 0 means that the vaccine is fully effective, 1 totally ineffective;
    \item $IFR_i$: age-dependant Infection Fatality Rate;
    \item The fatality function $f_i(S_i,V_i,W_i)$ changes accordingly to the immunological profile of the population, \textit{i.e.} taking into account the vaccination effect in reducing mortality:
    \begin{equation}
        f_i(S_i,V_i,W_i) =  \begin{cases} 
        IFR_i  \left ( \dfrac{S_i(t-t_{a}) + \theta_V \sigma_V V_i(t-t_{a})+ \theta_W \sigma_W W_i(t-t_{a})}{S_i(t-t_{a}) + \sigma_V V_i(t-t_{a})+ \sigma_W W_i(t-t_{a})} \right) & t > t_a \\
      IFR_i & 0\leq t \leq t_a \\
   \end{cases}
    \end{equation}
    where $t_{a}$ is the amount of days after the inoculation that we consider for reaching the complete vaccine effectiveness. Heuristically, the fatality function reduces the fatality rate $IFR_i$ of a proper factor that takes into account the reduced probability of getting severe symptoms when vaccines have been inoculated (this effectiveness is ruled by the parameters $\theta_V,\theta_W$);
    \item $C_{ik}$: $i,k$-th entry of the contact matrix, tracing back contacts  between ages starting from the POLYMOD surveys \cite{10.1371/journal.pmed.0050074};
    \item $r_i$: susceptibilities to infection depending on age;
    \item $\gamma$: recovery rate from the disease infection, which is maintained constant across ages. Since infectious individuals are supposed to exit from the correspondent compartment  with  flux $\gamma I_i$, the parameter $\gamma$ is interpreted as the inverse of the average time of recovery $t_R$;
    \item $N_i$: {Number} of individuals in the i-th age stratification;
    \item $\mu_R$: natural waning immunity rate, taking into account plausible reinfections coming from previously-recovered individuals;
    \item $I_{u,i}(t) = (1 - \delta(t))\, I_i(t),\, t \, \in \mathrm{I}$: approximate number of undetected individuals whose age falls in the $i$-th stratification at time $t$;
    \item $U_{1,i},\, U_{2,i},\, U_{R,i}$: daily amount of administered first doses, second doses and doses administered to the $i$-th age-class, respectively. The choice of these variables is coherent with actual implementation of the italian vaccination campaign: two consecutive doses to be administered for completing the vaccination cycle to Susceptible individuals, one single administration to Recovered ones. 
    {To reduce the computational complexity of the model we assume that the functions $U_{1,i}(t),\, U_{2,i}(t),\, U_{R,i}(t)$ are piecewise constant (constant on each week) and the weekly value of administrations  is supposed to be equally distributed among each day of the week.} 
\end{itemize}
{
\begin{remark}
Since the  horizon of interest of our simulations is reasonably short, we neglect migratory effects, births and deaths which are not COVID-related. Moreover, we neglect the possible presence of comorbidities fostering the onset of severe symptoms leading to death. Finally, there are other implicit assumptions that is worth underlining: there is no genetic mutation of the virus during the period of interest; vaccine effectiveness waning is not contemplated in the considered time interval. 
\end{remark}

{\color{black} The differential problem \eqref{ODEmodel} has been endowed with proper initial conditions for each of the considered age-state compartment. Throughout this work, we consider time-invariant parameters except for the transmission rate $\beta$ and the detection rate $\delta$, which are assumed to depend on time. Especially for the transmission rate, which embodies many different effects, carrying out a calibration process is of paramount importance when dealing with realistic scenarios. In Appendix B one can find the calibration procedure employed for the Italian setting during the first half of 2021. On the other hand, more details about the choices of the parameters and the initial conditions can be retrieved in Appendix A. 


}
{
\subsection{The Optimal Control problem}
\label{FormulationOC}
{In this section we assume the perspective of Public Health Authorities asking whether it is possible to plan an optimal vaccination campaign in order to minimize some specific goals (total number of infected, deceased or hospitalized) and understand  possible priority orders among age--classes}. 
Let us specify the main features of the Optimal Control Problem: 
\begin{itemize}
    \item[] \textbf{Control Variables}. We assume that we could control the weekly amount of doses to be distributed among susceptibles and individuals who already got a first dose, and each day of the same week we assume to distribute an equal amount of doses. The control variables are 
    \begin{equation}
        0 \leq U_{1,i}(t), \, U_{2,i}(t) \leq N_{week}(t), \; \forall i  \, \in \, N_{A}, \, t \, \in \, \mathrm{I}. 
    \end{equation}
    In particular,
    \begin{equation}
    \displaystyle\{U_{1,i}(t)\}_{i}, \, \displaystyle \{U_{2,i}(t)\}_{i}, \, N_{week}(t) \in  \mathbb{P}^0_c := \biggl \{ f(t) = \displaystyle \sum_{k = 0}^{N_{int}} a_k \mathbbm{1}_{[T_{k-1}, T_{k}]}(t), \displaystyle \{T_k\}_{k\leq N_{int}} \in \mathbb{R},\, a_k \in \mathbb{R},\, \forall t \in  \mathrm{I} \biggr \},
    \end{equation}
letting $N_{int} \in \mathbb{N}$ the number of weeks and $\displaystyle \{T_k\}_{k\leq N_{int}}$ the set of first days of each week. In Section \ref{Results}, unless otherwise specified the piecewise constant function $N_{week}(t)$ has been obtained summing \textit{per-weeks} the amount of actually administered doses in Italy during the {considered period}. {Moreover, for reducing the computational complexity of the control problem and to adhere as much as possible with typical realistic scenarios (e.g. the Italian one)}, we assume that the numbers of administered second doses are set equal to the number of first administrations with a time delay imposed by pharmaceutical properties of the vaccine, \textit{i.e.}
    \begin{equation}
        U_{2,i}(t) = U_{1,i}(t - \delta_w),\; \forall t \in (\delta_w, T_f), \, i \in N_{A},
        \label{constDoses}
    \end{equation}
    where $\delta_w$ is the elapsing time among subsequent administrations.
    
    Finally, the control variables need to fulfill a budget constraint imposed by the limited amount of available doses and due to other potential sanitarian restrictions, such as the limited personnel and infrastructures. In particular, this constraint reads as: {$\forall j=0, \ldots, N_{int}-1$}
    \begin{equation}
        \displaystyle \sum_{i \in N_A} \int_{T_j}^{T_{j+1}} \left (U_{1,i}(t) + U_{2,i}(t) + U_{R,i}(t) \right ) \, dt\leq \mathrm{min}(N_S, N_{week}(T_j)),
        \label{budgConstraint}
    \end{equation}
    with $N_S$ is the budget limit due to sanitarian capabilities. 
    \item[] \textbf{State Problem}. The state problem is a revised version of \eqref{ODEmodel} under assumption \eqref{constDoses} on the administrations of first and second doses. For each age class $i \, \in \, N_{A}$ and given initial conditions, the problem at each time instant $t \in \mathrm{I}$ reads as \eqref{ODEmodel} in which \eqref{constDoses} has been plugged in. 
    In a compact way, the state problem can be rewritten as 
    \begin{equation}\label{pb:compact}
        \dot{\textbf{x}}_i(t) = \textbf{F}_i(\textbf{x},t,U_{1,i}) = \textbf{f}_i(\textbf{x},t) +  \widehat{\textbf{b}}\, U_{1,i}(t) + \widetilde{\textbf{b}}\, U_{1,i}(t - \delta_m) \;\; \forall i \in N_A,\, t \in (0,T_f],
    \end{equation}
    with
    \begin{equation}
    \textbf{x}_i(t) = \begin{bmatrix} S_i(t) \\ I_i(t) \\ R_i(t) \\ D_i(t) \\ V_i(t) \\ W_i(t)
    \end{bmatrix}, \;
    \textbf{x}(t) = \begin{bmatrix}
    \textbf{x}_1(t) \\ \textbf{x}_2(t) \\ \textbf{x}_3(t) \\ \textbf{x}_4(t) \\ \textbf{x}_5(t)
    \end{bmatrix},\;
        \widehat{\textbf{b}} = \begin{bmatrix} -1\\0\\0\\0\\1\\0
\end{bmatrix}\;\; \textrm{and} \;\; \widetilde{\textbf{b}} =\begin{bmatrix}0\\ 0\\ 0\\ 0\\-1\\1 \end{bmatrix}.
    \end{equation}
    This compact form will be useful for the computation of the gradient to be adopted in the numerical optimization process. 
    \item[] \textbf{Cost functionals}. 
    Assuming the perspective of policy makers, we consider three different cases:
    \begin{enumerate}
        \item 
        \begin{equation}
            \mathcal{J}_D(\textbf{x}) \, = \, \sum_{i  \in N_A} \int_{0}^{T_f} D_i(t)^2 \, dt,
        \label{JD}
        \end{equation}
        \textit{i.e.} the number of total deceased individuals during the whole process;
        \item 
        \begin{equation}
            \mathcal{J}_I(\textbf{x}) \, = \, \sum_{i  \in N_A} \int_{0}^{T_f} I_i(t)^2 \, dt,
        \label{JI}
        \end{equation}
        \textit{i.e.} the number of total infected individuals during the whole process;
        
        \item 
        \begin{equation}
            \mathcal{J}_H(\textbf{x}) \, = \, \sum_{i  \in N_A} \int_{0}^{T_f} H_i(t)^2 \, dt,
            \label{JH}
        \end{equation}
        \textit{i.e.} the number of total hospitalized individuals during the whole process. Indeed, one of the main issues to be faced during the first waves of the pandemic was dictated by the limited amount of beds and medical equipment for individuals affected by COVID19 exhibiting severe symptoms, especially in Intensive Care Units (ICU). However, we need an explicit computation of the hospitalized individuals as in \cite{marziano2021retrospective}, since the $SIRDVW$ does not provide directly the amount of hospitalized individuals as in \cite{parolini2022modelling}. In our case, for each age class we determine the amount of hospitalized individuals as
        \begin{equation}
            H_i(t) = h\, \kappa_i I_i(t) \left ( \dfrac{S_i(t-t_{a}) + \theta_V \sigma_V V_i(t-t_{a})+ \theta_W \sigma_W W_i(t-t_{a})}{S_i(t-t_{a}) + \sigma_V V_i(t-t_{a})+ \sigma_W W_i(t-t_{a})} \right ), 
        \end{equation}
        where $\kappa_i$ is the age-dependent propensity for severe respiratory symptoms and $h$ is an estimated fraction of hospitalized individuals with respect to the amount of infected (values retrieved from \cite{marziano2021retrospective}). 
    \end{enumerate}
\end{itemize}
{We are now ready to formulate the problem of the optimal vaccination campaign as an optimal control problem.}\\
\begin{ocprob}
Find the doses administrations $\widehat{U}_{1,i} \in \mathbb{P}_c^0 \; \forall i \in N_A$ such that it minimizes the chosen cost functional  $\mathcal{J}_X$, \textit{i.e.}
\begin{equation}\label{pb:OCP}
    \{\widehat{U}_{1,i}\}_{i \in N_{A}}  = \underset{\{U_{1,i}\}_{i \in N_A}\in \mathbb{P}^0_c}{\mathrm{argmin}}\mathcal{J}_X(\textbf{x}),
\end{equation}
subject to the state problem \eqref{pb:compact} under the constraints \eqref{constDoses} and \eqref{budgConstraint}.
\end{ocprob}
\subsection{Numerical procedure}
\label{NumericalProcedure}
To solve the optimal control problem  \eqref{pb:OCP} we employ a Projected Gradient Descent algorithm \cite{calamai1987projected} whose main steps are hereafter summarized:
\begin{itemize}
    \item[] \textbf{Choice of the initial guess}. Set the initial guess of the control functions.
    \item[] \textbf{Optimization cycle}. Fix a number of maximum iterations and at each iteration perform the following steps :
    \begin{enumerate}
        \item Solve the direct problem with the current control variable; 
        \item {Solve the adjoint time-reverse problem
        \begin{equation}
            \dot{\textbf{p}}_i(t) = - \dfrac{\partial \mathcal{H}_X}{\partial \textbf{x}_i}(t), \; \textbf{p}_i(T_f)|_k = -2 X_i(T_f) \delta_{k,X} \;\; \forall i \in N_A
            \label{eq:hamiltonianSys}
        \end{equation}
        with 
        \begin{equation}
            \mathcal{H}_X(t) := \displaystyle\sum_{i  \in N_A}X_i^2(t) + \textbf{p}_i(t)^T \textbf{F}_i(\textbf{x}, U_{1,i}),
        \end{equation}
        where $X \, \in \, \{D,\, I,\, H \}$ , $\textbf{p}_i$ the adjoint vector referring to the $i$-th age class.} 
        \item Compute the descending direction for the numerical optimization scheme:
        \begin{equation}\label{eq:descend}
         \delta U_{1,D_i}(t)= -(\textbf{p}_i(t) \widehat{\textbf{b}} + \textbf{p}_i(t + \delta_m) \widetilde{\textbf{b}} \mathbbm{1}_{[0; T-\delta_m]}),\, \forall t \in \mathrm{I}.
        \end{equation}
        In the sequel we briefly motivate this choice. Following the KKT conditions \cite{wright1999numerical}, define the Lagrangian function
        \begin{equation}
            \mathcal{L}_h(\textbf{x}, \dot{\textbf{x}}, \{U_{1,i}\}_{i\in N_A}, \textbf{p}) = \sum_{i  \in N_A} \mathcal{L}_{h,i}(\textbf{x}_i, \dot{\textbf{x}}_i,U_{1,i}, \textbf{p}_i),
        \end{equation}
        where,
        \begin{equation}
            \mathcal{L}_{h,i}(\textbf{x}_i, \dot{\textbf{x}}_i,U_{1,i}, \textbf{p}_i) = \int_0^T X_i(t)^2 \, dt + \int_0^T \textbf{p}_i^T(t) \, (\textbf{f}_i(\textbf{x}, t) + \widehat{\textbf{b}} \,U_{1,i}(t) + \widetilde{\textbf{b}} \, U_{1,i}(t - \delta_m) - \dot{\textbf{x}} ) \;dt.
        \end{equation}
    Compute the gradient of the Lagrangian at the continuous level in order to compute a descending direction for the optimal control algorithm. We obtain
\begin{equation}\label{aux:1}
    D_{U_{1,i}}\mathcal{L}_{h,i}(\delta U_{1_i}) = \langle \nabla_{U_{1,i}}\mathcal{L}_{h,i}, \delta U_{1_i} \rangle_{L^2(0;T)} = \langle \textbf{p}_i(\cdot) \widehat{\textbf{b}} + \textbf{p}_i(\cdot + \delta_m) \widetilde{\textbf{b}} \mathbbm{1}_{[0; T-\delta_m]}\,,\; \delta U_{1_i} \rangle_{L^2(0;T)}.
\end{equation}
Indeed,
\begin{equation}
    \begin{split}
        D_{U_{1,i}}\mathcal{L}_{h,i}(\delta U_{1,i}) &= \lim_{\epsilon \rightarrow 0}\dfrac{\mathcal{L}_{h,i}(U_{1,i} + \epsilon\,\delta U_{1,i}) - \mathcal{L}_{h,i}(U_{1,i})}{\epsilon} =\\
        &=\lim_{\epsilon \rightarrow 0}\dfrac{1}{\epsilon} \left ( \int_0^{T} X_i(t)^2 \, dt - \int_0^T X_i^2(t) \, dt + \int_0^T \textbf{p}_i(t)^T \, \dot{\textbf{x}_i}(t) \, dt -  \int_0^T \textbf{p}_i(t)^T \, \dot{\textbf{x}_i}(t) \right. \, dt +\\
        &+ \int_0^T \textbf{p}_i^T(t) \, \textbf{f}(\textbf{x}, t)\, dt - \int_0^T \textbf{p}_i^T(t) \, \textbf{f}(\textbf{x}, t)\, dt + \\
        & + \int_0^T \textbf{p}_i(t)^T \widehat{\textbf{b}}\, U_{1,i}(t) + \epsilon \, \textbf{p}_i(t)^T \widehat{\textbf{b}}\,  \delta U_{1,i}(t) - \textbf{p}_i(t)^T \widehat{\textbf{b}}\, U_{1,i}(t) \, dt  + \\
        & + \left. \int_0^T \textbf{p}_i(t)^T \widetilde{\textbf{b}}\, U_{1,i}(t - \delta_m) + \epsilon \,  \textbf{p}_i(t)^T \widetilde{\textbf{b}}\,  \delta U_{1,i}(t - \delta_m) - \textbf{p}_i(t)^T \widetilde{\textbf{b}}\, U_{1,i}(t - \delta_m) \, dt  \right ) = \\
        & = \int_0^T \textbf{p}_i(t)^T \widehat{\textbf{b}}\,  \delta U_{1,i}(t) \, dt + \int_0^T \textbf{p}_i(t)^T \widetilde{\textbf{b}}\,  \delta U_{1,i}(t - \delta_m) \, dt = \\
        & = \int_0^T \textbf{p}_i(t)^T \widehat{\textbf{b}}\,  \delta U_{1,i}(t) \, dt + \int_{\delta_m}^T \textbf{p}_i(t)^T \widetilde{\textbf{b}}\,  \delta U_{1,i}(t - \delta_m) \, dt= \\
        & = \int_0^T \textbf{p}_i(t)^T \widehat{\textbf{b}}\,  \delta U_{1,i}(t) \, dt + \int_0^{T - \delta_m} \textbf{p}_i(t + \delta_m)^T \widetilde{\textbf{b}}\,  \delta U_{1,i}(t) \, dt= \\
        &= \langle \textbf{p}_i(\cdot) \widehat{\textbf{b}} + \textbf{p}_i(\cdot + \delta_m) \widetilde{\textbf{b}} \mathbbm{1}_{[0; T-\delta_m]}\,,\; \delta U_{1_i} \rangle_{L^2(0;T)} = \langle \nabla_{U_{1,i}}\mathcal{L}_{h,i}, \delta U_{1_i} \rangle_{L^2(0;T)}.
    \end{split}
\end{equation}
{Plugging \eqref{eq:descend} into \eqref{aux:1} ensures $ D_{U_{1,i}}\mathcal{L}_{h,i}(\delta U_{1,D_i})<0$ which implies, in view of \eqref{pb:compact}, that \eqref{eq:descend} is a descent direction for the cost functional $\mathcal{J}_X$.} 

    \item Update the control variables at the current step as in the Projected Gradient Descent method with Armijo adaptive learning step:
    {
    $$ 
    U_{1_i}^{new}=\Pi(U_{1_i}^{old} - \alpha\,\delta U_{1_i})
    $$
    where $\Pi$ stands for the projection operator on the space of admissible controls (i.e. satisfying the constraints).}
    At the discrete level, all the constraints are linear, therefore we applied the Shalev-Schwarz method \cite{shalev2006efficient} for linear projection over a $d$-dimensional simplex.
\end{enumerate}
    \item[] \textbf{Stopping criterion}. {The algorithm stops when the difference between two successive values of the cost functional is lower than a fixed tolerance $tol$.} 
\end{itemize}

\section{Results and Discussion}
\label{Results}
In this section we present and discuss the numerical results of the solution of the optimal control problem \eqref{pb:OCP} in the context of the Italian third epidemic wave of SARS-CoV-2 (first half of 2021) corresponding to the beginning of the vaccination campaign. We are aware of the main limitations of the considered model (see Section \ref{MathematicalModel}), hence we do not intend to propose a critical retrospective analysis of the actual implemented vaccination policy. {\color{black} More realistically, we aim at: (1) extracting the main differences between the optimal solutions obtained minimizing different cost functionals; (2) highlighting the structural features (e.g. in terms of age stratification) of the obtained vaccination strategies. In the following, the population is split in five compartments depending on age: 
\begin{equation}
    N_A = \{ (0\div19), \,(20\div39), \,(40\div59), \,(60\div79), \,(80+) \}.
\end{equation}
We set the optimization process in the period from February 12th, 2021, to June 1st, 2021, \textit{i.e.} we consider 15 weeks. However, we carried out a calibration process of the model parameters in the six months starting on January 1st, 2021. The first month of the Italian vaccination campaign has been neglected in the optimization framework since it has been devoted to the immunisation of the sanitarian personnel. In this period the vaccination priority order has been dictated by a very specific political choice, impossible to be embodied in our age-stratified compartmental model that does not take into account for working classes. All the choices on the parameters, the results of the calibration process and the sensitivity analyses on the reproduction number can be found in the Appendices A and B. We recover the initial conditions of each compartment by running a direct simulation of \eqref{pb:compact} starting on January 1st, 2021 and ending on February 12th, 2021. Indeed, the assumption of initializing to 0 all the vaccination states is reasonable on January 1st, 2021, and it helps in limiting uncertainty on the initial conditions to the other states (see Appendix B). The implemented code is available in the dedicated GitHub repository \cite{ziarelli2022github}. Both the state and the adjoint problems have been solved employing a Runge-Kutta method of order 4, with a time step of 1 day. The gradient of the Hamiltonian function with respect to each state variable, necessary in \eqref{eq:hamiltonianSys}, has been retrieved through automatic differentiation. In all our simulations the tolerance of the optimization iterative scheme has been fixed to 5e-4.
}

Each analysis is supported by quantitative results as the evolution of the cost functional, the optimized control variables by age-stratification and the discrepancies in terms of deceased, infected and hospitalized between the simulated optimal solution and the simulated solution with the starting policy. The discrepancies (with sign) are measured as follows: 
\begin{enumerate}
    \item \textbf{Daily Variation of Infected}. $\displaystyle \Lambda_I : \mathbb{N} \rightarrow \mathbb{R}$
    \begin{equation}\label{eq:DVI}
        \Lambda_I(d) = \displaystyle \sum_{i \in N_A} [\textbf{I}_{i,IC}]_d - [\textbf{I}_{i,OC}]_d, \; \forall d \leq N_{days}, \, d \in \mathbb{N}.
    \end{equation}
    where $\textbf{I}_{i, IC}$ is the vector of daily infected belonging to the $i$-th age class, obtained from the solution of the state problem \eqref{pb:compact} where the daily amount of administered doses is set equal to the value actually administered in Italy. The vector $\textbf{I}_{i, OC}$ contains the same type of informations computed by solving \eqref{pb:compact} with the optimal vaccination campaign coming from the solution of \eqref{pb:OCP};
    \item \textbf{Daily Variation of Hospitalized}. $\displaystyle \Lambda_H : \mathbb{N} \rightarrow \mathbb{R}$
    \begin{equation} \label{eq:DVH}
        \Lambda_H(d) = \displaystyle \sum_{i \in N_A} [\textbf{H}_{i,IC}]_d - [\textbf{H}_{i,OC}]_d, \; \forall d \leq N_{days}, \, d \in \mathbb{N}.
    \end{equation}
    where $\textbf{H}_{i, IC}$ is the vector of daily hospitalized of the $i$-th age class of the solution of the state problem \eqref{pb:compact} with vaccination imposed as the vaccinations actually implemented in Italy, and $\textbf{H}_{i, OC}$ is the amount of hospitalized computed with the optimal control variables;
    \item \textbf{Daily Variation of Deceased}. $\displaystyle \Lambda_D : \mathbb{N} \rightarrow \mathbb{R}$
    \begin{equation} \label{eq:DVD}
        \Lambda_D(d) = \displaystyle \sum_{i \in N_A} [\textbf{D}_{i,IC}]_d - [\textbf{D}_{i,OC}]_d, \; \forall d \leq N_{days}, \, d \in \mathbb{N}.
    \end{equation}
    where $\textbf{D}_{i, IC}$ represents the vector of daily deaths of the $i$-th age class of the solution of the state problem \eqref{eq:hamiltonianSys} with vaccination imposed as the vaccinations actually implemented in Italy, $\textbf{D}_{i, OC}$ is the respective counterpart computed with the optimal control variables;
\end{enumerate}

{\color{black}
The results are organized as follows. In Sections \ref{MinDec}-\ref{MinHosp} we present and discuss the optimal solutions obtained by minimizing deceased, infected and hospitalized, respectively. In Section \ref{optDiffIGITA} we explore the impact  on the optimization of the vaccination policy of different choices for the initial guess, while in Section \ref{sec:constraints}  we explore the impact of different constraints on the available amount of doses to be administered during the period of interest.}
\subsection{Minimization of deceased}
\label{MinDec}
In this section we optimize the vaccination campaign (\textit{i.e.} the administration doses by ages) with respect to the total amount of deceased caused by SARS-CoV-2 infections, cf. \eqref{JD}. The initial guess is chosen equal to the vaccination campaign actually implemented at the national level. In particular, from the DPC data \cite{dpcdatavax} we reconstruct the vaccination policies of first doses applied in Italy split by age-classes and we average the respective values on a weekly basis as shown in Figure A.3.
From Figure \ref{costDec} we notice that during the optimization process the cost functional has a descending behaviour corresponding to a reduction in terms of deceased of nearly 498 individuals at the final time. The algorithm stops after 951 iterations reaching the desired tolerance. 
\begin{figure}[H]
    \centering
    \includegraphics[scale=0.4]{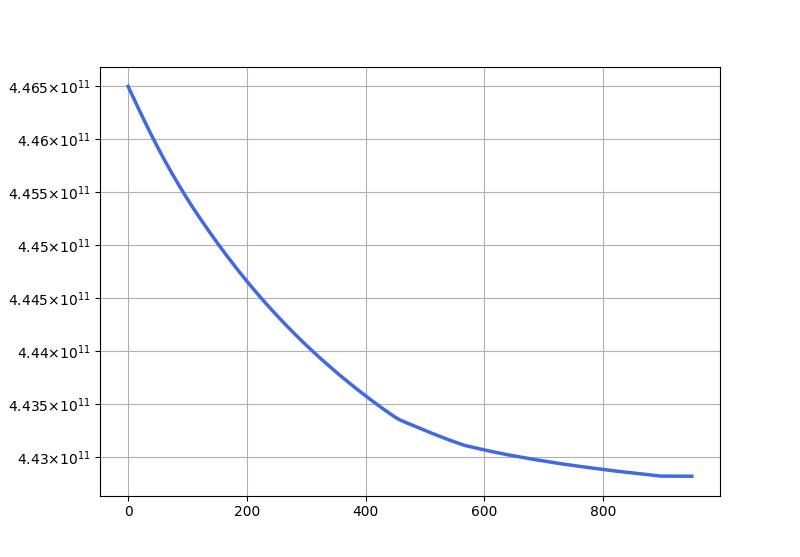}
    \caption{Minimization of deceased. Evolution of the cost functional during the optimization process. X-axis reports the iterations of the scheme.}
    \label{costDec}
\end{figure}
Figure \ref{costDecUS} compares the amount of first doses administered 
in the optimal solution with the ones administered in the initial policy, while Figure \ref{costDecPerc} shows the percentage repartition of doses across ages in the initial policy  and in the optimal solution. In Figure \ref{costDecUS} we report the weekly total initial doses assigned to each age-stratification retrieved by data available from Dipartimento di Protezione Civile Italiana (dashed line) and obtained from the solution of the optimal problem (solid line). As one can notice from both Figures, the optimal strategy to minimize deceased suggests to increase the amount of doses to the over 80s while reducing the amount of doses to the other age classes. In particular, the younger is the age class, the higher is the relative amount of the reduction of doses. These results are not completely unexpected since they point out how the strategy to minimize deceased is to administer vaccinations to the elderly as much as possible (the highest increment in a week is of more than 150 thousands doses). Indeed, the over 80s represent the most fragile in terms of probability to contract the infection in a severe and, consequently, fatal way (see the values of the $IFR_i$ in Table A.1). Moreover, Figure \ref{costDecPerc} compares the percentage repartition across age-classes of the vaccinations of the initial with the ones of the optimal policy. Notice that the optimal policy addresses to the over 60s from the 69\% to the 88\% of doses in each week, whilst the same quantity ranges from 43\% to 81\% with the initial policy.
\begin{figure}[H]
    \centering
    \includegraphics[trim={5.5cm 0 5.5cm 0},width=\textwidth]{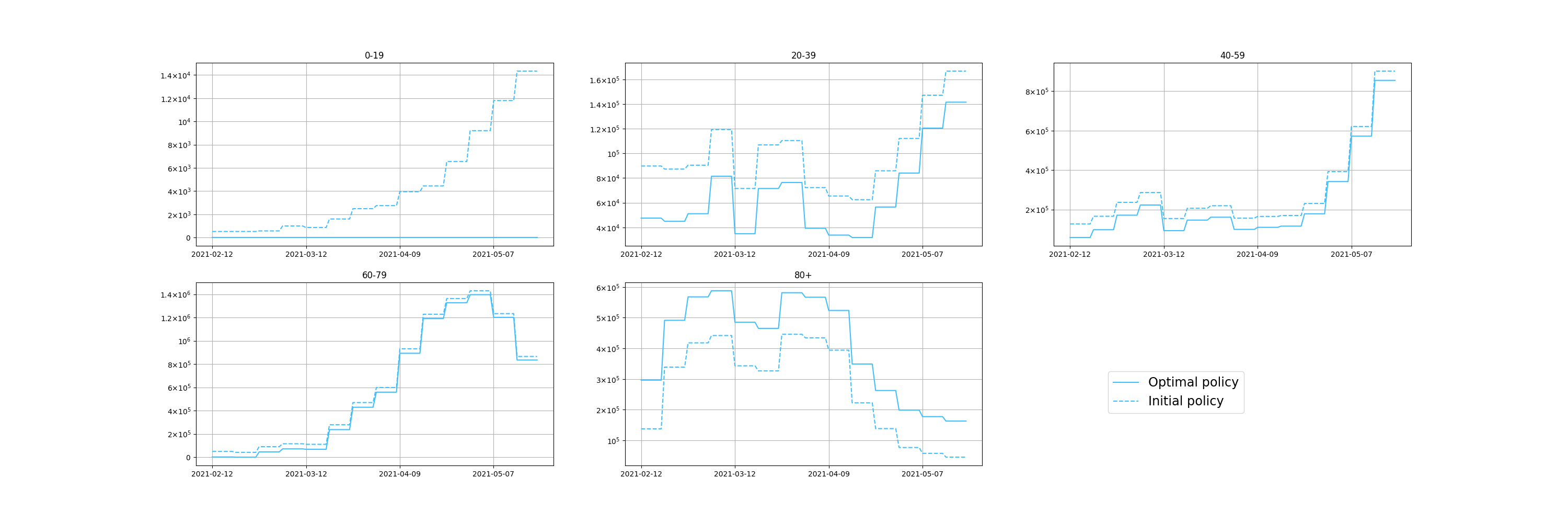}
    \caption{Minimization of deceased. Weekly amount of doses delivered for each age-stratification.}
    \label{costDecUS}
\end{figure}

\begin{figure}[H]
    \centering
    \includegraphics[scale=0.6]{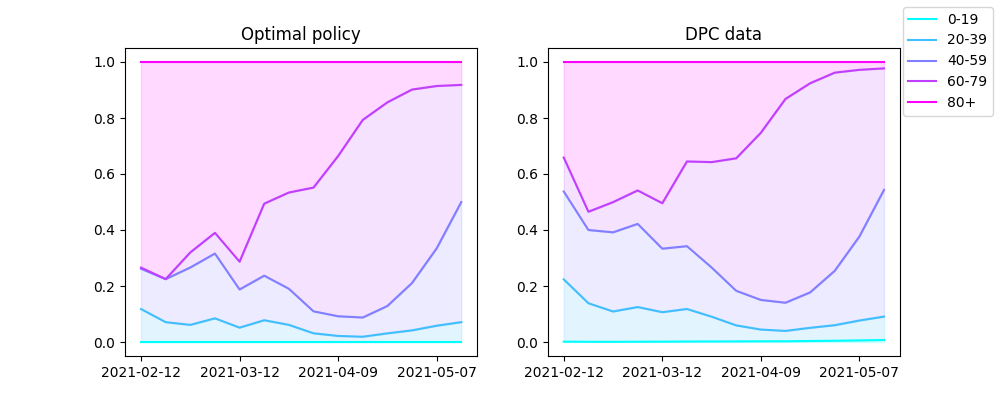}
    \caption{Minimization of deceased. Percentage repartition of doses across ages in the initial policy (left) and in the optimal solution (right).}
    \label{costDecPerc}
\end{figure}
However, as one can expect, this is not the best solution in terms of other quantities of interest such as infected and hospitalized. Indeed, from Figures \ref{costDecDec}-\ref{costDecHosp} it turns out that that the optimal solution obtained from the minimization of deceased produces an increase in terms of infected and hospitalized especially at the final time (red curves assess a negative Variation of Infected and Hospitalized, while green curve assess a positive Variation of Deceased, cf. \eqref{eq:DVI}-\eqref{eq:DVD}). 
\begin{figure}[H]
\minipage{0.32\textwidth}
  \includegraphics[trim={0 0 0 1cm},clip, width=\linewidth]{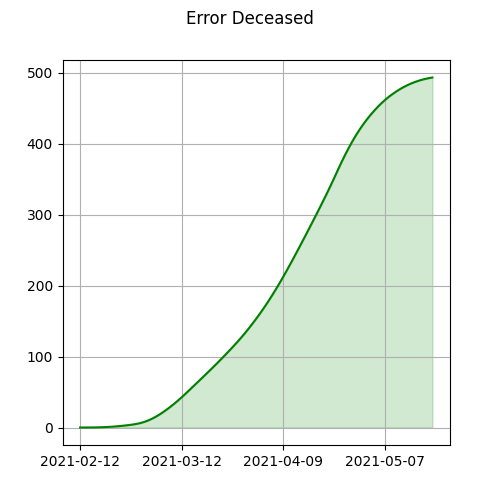}
  \caption{Absolute value of the Variation of Deceased measure ($\Lambda_D$). }\label{costDecDec}
\endminipage\hfill
\minipage{0.32\textwidth}
    \includegraphics[trim={0 0 0 1cm},clip,width=\linewidth]{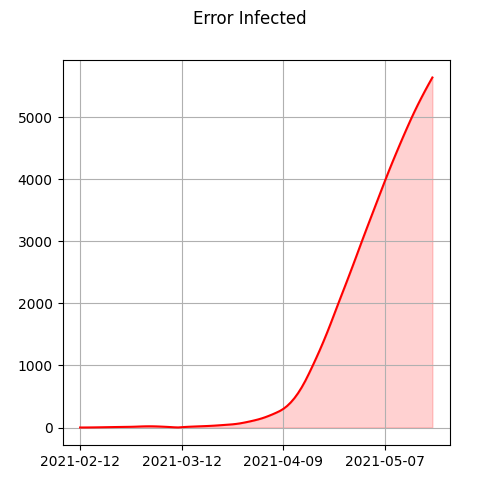}
  \caption{Absolute value of the Variation of Infected measure ($\Lambda_I$).}\label{costDecInf}
\endminipage\hfill
\minipage{0.32\textwidth}%
  \includegraphics[trim={0 0 0 1cm},clip,width=\linewidth]{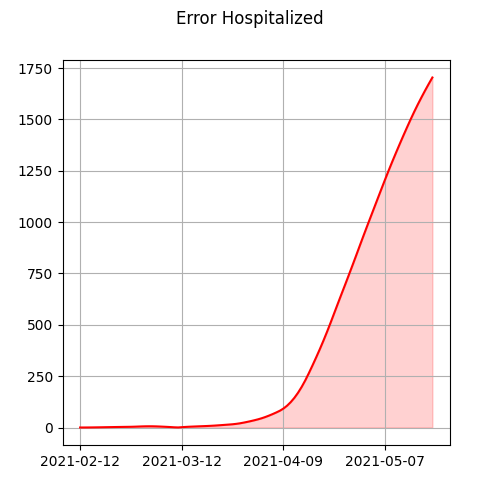}
  \caption{Absolute value of the Variation of Hospitalized measure ($\Lambda_H$).}\label{costDecHosp}
\endminipage

\end{figure}

\subsection{Minimization of infected}
\label{MinInf}
In this section we present the results of the optimal vaccination policy minimizing infected individuals (cf. \eqref{JI})during the whole time frame. In Figure \ref{costInf} we report the history of the cost functional in terms of the number of iterations: we observe a reduction of the number of infected individuals from 58.5 thousands (initial value) to 55.5 thousands (final value).
The optimization process stops after 405 iterations.
\begin{figure}[H]
    \centering
    \includegraphics[scale=0.4]{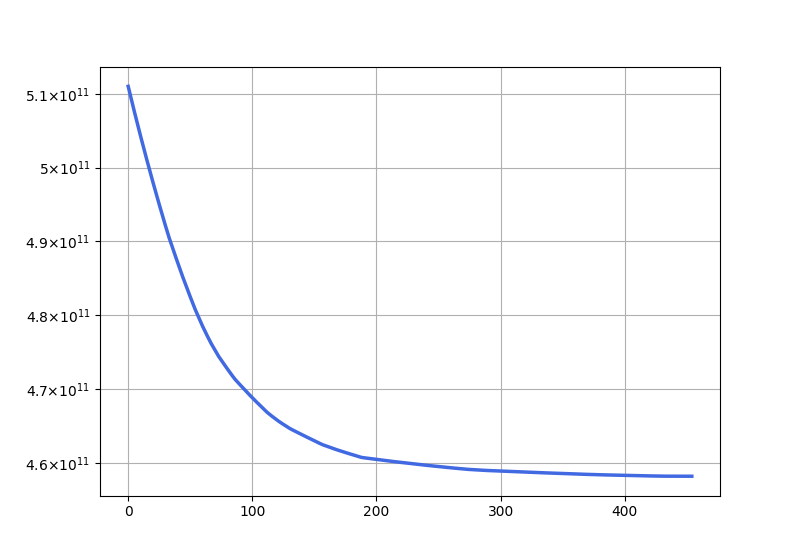}
    \caption{Minimization of infected. Evolution of the cost functional during the optimization process. X-axis reports the iterations of the scheme.}
    \label{costInf}
\end{figure}
From Figures \ref{costInfUS} and \ref{costInfPerc} we can extract the main features, in terms of age-class dose repartition, of the optimal vaccination campaign for minimizing infected. The reduction implied by the cost functional is guided by an unattended decrease in the amount of doses to the youngest and to the oldest age-classes against an increase of the amount of doses in the (20$\div$59) age-class. This result can be interpreted in light of the contact matrix weighted by the age-dependant susceptibility to the virus (Figure \ref{costInfMat}). 
This matrix has been computed multiplying each column of the POLYMOD contact matrix by the age-dependant susceptibilities $r_i$. In this way each row takes into account the absolute amount of high-infection-risk contacts that one individual belonging to a specific age class has with individuals belonging to the others. Then, its values has been normalized through the maximum norm. In particular, the two rows corresponding to the age classes (20$\div$39) and (40$\div$59) achieve the highest values, meaning that individuals ageing (20$\div$59) are the ones having the highest probability of transmitting the virus heterogeneously across ages.  Figure \ref{costInfPerc} highlights how almost the totality of the delivered doses has to be administered at these age classes to the detriment of (0$\div$19), (60$\div$79) and (80+) classes. One may notice that the entry $((0\div19), (0\div19))$ in the weighted contact matrix is also relatively high with respect to the others, even though vaccinations to this age-stratification are not promoted by the optimal controlled solution. 
\begin{figure}[H]
    \centering
    \includegraphics[trim={5.5cm 0 5.5cm 0},width=\textwidth]{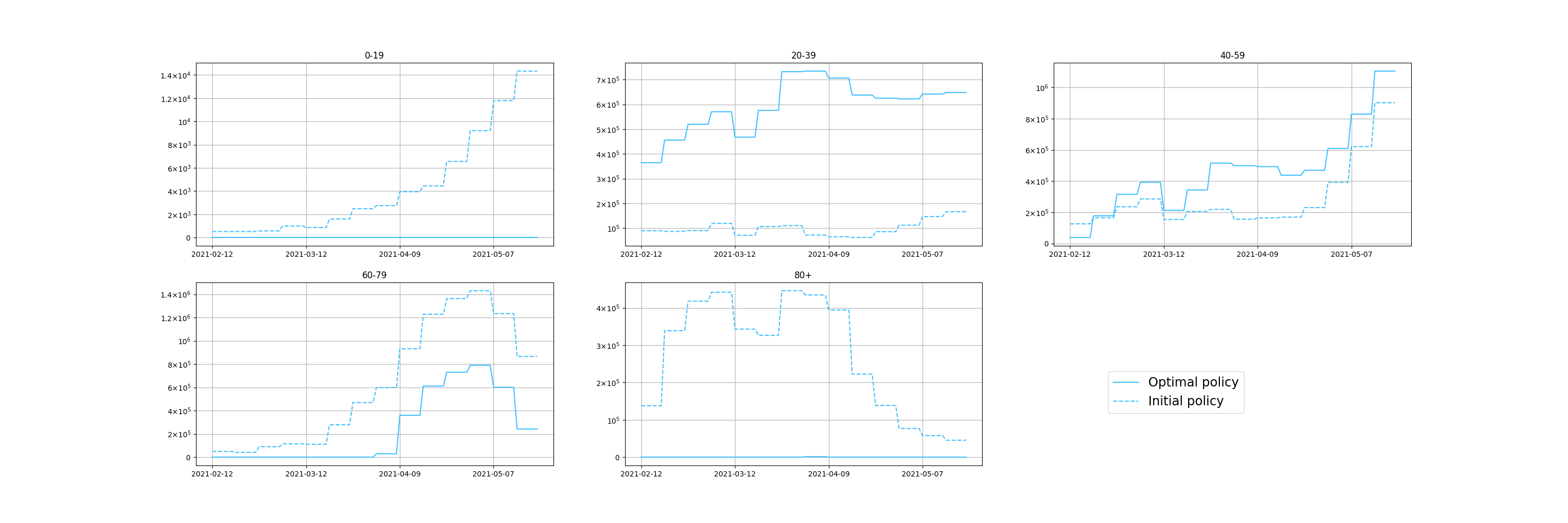}
    \caption{Minimization of infected. Weekly amount of doses delivered for each age-stratification.}
    \label{costInfUS}
\end{figure}

\begin{figure}[H]
    \centering
    \includegraphics[scale=0.6]{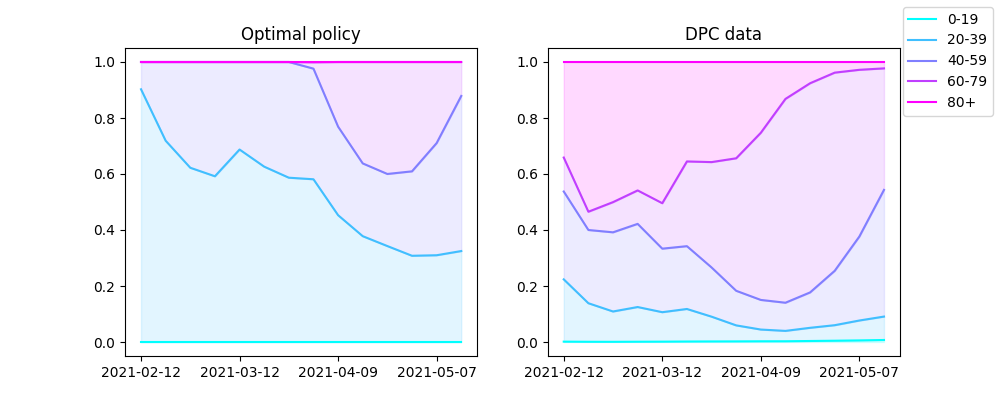}
    \caption{Minimization of infected. Percentage repartition of doses across ages in the initial policy (left) and in the optimal solution (right).}
    \label{costInfPerc}
\end{figure}

\begin{figure}[H]
    \centering
    \includegraphics[scale=0.3]{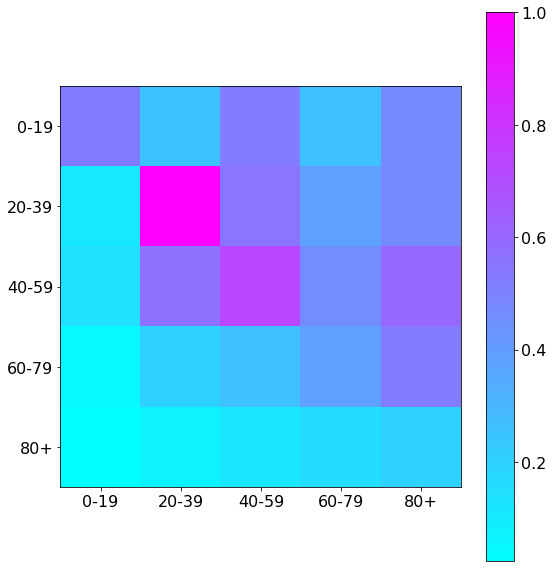}
    \caption{Normalized contacts weighted by age-dependant susceptibilities $r_i$.}
    \label{costInfMat}
\end{figure}

Finally, let us notice that, similarly to the previous section, also in this case the optimal solution obtained by minimizing infected does not necessarily imply the minimization of deceased. Indeed, without protecting from the illness those who are more prone to severe outcomes, this optimal solution prescribes an increase in terms of deaths of more than 1750 units. On the other hand, hospitalized, which are strictly linked to the infected individuals, are correspondingly reduced even though their reduction was not directly contemplated directly in the definition of the cost functional (see Figures \ref{costInfDec}-\ref{costInfHosp}).

\begin{figure}[H]
\minipage{0.32\textwidth}
  \includegraphics[trim={0 0 0 1cm},clip,width=\linewidth]{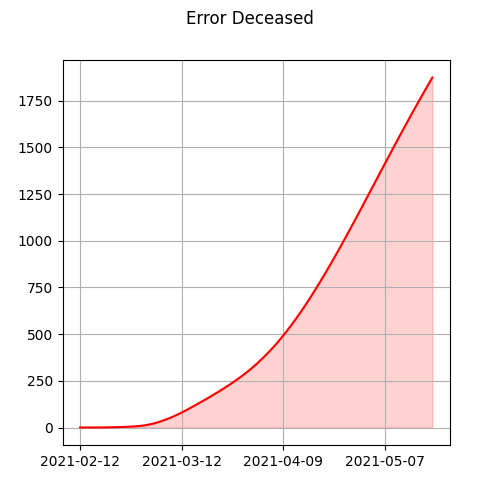}
  \caption{Absolute value of the Variation of Deceased measure ($\Lambda_D$).}\label{costInfDec}
\endminipage\hfill
\minipage{0.32\textwidth}
    \includegraphics[trim={0 0 0 1cm},clip,width=\linewidth]{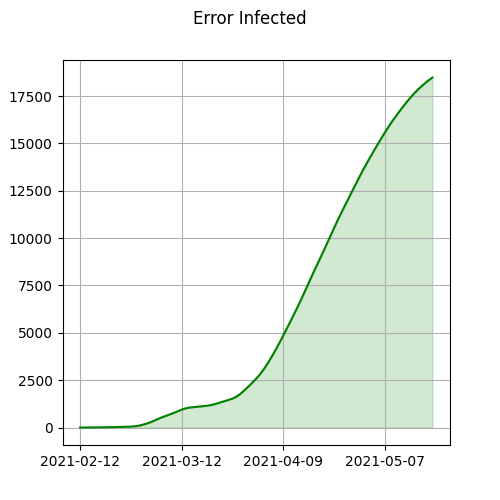}
  \caption{Absolute value of the Variation of Infected measure ($\Lambda_I$).}\label{costInfInf}
\endminipage\hfill
\minipage{0.32\textwidth}%
  \includegraphics[trim={0 0 0 1cm},clip,width=\linewidth]{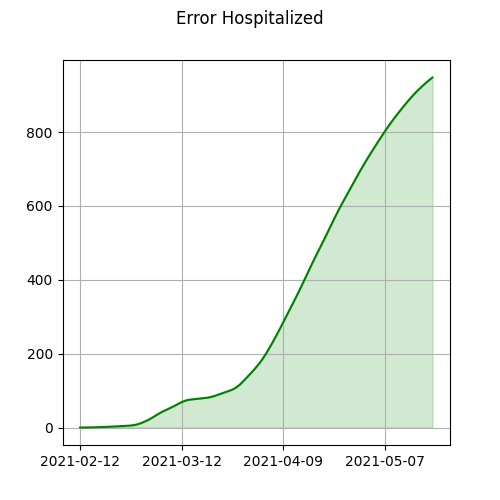}
  \caption{Absolute value of the Variation of Hospitalized measure ($\Lambda_H$).}\label{costInfHosp}
\endminipage
\end{figure}

\subsection{Minimization of hospitalized}
\label{MinHosp}
In this section we collect the results of the optimal vaccination policy obtained by minimizing hospitalized individuals (cf. \eqref{JH}) during the period of interest. After 1947 iterations the algorithm stops statisfying the prescribed stopping criterion. During the whole time horizon, the amount of hospitalized is reduced by 70.3 thousands individuals.
This is an important remark, since during the first two waves of COVID19 epidemic one of the greatest issues was dealing with the limited amount of ICUs (Intensive Care Units) in hospitals and the extremely high level of beds occupancy due to people with respiratory symptoms linked to SARS-CoV-2. Reducing hospitalized individuals is one of the main objectives to be taken into account assuming a social perspective.
\begin{figure}[H]
    \centering
    \includegraphics[scale=0.4]{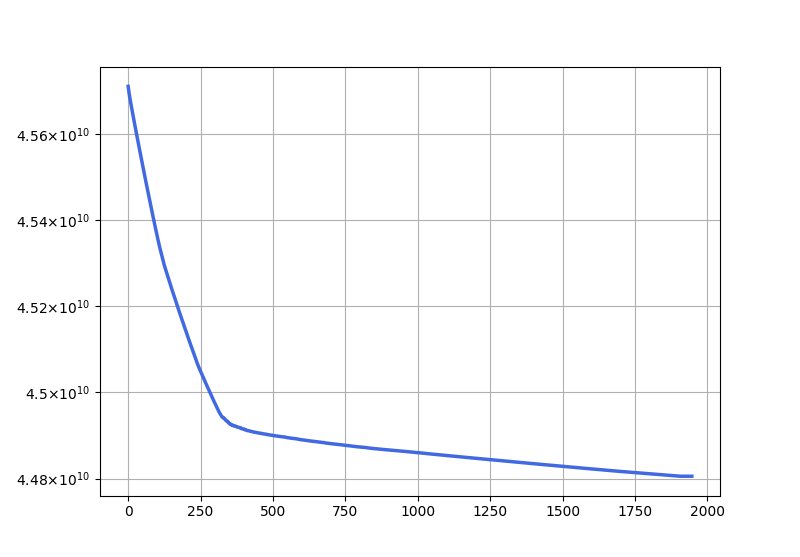}
    \caption{Minimization of hospitalized. Evolution of the cost functional during the optimization process. X-axis reports the iterations of the scheme.}
    \label{costHosp}
\end{figure}
Figure \ref{costHospUS} shows which is the optimal strategy to reduce the amount of hospitalized individuals. It is suggested to increase the amount of doses to administer to the age stratification (20$\div$59), without neglecting the administrations to the over 60s, but reducing the administrations to the (0$\div$19). Indeed, the age-classes to which an higher absolute amount of doses has to be provided are the ones more prone to host the virus in a more severe way. This policy is confirmed by Figure \ref{costHospPerc} representing the age repartition of doses.
\begin{figure}[H]
    \centering
    \includegraphics[trim={5.5cm 0 5.5cm 0},width=\textwidth]{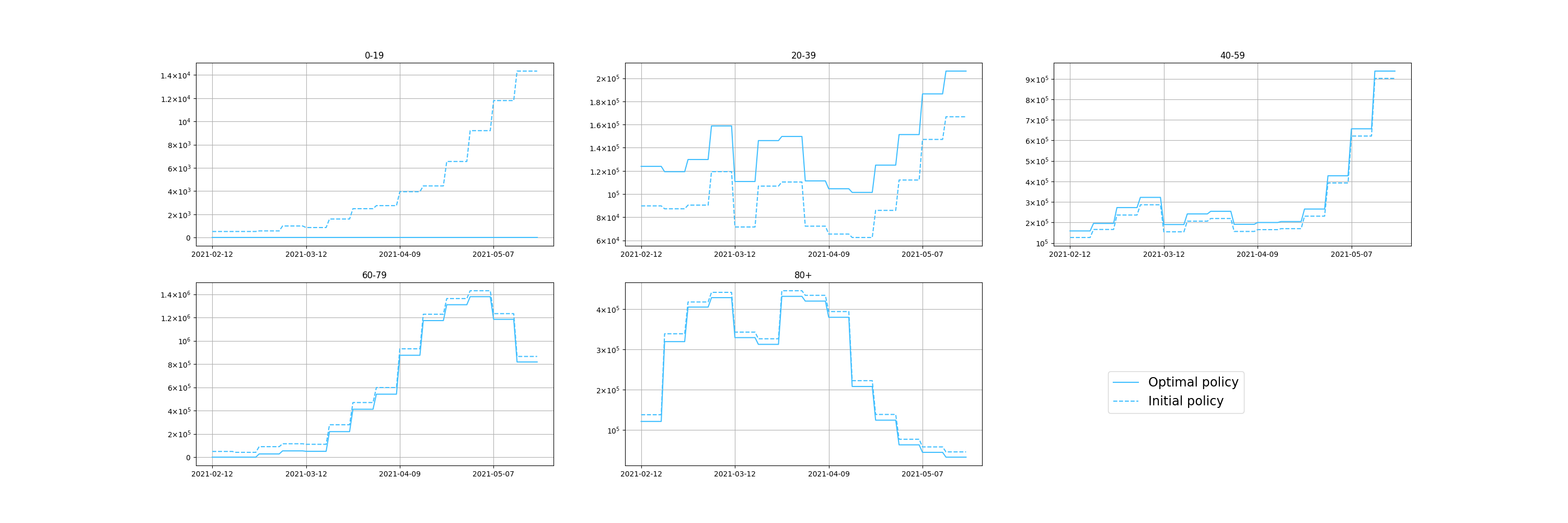}
    \caption{Minimization of hospitalized. Weekly amount of doses delivered for each age-stratification.}
    \label{costHospUS}
\end{figure}

\begin{figure}[H]
    \centering
    \includegraphics[scale=0.6]{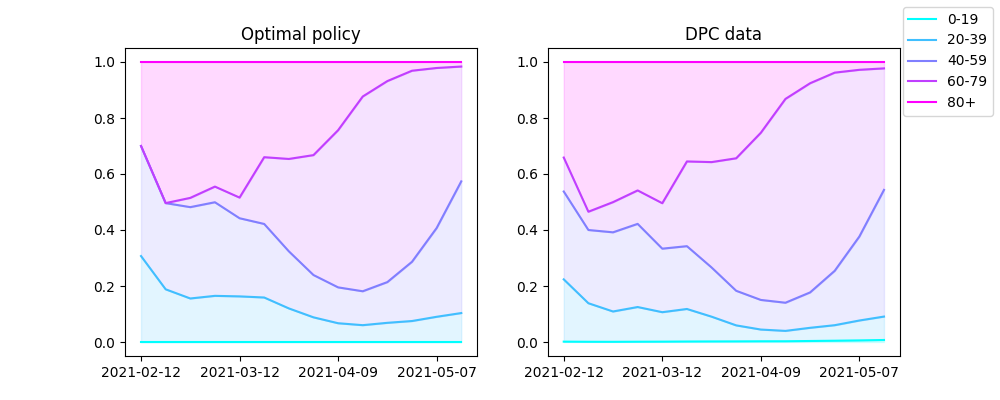}
    \caption{Minimization of hospitalized. Percentage repartition of doses across ages in the initial policy (left) and in the optimal solution (right).}
    \label{costHospPerc}
\end{figure}
Lastly, we note that the optimal solution that reduces the most the amount of hospitalized individuals during the whole time frame is actually a solution increasing the amount of deceased. Indeed, in the optimal policy the amount of doses delivered to the most fragile people (the elderly) is slightly reduced with respect to the ones administered in the initial guess. Moreover, even though we have reduced the total amount of hospitalized, this solution does not improve the amount of infected with respect to the solution reducing the infectious in the interval of interest (Figure \ref{costInfUS}). Although the total amount of infected in this solution is increased with respect to the one minimizing infectious, the value of hospitalized is lower due to the different repartition of infected across ages (notice that in \eqref{JH} the infected belonging to different age-classes are weighted differently by age). 
\begin{figure}[H]
\minipage{0.32\textwidth}
  \includegraphics[trim={0 0 0 1cm},clip,width=\linewidth]{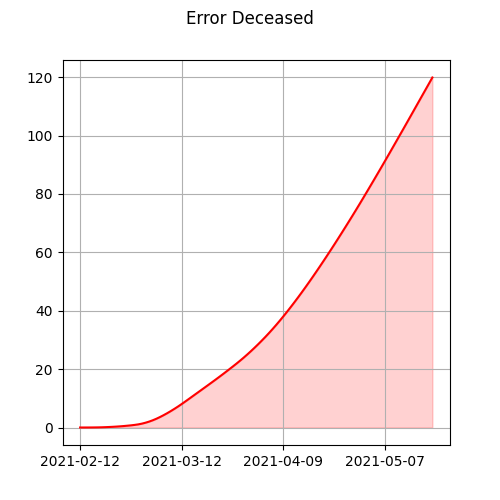}
  \caption{Absolute value of the Variation of Deceased measure ($\Lambda_D$).}\label{costHospDec}
\endminipage\hfill
\minipage{0.32\textwidth}
    \includegraphics[trim={0 0 0 1cm},clip,width=\linewidth]{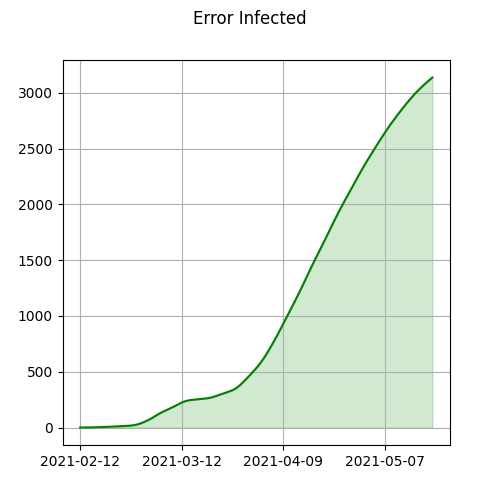}
  \caption{Absolute value of the Variation of Infected measure ($\Lambda_I$).}\label{costHospInf}
\endminipage\hfill
\minipage{0.32\textwidth}%
  \includegraphics[trim={0 0 0 1cm},clip,width=\linewidth]{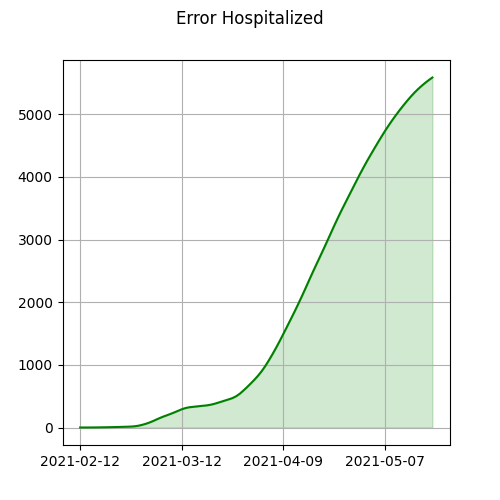}
  \caption{Absolute value of the Variation of Hospitalized measure ($\Lambda_H$).}\label{costHospHosp}
\endminipage
\end{figure}

\subsection{Optimal solution with different initial guesses}
\label{optDiffIGITA}
{ As the the optimal control problem \eqref{pb:OCP} is solved via a Projected Gradient Descent method, the solution is typically a local minimum and it is influenced by the specific choice of the initial policy. In view of this remark, in the present section we investigate the dependence on the chosen initial policy of the  solution of \eqref{pb:OCP}.
Particularly, we compare the results obtained in the previous sections (where the initial guess was set equal to the actual implemented national vaccination campaign) with the ones obtained starting from a different initial guess. The new initial guess that has been considered as starting policy is the homogeneous allocation of doses proportionally to the populosity of each age stratification. Specifically, let $N_{week}(d)$ be the total amount of doses that have been distributed during the week at which day $d$ belongs. Then, the initial guess of the amount of doses to be distributed each day of the week correspondent to $d$ can be computed as
\begin{equation}\label{new_ini}
    U_{1,i}^{(0)}(d) = \displaystyle \dfrac{N_{week}(d)}{7 \,\mathrm{(days)}} \dfrac{N_i}{N} = \dfrac{N_{week}(d)}{7} \displaystyle \dfrac{N_i}{ \displaystyle \sum_{k \in N_A} N_k},\; \forall d \in \mathrm{I},\,i \in N_A.
\end{equation}
Figure \ref{InfIC} shows the comparison between infected individuals generated by the optimal policy starting from the DPC-oriented policy (DPC-IG) and the one from the homogeneous doses initial guess (Homogeneous-IG). Instead, Figure \ref{DecIC} displays the comparison of deceased between the two optimal solutions minimizing deceased, and lastly, Figure \ref{HosIC} the trends of hospitalized when optimizing hospitalizations.
The states have been simulated in a supplemental interval of 40 days, considering the calibrated values of the transmission rates and assuming that no first-dose administrations are allocated to each age class after the June 1st, 2021, this in order to have fair comparison among the solutions. The optimal policies retrieved in each case, together with their specific initial guesses (dashed lines), are presented in Figures \ref{costInfUSIC}--\ref{costHosUSIC}. For the solution optimizing the number of infected individuals, one can notice from Figure \ref{InfIC} that the two solutions overlap almost completely during the whole time interval. However, the respective policies in Figure \ref{costInfUSIC} do not coincide: they both tend to allocate more than the 80\% of the weekly available doses to individuals of age comprised between 20 and 59, neglecting completely the doses to administer to the older and the younger age-classes. However, the two solutions are actually different starting from April to the end of the simulation considering the doses to allocate to the (60$\div$79) years old. Hence, the solution starting with the DPC-based initial guess stops at a local minimum with non-zero administrations to this age-class; this is probably due to the ifluence of the chosen initial guess. 
On the other hand, the solutions minimizing deceased agree in administering more doses to the (80+) category, as we previously commented in Section \ref{MinDec} (see Figure \ref{costDecUSIC}), although the curves of the optimal policies keep the same trend of the initial guesses which are different from the two simulations. The (Homogeneous-IG) optimal solution starts from a value of total deceased of 106244 and reduces deceased to the optimal value of 102657, whilst the other starts with an IG value of 103000 up to 102213 deceased. Therefore, the trend suggested by the initial policy proportional to the populosity of the compartment generate an higher reduction of deceased with respect to the initial guess, although the DPC-oriented one reaches an improved optimal value in terms of absolute amount of deceased. 
Finally, the solutions for the minimization of hospitalized (Figures \ref{HosIC} and \ref{costHosUSIC}) concur in increasing the amount of doses to the (20$\div$59) accordingly to the specific trend of the initial guess, and decreasing allocations corresponding to the other age-classes. 
However, the optimal policy starting from homogeneous initial guess allocates more doses to the over 80 with respect to the initial guess. This is not unexpected since individuals belonging to this age-class are the ones more prone to contract the disease in its most severe form, and so they require hospitalization more often than the other age-states. Hence, in the hospitalized-reduction case the social interactions responsible for the spread of infections and age-dependent frailties are two fundamental components that do not prevail with each other and have to be both taken into account during the administration process. 

Summarizing, starting from different initial guesses, the optimal solutions return the same minimisation values when minimising infected and hospitalised, while there is a significant difference of about 444 deaths in the case of minimising deceased. On the other hand, in all cases the optimal policies are different confirming the local nature of the optimisation method used (PGD) and the impact of the initial guess. However, the increasing and decreasing trends by age of the vaccine doses in the optimal solutions are preserved in the case of minimisation of the infected and the deceased, while for the hospitalised the doses to the elderly are different between the two solutions (increasing doses for the homogeneous solution with respect to the initial policy and decreasing the same amount for the solution starting with the DPC-oriented policy).
\begin{figure}[H]
\minipage{0.32\textwidth}
  \includegraphics[width=\linewidth]{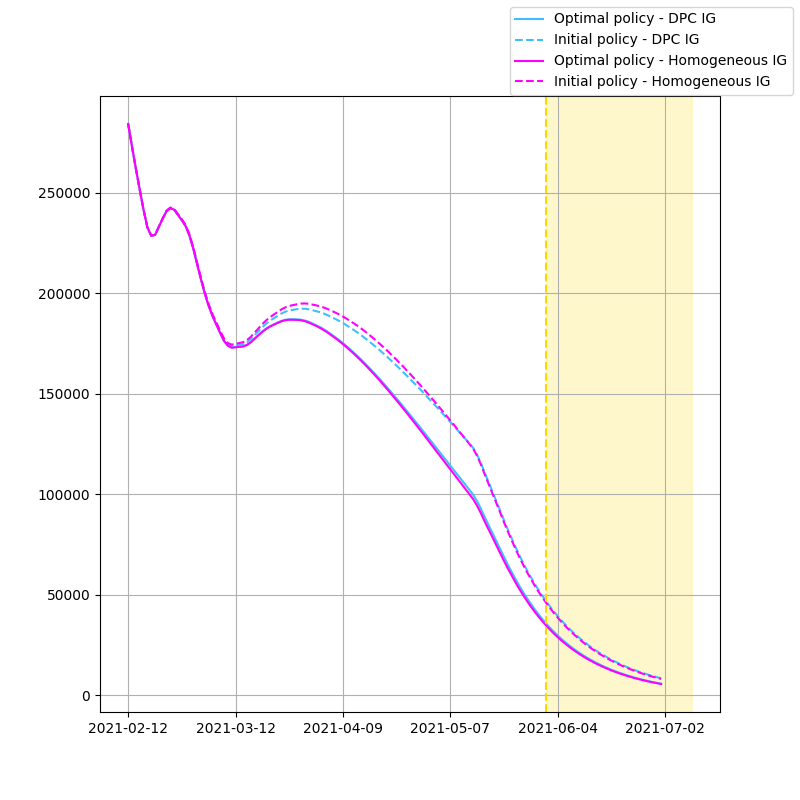}
  \caption{Evolution of infected of the two solutions obtained minimizing infected starting from the two distinct initial guesses.}\label{InfIC}
\endminipage\hfill
\minipage{0.32\textwidth}
    \includegraphics[width=\linewidth]{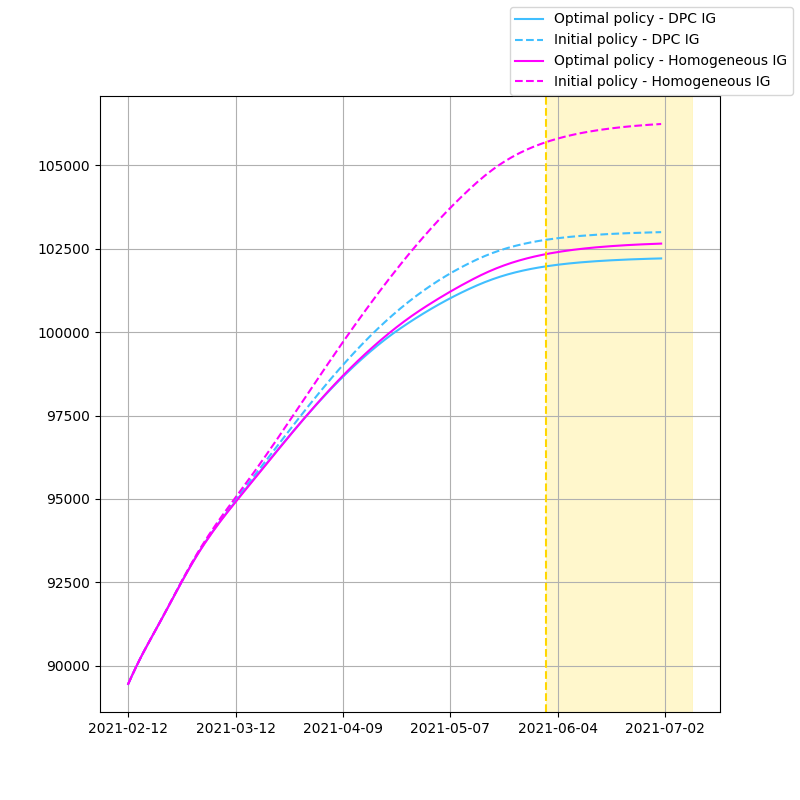}
  \caption{Evolution of deceased of the two solutions obtained minimizing deceased starting from the two distinct initial guesses.}\label{DecIC}
\endminipage\hfill
\minipage{0.32\textwidth}%
  \includegraphics[width=\linewidth]{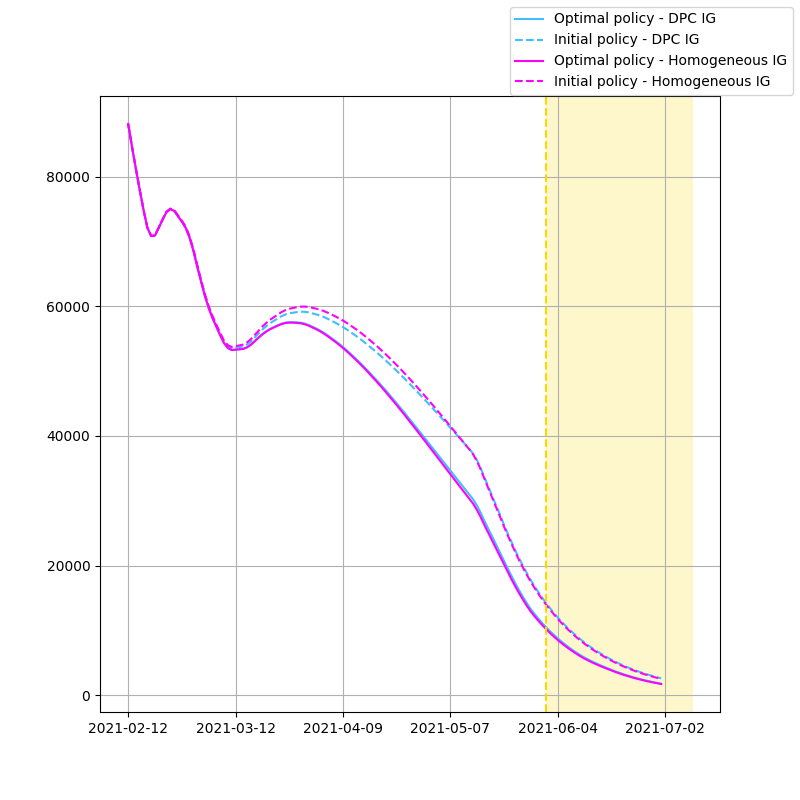}
  \caption{Evolution of hospitalized of the two solutions obtained minimizing hospitalized starting from the two distinct initial guesses.}\label{HosIC}
\endminipage
\end{figure}

\begin{figure}[H]
    \centering
    \includegraphics[trim={5.5cm 0 5.5cm 0},width=\textwidth]{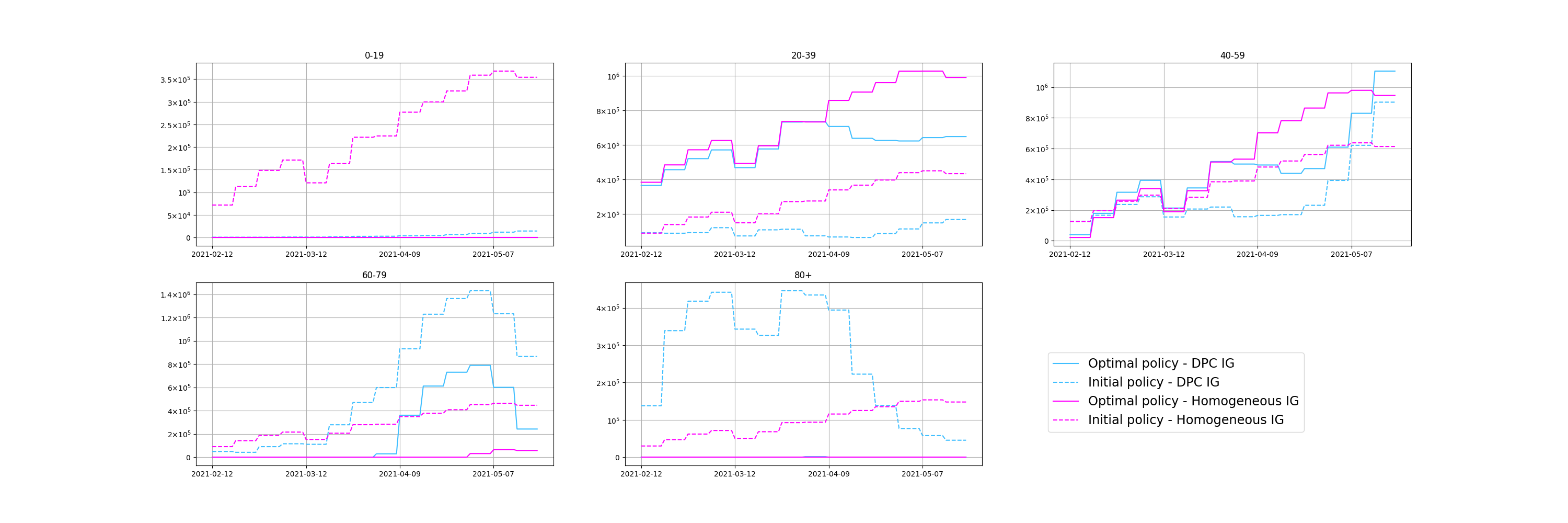}
    \caption{Weekly amount of doses delivered for each age-stratification in the solutions minimizing infected starting from DPC-IG and Homogeneous-IG.}
    \label{costInfUSIC}
\end{figure}

\begin{figure}[H]
    \centering
    \includegraphics[trim={5.5cm 0 5.5cm 0},width=\textwidth]{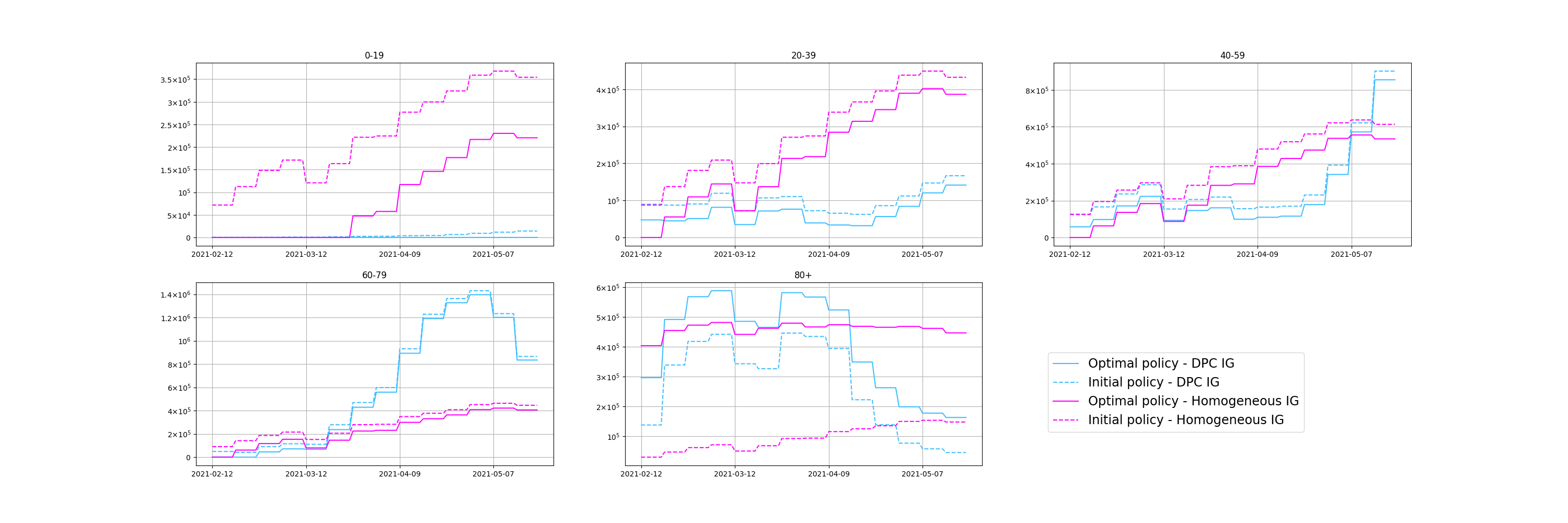}
    \caption{Weekly amount of doses delivered for each age-stratification in the solutions minimizing deceased starting from DPC-IG and Homogeneous-IG.}
    \label{costDecUSIC}
\end{figure}
\begin{figure}[H]
    \centering
    \includegraphics[trim={5.5cm 0 5.5cm 0},width=\textwidth]{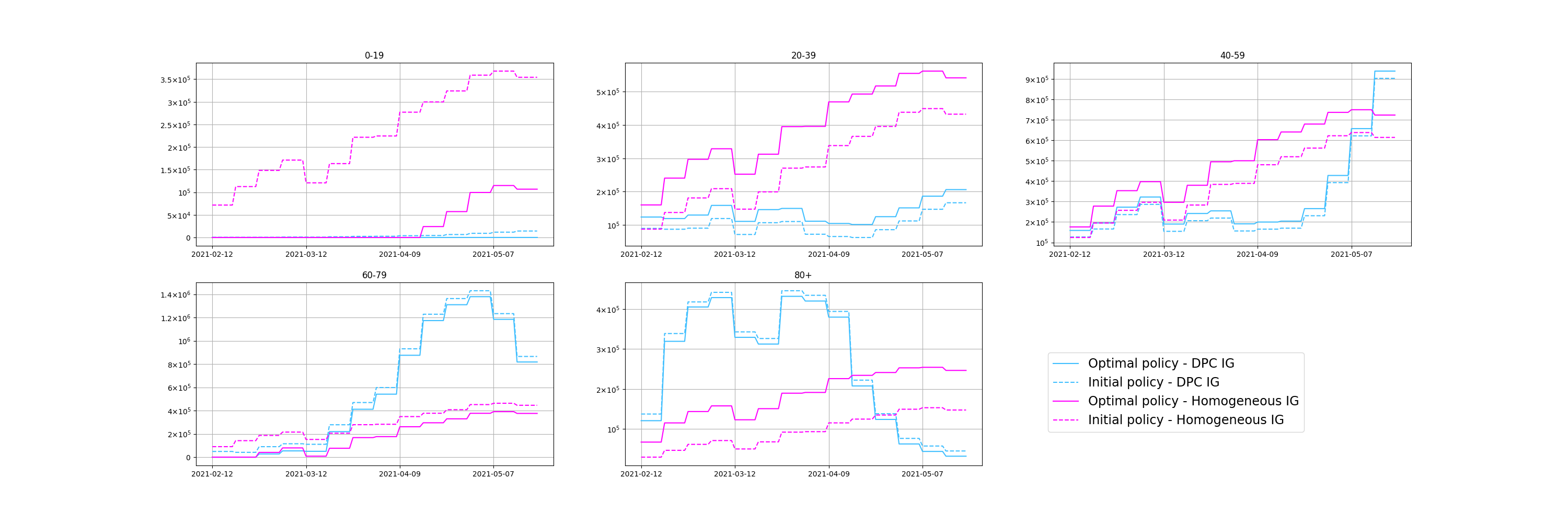}
    \caption{Weekly amount of doses delivered for each age-stratification in the solutions minimizing hospitalized starting from DPC-IG and Homogeneous-IG.}
    \label{costHosUSIC}
\end{figure}




\subsection{Optimal solution subject to different initial reproduction numbers}\label{sec:constraints}
In this section we consider different scenarios where the constraints on the specific amount of doses administered in Italy in the first half of 2021 has been relaxed. In addition, we fixed the transmission rate during the whole time frame, the value having  been selected so to deal with paradigmatic values of the initial reproduction numbers (i.e. smaller, slightly larger, larger than $1$).
Specifically, the right-hand side in constraint \eqref{budgConstraint} is a constant function of time, imposed at the constant value of $N_{week} =2.1$ millions doses per week, which is plausible approximation of the average amount of doses administered in Italy when the sanitarian capacity has reached its maximum. Moreover, following the scenario-analyses proposed in \cite{giordano2021modeling}, we fix the transmission rate during the whole time frame at three different values in order to achieve three different initial reproduction numbers, \textit{i.e.} $\mathcal{R}_0 = \{0.72, 1.01, 1.30 \}$.  In this way we consider optimal policies in presence of a minor epidemic (Case 1, $\mathcal{R}_0 \simeq 0.7$), at the bifurcation value (Case 2, $\mathcal{R}_0, \simeq 1$) and in presence of a major outbreak (Case 3, $\mathcal{R}_0, \simeq 1.3$). The values of the transmission rate corresponding to the three reproduction numbers have been computed from the definition of the initial reproduction number for the $SIRDVW$ model, that is 
\begin{equation}
    \mathcal{R}_0 := \dfrac{\beta}{\gamma}.
\end{equation}
Each scenario has been simulated starting from three different initial guesses, representative of possible different approaches to the vaccination campaign:
\begin{itemize}
    \item \textbf{IG1}: each age class receives the proportion of total doses correspondent to the proportion of population of the same age-class on the total, \textit{i.e.}
    \begin{equation}
    U_{1,i}^{(0)}(d) = \displaystyle \dfrac{1}{2}\,\dfrac{N_{week}}{7\, \mathrm{(days)}} \dfrac{N_i}{N_t} = \dfrac{N_{week}}{14} \displaystyle \dfrac{N_i}{ \displaystyle \sum_{k \in N_A} N_k},\; \forall d \in I,\, i \in N_A.
    \label{eq:secResDoses1}
\end{equation}
    The total amount of doses are halved between first and second administrations. This explains the coefficient $\frac{1}{2}$ appearing in \eqref{eq:secResDoses1} (and in \eqref{eq:secResDoses2});
    \item \textbf{IG2}: the amount of doses assigned to each age-class are proportional to the correspondent Infectious Fatality Rate, \textit{i.e.}
    \begin{equation}
    U_{1,i}^{(0)}(d) = \displaystyle \dfrac{1}{2}\,\dfrac{N_{week}}{7\, \mathrm{(days)}} \dfrac{IFR_i}{IFR_t} = \dfrac{N_{week}}{14} \displaystyle \dfrac{IFR_i}{\displaystyle \sum_{k \in N_A} IFR_k},\; \forall d \in I,\,i \in N_A;
    \label{eq:secResDoses2}
\end{equation}
\item \textbf{IG3}: we deal with a total amount of first doses administrations as a square wave with phase equal to three weeks and peak equal to the total weekly amount of doses. In each age class we distribute the amount of doses proportionally to the populosity as in IG1, \textit{i.e.}
\begin{equation}
    \begin{split}
        &N_{sw,doses}(t) = \begin{cases} 
        N_{week} & t \in [21\, (2k), 21 \, (2k+1))\; \forall k \in \mathbb{N} \;, 0 \leq t \leq T_f, \\
      0 & t \in [21\, (2k +1), 21 \, (2(k+1)))\; \forall k \in \mathbb{N} \;, 0 \leq t \leq T_f; \\
   \end{cases} \\
   &U_{1,i}^{(0)}(d) = \displaystyle N_{sw,doses}(d) \,\dfrac{N_i}{N_t},\; \forall d \in I,\,i \in N_A.
    \end{split}
\end{equation}
The amount of second doses is analogously a square wave in counter phase with respect to the amount of first doses.
\end{itemize}
Firstly, we compare the solutions employing  the initial guess IG1. As for the minimization of infected (refer to Figures \ref{InfIG2} and \ref{poliInfDifferentR0IG2}) the optimal solutions concur in reducing the amount of administrations to the oldest (over sixties) and to the youngest (under twenties) age-classes. Moreover, as far as the initial reproduction number increases, the amount of administrations decreases. However, the amount of doses tends to increase with time as far as the outbreak runs out. In presence of a severe outbreak (Case 3, $\mathcal{R}_0 = 1.3$), the optimal vaccination roll-out indicates to administer over the 60\% of doses to the (20$\div$39) age-classes, and the remaining part to the (40$\div$59), leading to the same conclusions of Subsection \ref{MinInf}. However, the amount of doses allocated to the (20$\div$39) decreases with the depletion of the epidemic. In this case, with the optimal administrations the solution reaches a peak of infected which is dampened of approximately the half with respect to the initial guess, and it is advanced of nearly 23 days (see Figure \ref{InfIG2}). Instead, considering the minimization of deceased (Figures \ref{DecIG2} and \ref{poliDecDifferentR0IG2}) the optimal roll-out suggests to decrease the administrations to the (0$\div$19) and increase the amount of administrations to the oldest ones. However, in case of a sever outbreak a significant amount of allocations is associated to the (20$\div$39) and to the (40$\div$59), meaning that in this case it is difficult to contain the amount of deaths without reducing the spread of the disease carried out by the  socially active population. This result is unexpected, since the amount of doses distributed to those age-stratifications deficits the allocations to the most fragile, \textit{i.e.} the over eighties. The optimal solutions related to the severe outbreak allows to reduce deaths in the initial solution of nearly 23 thousands units in the major outbreak case, 8 thousands in the bifurcation one and almost 1.8 thousands in the case of a minor outbreak. Concerning the minimization of hospitalized (Figures \ref{HosIG2} and \ref{poliHosDifferentR0IG2}) in the three cases the optimal vaccination campaigns agree in allocating the majority of available doses to the (20$\div$59) and to the elderly as in Subsection \ref{MinHosp}. However, the administrations to the most socially active population (20$\div$59) increases with the initial reproduction number. With the optimal solution the evident peak of hospitalizations of the initial guess in the case of a major outbreak is anticipated of nearly twenty days, whilst in the case of a minor outbreak the different vaccination policy does not lead to significant reduction in the curve of hospitalizations.
\begin{figure}[H]
\minipage{0.32\textwidth}
  \includegraphics[width=\linewidth]{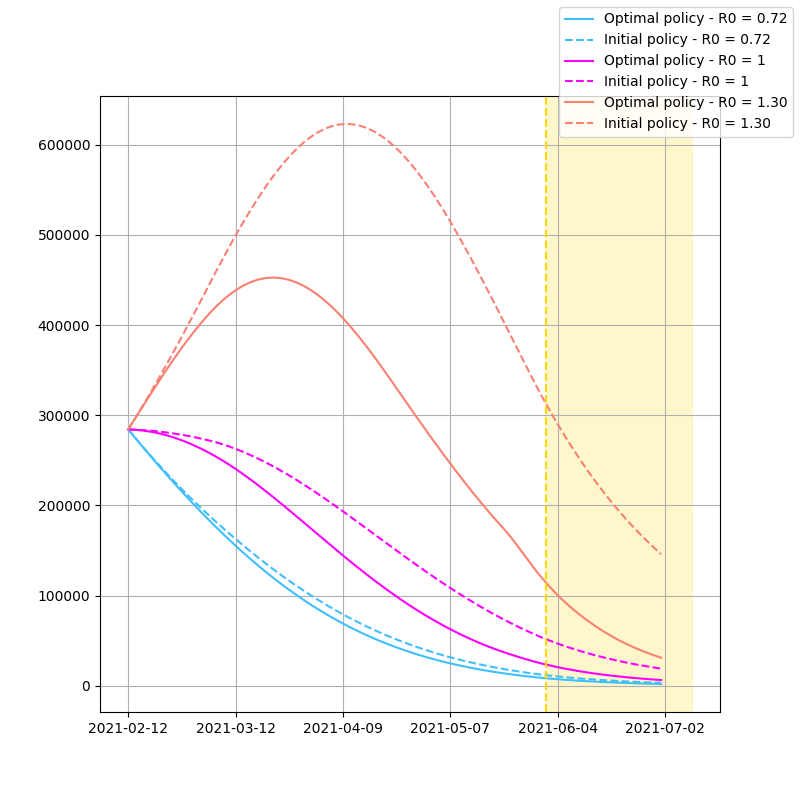}
  \caption{Evolution of infected of the two solutions obtained minimizing infected starting from IG1 with three different $\mathcal{R}_0$.}\label{InfIG2}
\endminipage\hfill
\minipage{0.32\textwidth}
    \includegraphics[width=\linewidth]{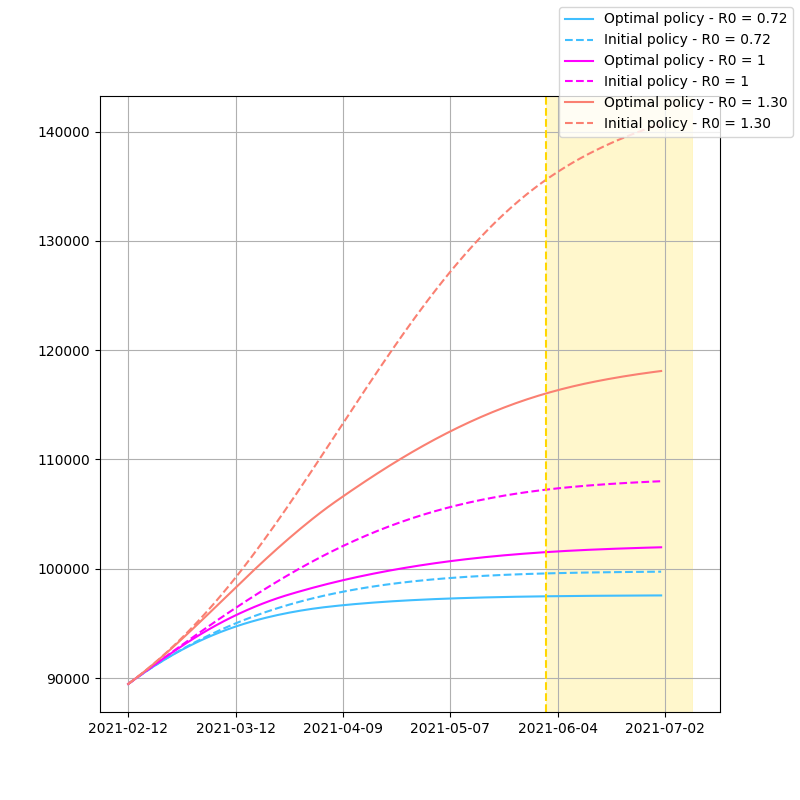}
  \caption{Evolution of deceased of the two solutions obtained minimizing deceased starting from IG1 with three different $\mathcal{R}_0$.}\label{DecIG2}
\endminipage\hfill
\minipage{0.32\textwidth}%
  \includegraphics[width=\linewidth]{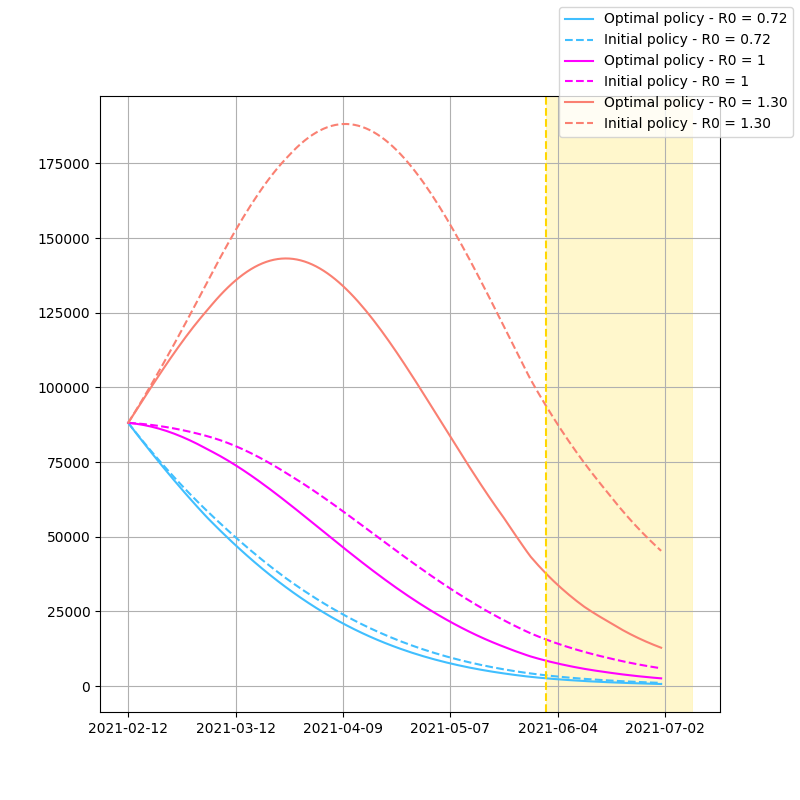}
  \caption{Evolution of hospitalized of the two solutions obtained minimizing hospitalized starting from IG1 with three different $\mathcal{R}_0$.}\label{HosIG2}
\endminipage
\end{figure}

\begin{figure}[H]
    \centering
    \includegraphics[trim={5.5cm 0 5.5cm 0},width=\textwidth]{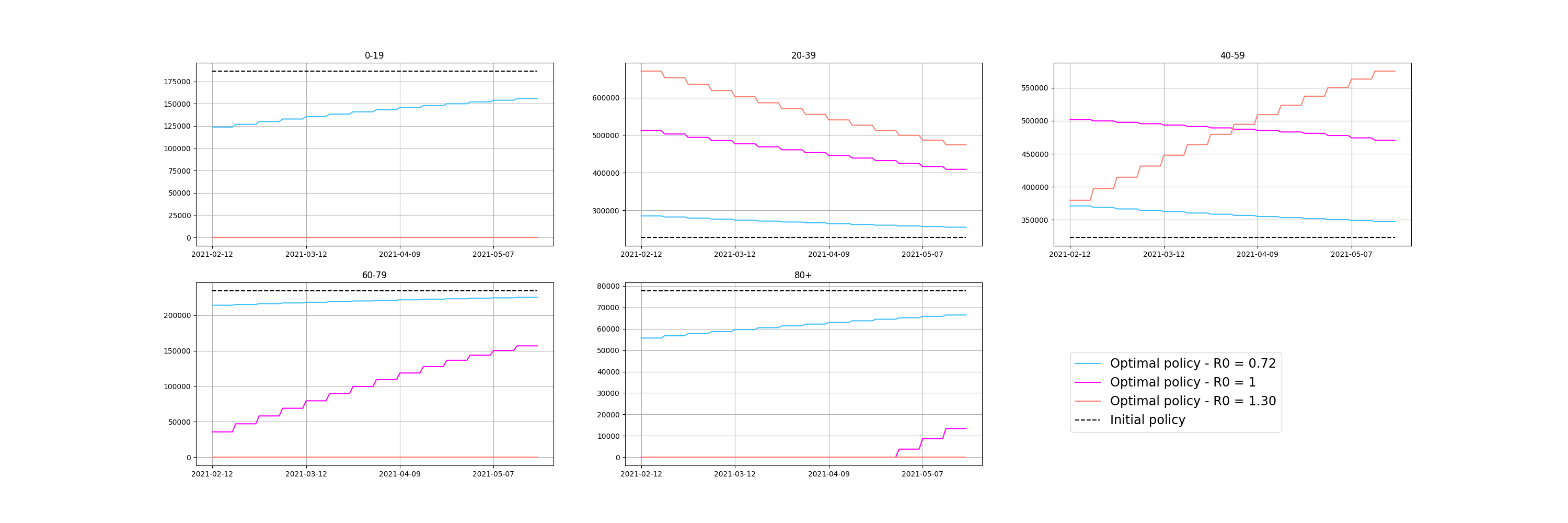}
    \caption{Weekly amount of doses delivered for each age-stratification in the solutions minimizing infected starting from IG1 and considering the three different outbreaks.}
    \label{poliInfDifferentR0IG2}
\end{figure}
\begin{figure}[H]
    \centering
    \includegraphics[trim={5.5cm 0 5.5cm 0},width=\textwidth]{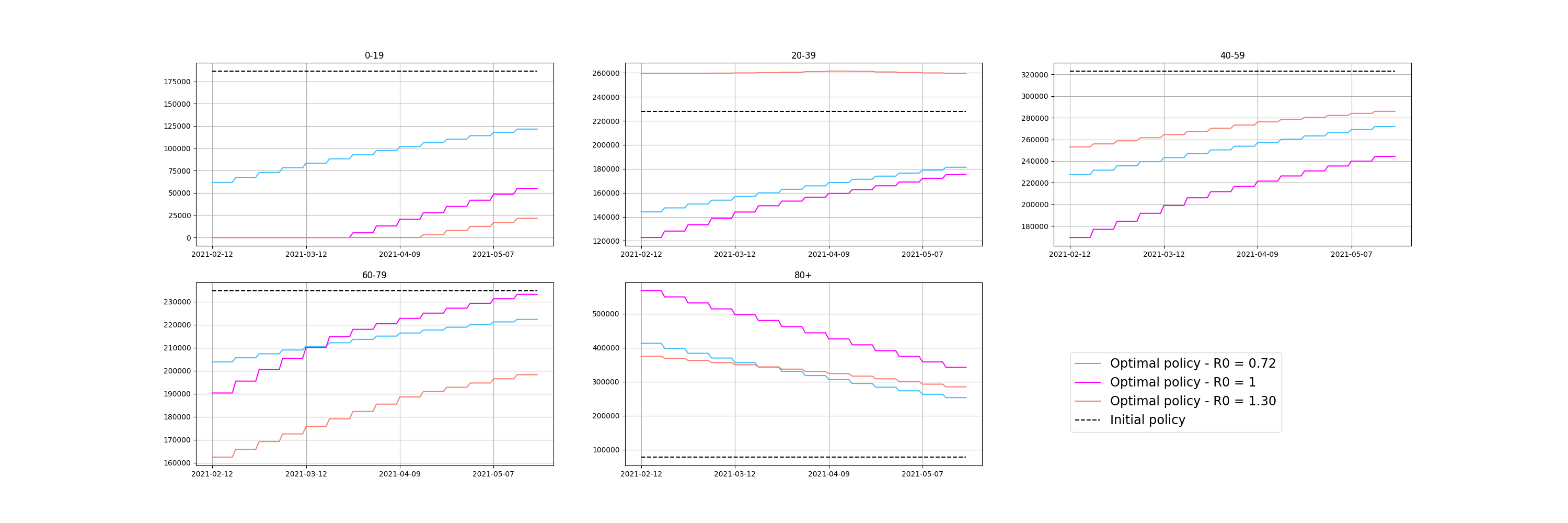}
    \caption{Weekly amount of doses delivered for each age-stratification in the solutions minimizing deceased starting from starting from IG1 and considering the three different outbreaks.}
    \label{poliDecDifferentR0IG2}
\end{figure}
\begin{figure}[H]
    \centering
    \includegraphics[trim={5.5cm 0 5.5cm 0},width=\textwidth]{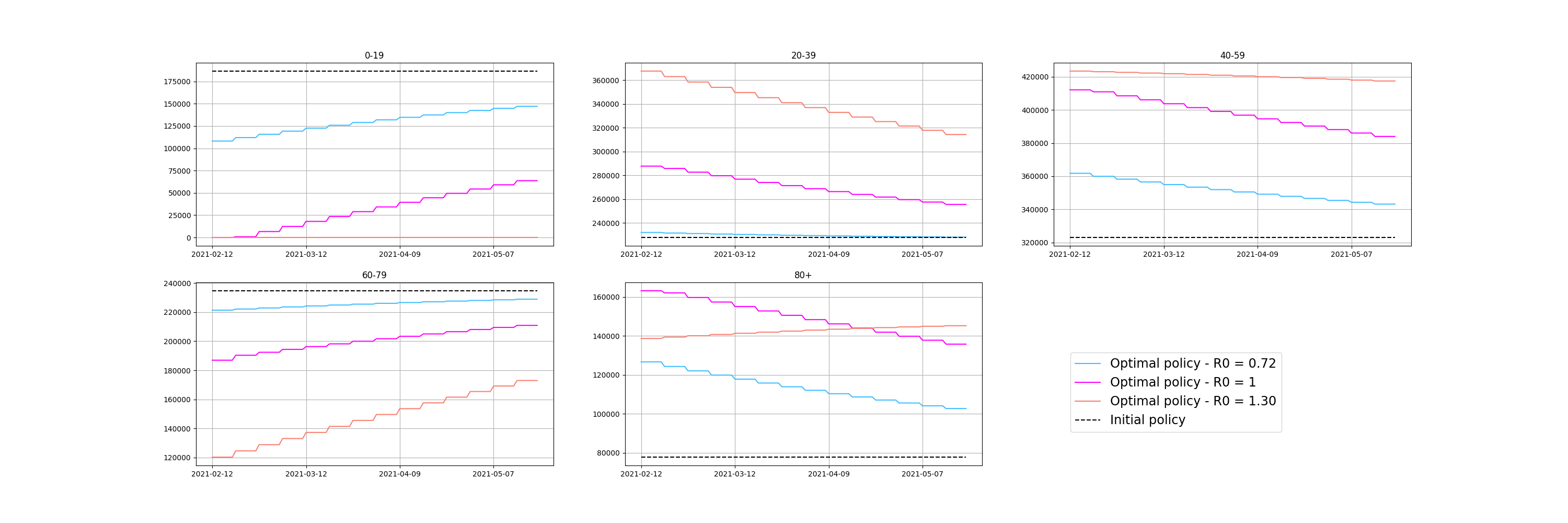}
    \caption{Weekly amount of doses delivered for each age-stratification in the solutions minimizing hospitalized starting from starting from IG1 and considering the three different outbreaks.}
    \label{poliHosDifferentR0IG2}
\end{figure}

We now consider the initial guess IG2. As a general remark, we expect that policies maximizing the amount of doses to those who are more likely to contract the disease with severe symptoms, as the $IFR$-based initial guess IG2, is closer to the solution minimizing deceased with respect to the former. We first focus on the case of minimization of infected (see Figures \ref{InfIG3} and \ref{poliInfDifferentR0IG3}), and we notice that the optimal policies consists in diminishing the amount of doses to the over sixties, increasing the (20$\div$59)'s allocations, and that the growth (respectively reduction) amounts depends on the respective value of the initial reproduction number. Indeed, administrations to the (20$\div$59) age classes increase as far as the reproduction number grows. In the case of a severe outbreak with the optimal solution we attain an early peak of nearly 20 days, and also the value at the peak has been lowered up to 451 thousands. On the other hand, the optimal allocation of doses for minimizing deaths tends to increase the amount of first doses to the (20$\div$59) to the detriment of the over sixties. Indeed, in the initial policy administrations to individuals in the (20$\div$39) and (40$\div$59) age-classes have been almost completely neglected due to the proportionality of doses with the $\{IFR_i\}_i$.
The only optimal solution which significantly decreases the amount of deceased is the one related to the severe outbreak (reduction of  nearly 4330 deaths). Lastly, for what concerns the optimal administrations for minimizing hospitalizations the solution confirms that a sufficient amount of doses (almost 68\%) has to be allocated to the (20$\div$59). In this way, for the case of a major outbreak there is the possibility of anticipating and reducing the peak of nearly two months. Actually, in the case of a minor outbreak the algorithm hardly moves from the initial guess, which seems to be almost optimal.
\begin{figure}[H]
\minipage{0.32\textwidth}
  \includegraphics[width=\linewidth]{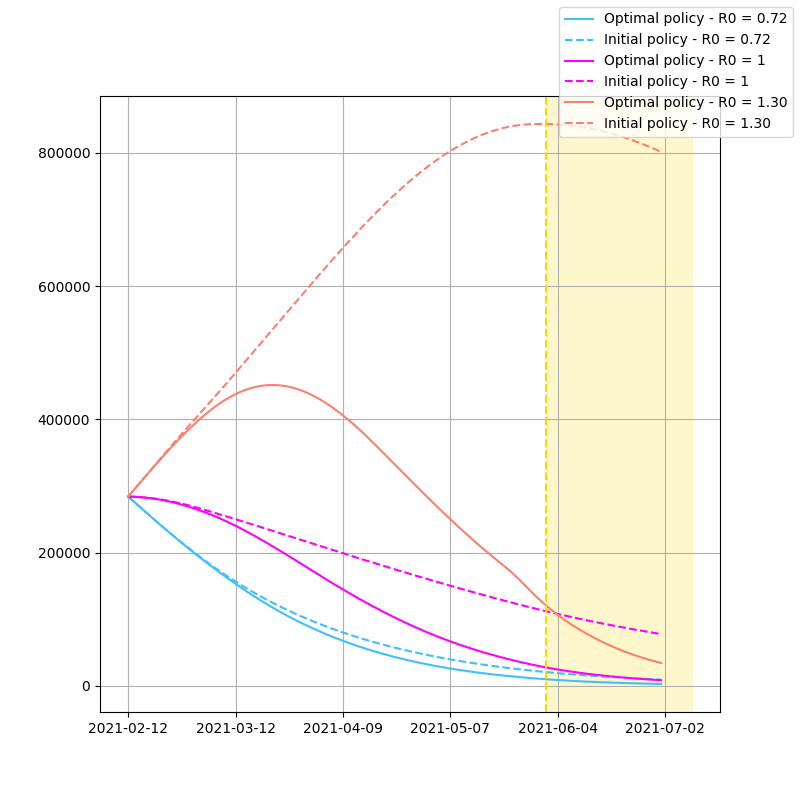}
  \caption{Evolution of infected of the two solutions obtained minimizing infected starting from IG2 with three different $\mathcal{R}_0$.}\label{InfIG3}
\endminipage\hfill
\minipage{0.32\textwidth}
    \includegraphics[width=\linewidth]{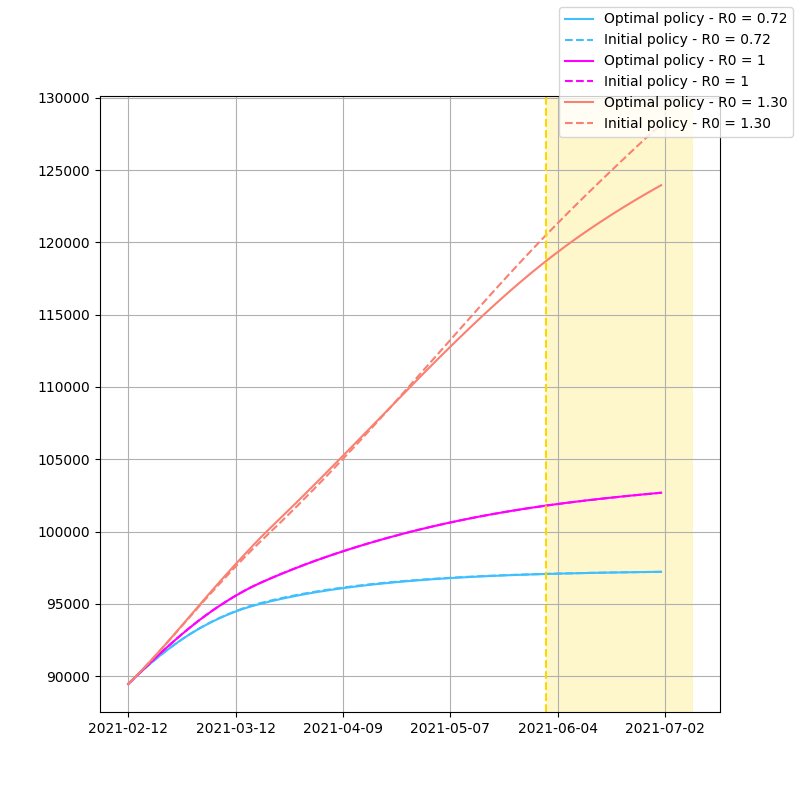}
  \caption{Evolution of deceased of the two solutions obtained minimizing deceased starting from IG2 with three different $\mathcal{R}_0$.}\label{DecIG3}
\endminipage\hfill
\minipage{0.32\textwidth}%
  \includegraphics[width=\linewidth]{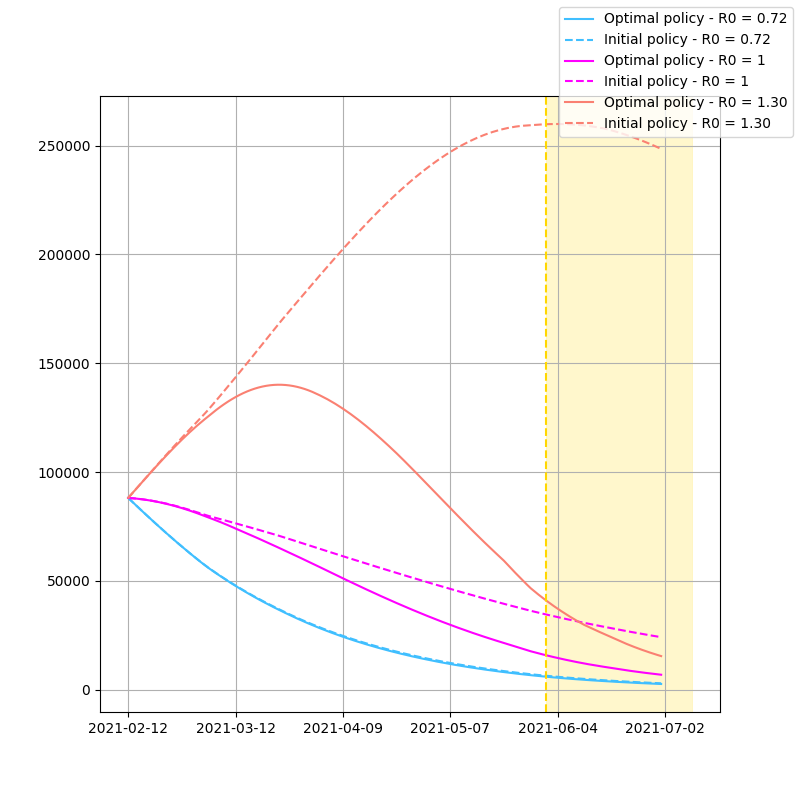}
  \caption{Evolution of hospitalized of the two solutions obtained minimizing hospitalized starting from IG2 with three different $\mathcal{R}_0$.}\label{HosIG3}
\endminipage
\end{figure}

\begin{figure}[H]
    \centering
    \includegraphics[trim={5.5cm 0 5.5cm 0},width=\textwidth]{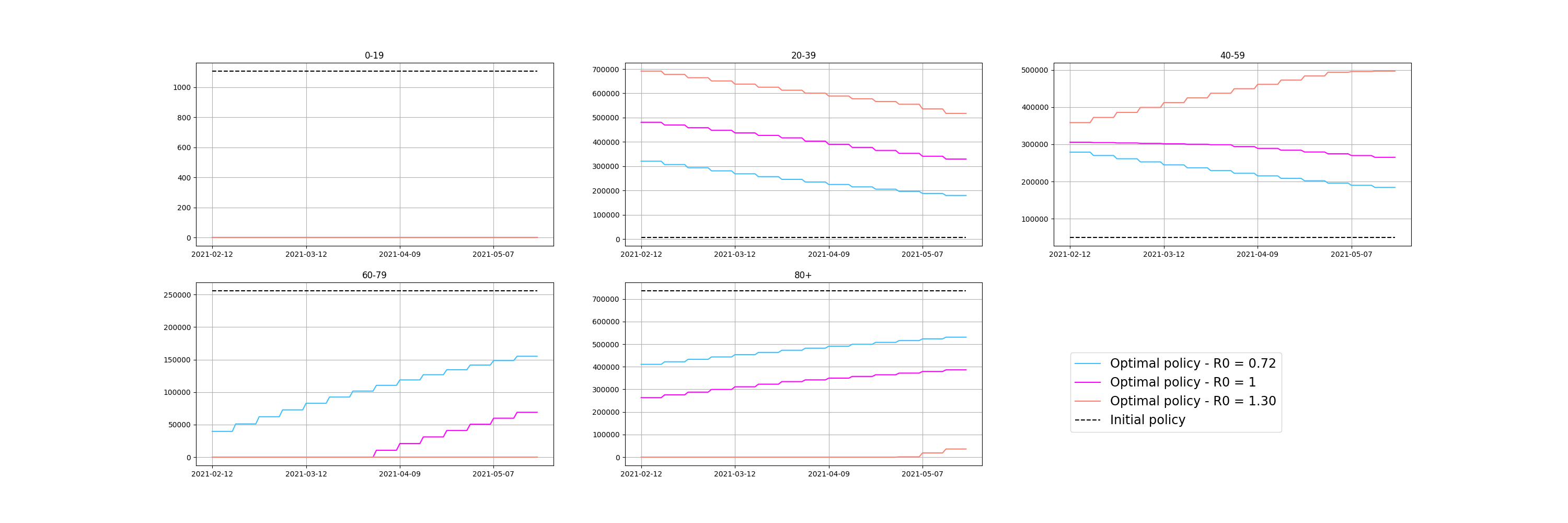}
    \caption{Weekly amount of doses delivered for each age-stratification in the solutions minimizing infected starting from IG2 and considering the three different outbreaks.}
    \label{poliInfDifferentR0IG3}
\end{figure}
\begin{figure}[H]
    \centering
    \includegraphics[trim={5.5cm 0 5.5cm 0},width=\textwidth]{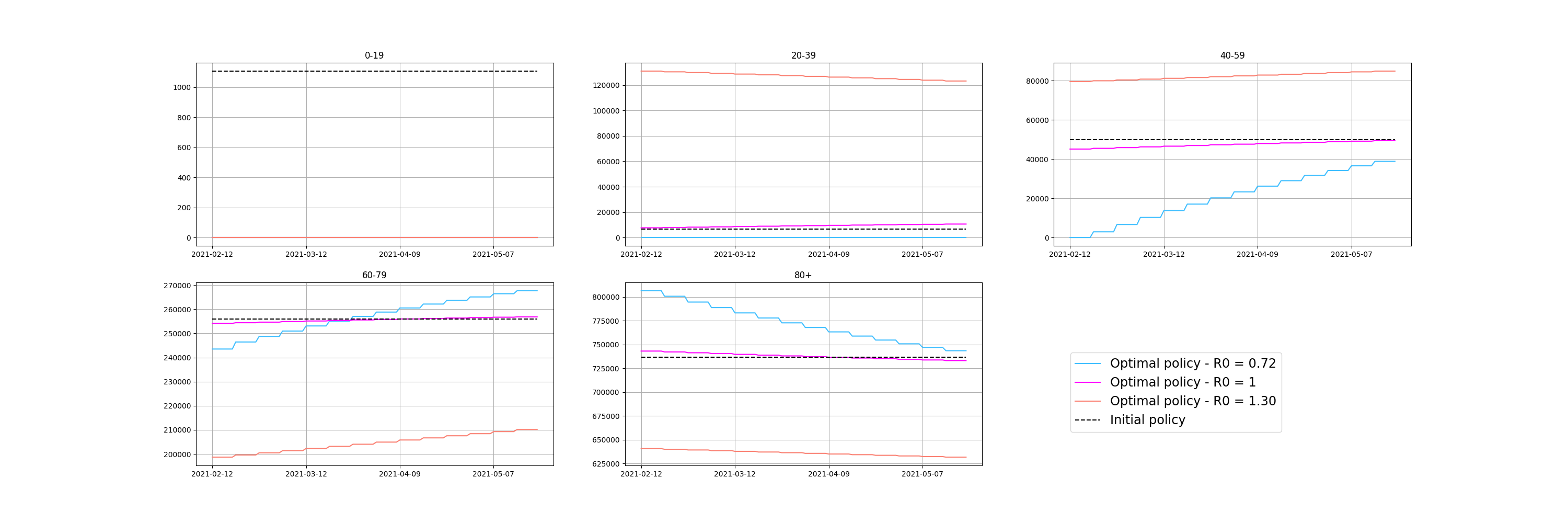}
    \caption{Weekly amount of doses delivered for each age-stratification in the solutions minimizing deceased starting from IG2 and considering the three different outbreaks.}
    \label{poliDecDifferentR0IG3}
\end{figure}
\begin{figure}[H]
    \centering
    \includegraphics[trim={5.5cm 0 5.5cm 0},width=\textwidth]{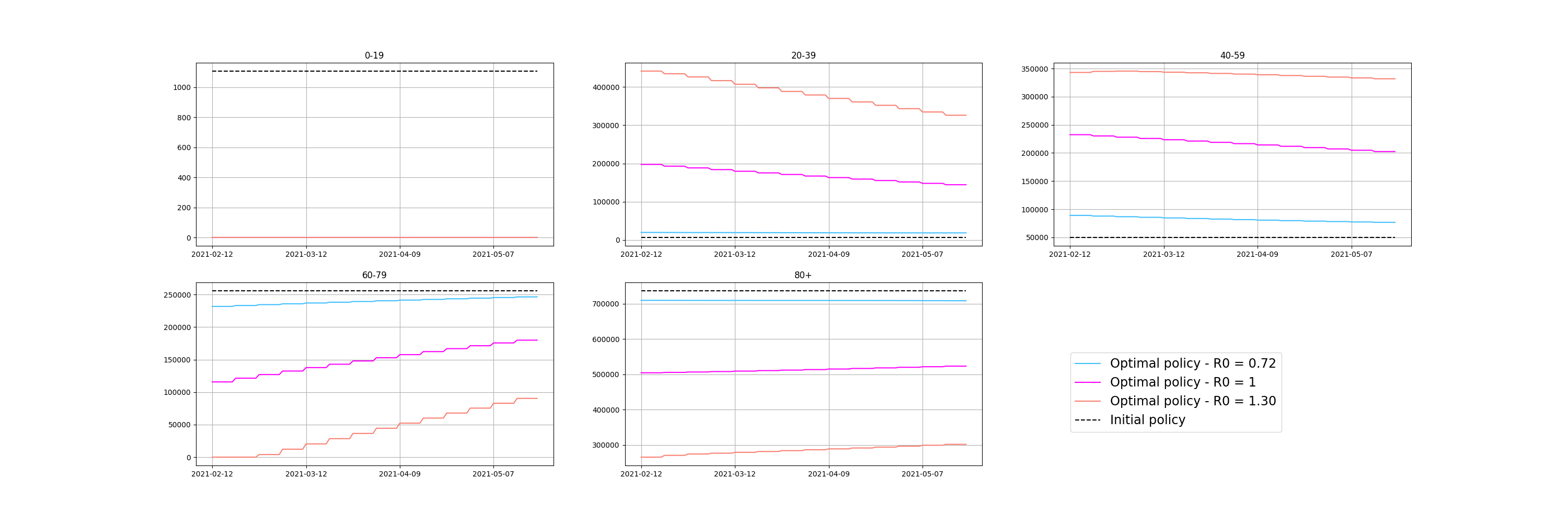}
    \caption{Weekly amount of doses delivered for each age-stratification in the solutions minimizing hospitalized starting from IG2 and considering the three different outbreaks.}
    \label{poliHosDifferentR0IG3}
\end{figure}

Finally, we run the numerical tests employing the initial guess IG3. From Figures \ref{InfSquare}-\ref{HosSquare} we notice that all the optimal vaccination campaigns regardless of cost functionals tend to retain the square wave behavior of the initial guess. However, the absolute amount of doses is different depending on the cost functional and on the value of the initial reproduction number. Specifically, the optimal policy obtained from the minimization of infected suggests to allocate more doses to the (20$\div$59), and in the case of a severe outbreak to privilege administrations to the (20$\div$39) with almost 700 thousands doses weekly (see Figure \ref{InfSquare}). On the other hand, the optimal solution from the minimization of  deceased increases the amount of administrations to the most fragile, and even to the (20$\div$39) in the case of a major outbreak, in agreement with what has been previously remarked (Figure \ref{DecSquare}). Finally, the optimal vaccination campaign for minimizing hospitalisations counsels to administer more doses to the over eighties and to the (20$\div$59), as far as the initial reproduction number grows. The solution minimizing hospitalized individuals is the one closer to the initial policy with respect to the other cases starting with the initial square wave policy. 
\begin{figure}[H]
\minipage{0.32\textwidth}
  \includegraphics[width=\linewidth]{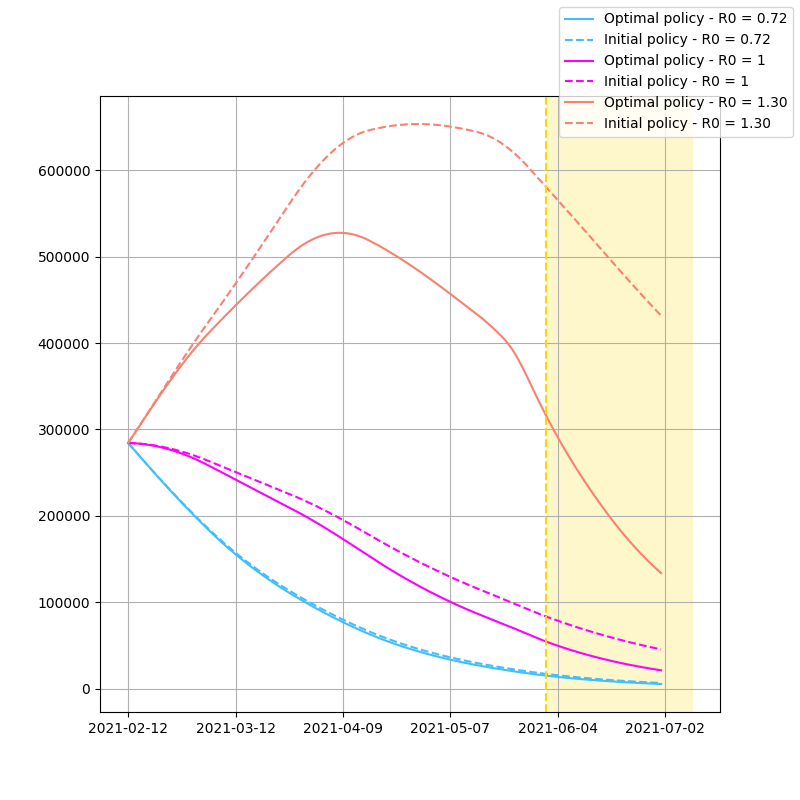}
  \caption{Evolution of infected of the two solutions obtained minimizing infected starting from IG3 with three different $\mathcal{R}_0$.}\label{InfSquare}
\endminipage\hfill
\minipage{0.32\textwidth}
    \includegraphics[width=\linewidth]{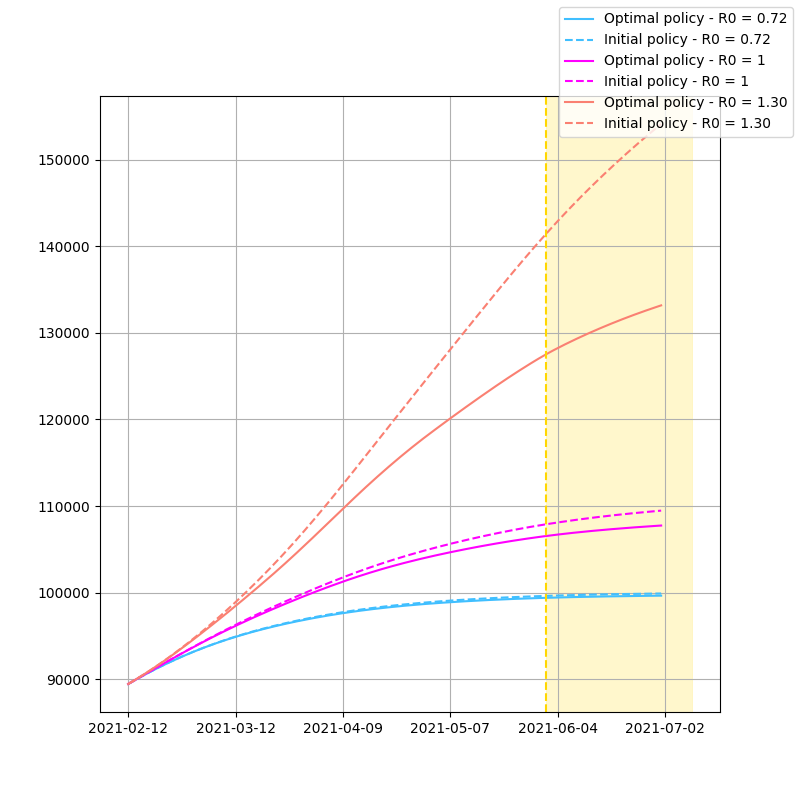}
  \caption{Evolution of deceased of the two solutions obtained minimizing deceased starting from IG3 with three different $\mathcal{R}_0$.}\label{DecSquare}
\endminipage\hfill
\minipage{0.32\textwidth}%
  \includegraphics[width=\linewidth]{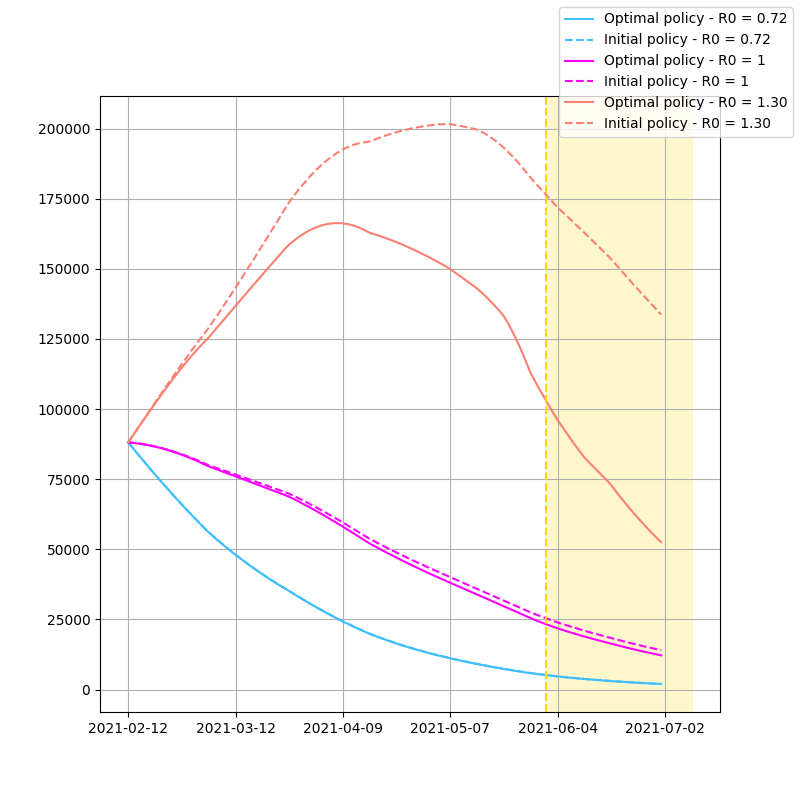}
  \caption{Evolution of hospitalized of the two solutions obtained minimizing hospitalized starting from IG3 with three different $\mathcal{R}_0$.}\label{HosSquare}
\endminipage
\end{figure}

\begin{figure}[H]
    \centering
    \includegraphics[trim={5.5cm 0 5.5cm 0},width=\textwidth]{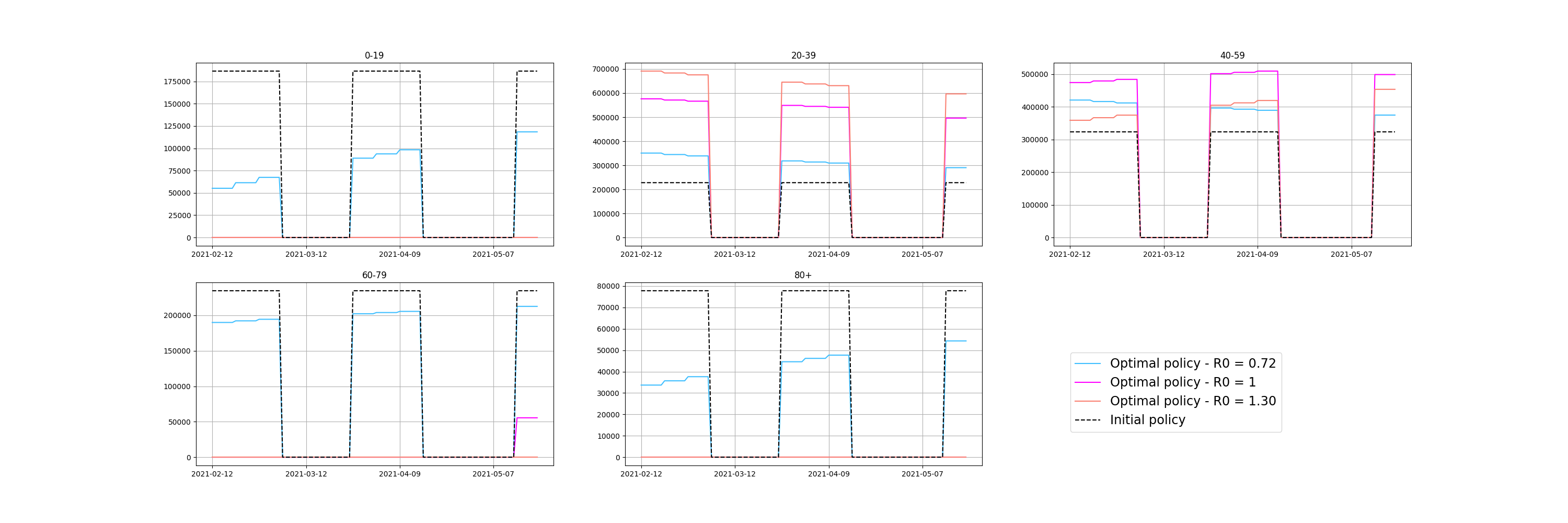}
    \caption{Weekly amount of doses delivered for each age-stratification in the solutions minimizing infected starting from IG3 and considering the three different outbreaks.}
    \label{poliInfDifferentR0square}
\end{figure}
\begin{figure}[H]
    \centering
    \includegraphics[trim={5.5cm 0 5.5cm 0},width=\textwidth]{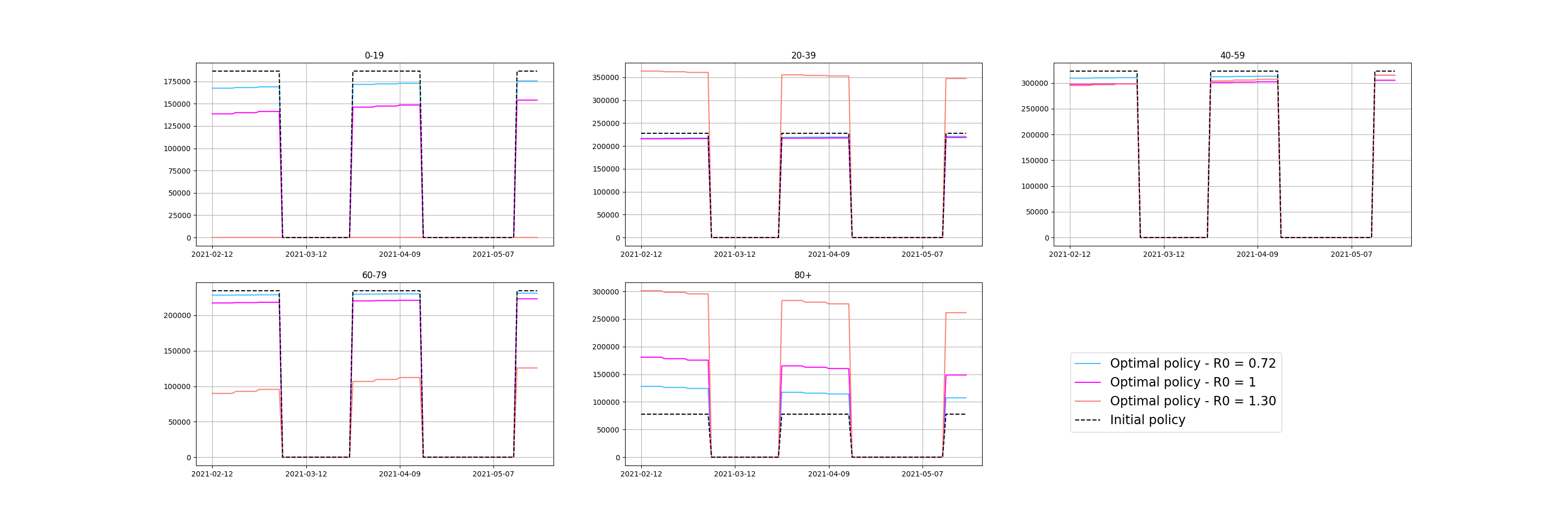}
    \caption{Weekly amount of doses delivered for each age-stratification in the solutions minimizing deceased starting from IG3 and considering the three different outbreaks.}
    \label{poliDecDifferentR0square}
\end{figure}
\begin{figure}[H]
    \centering
    \includegraphics[trim={5.5cm 0 5.5cm 0},width=\textwidth]{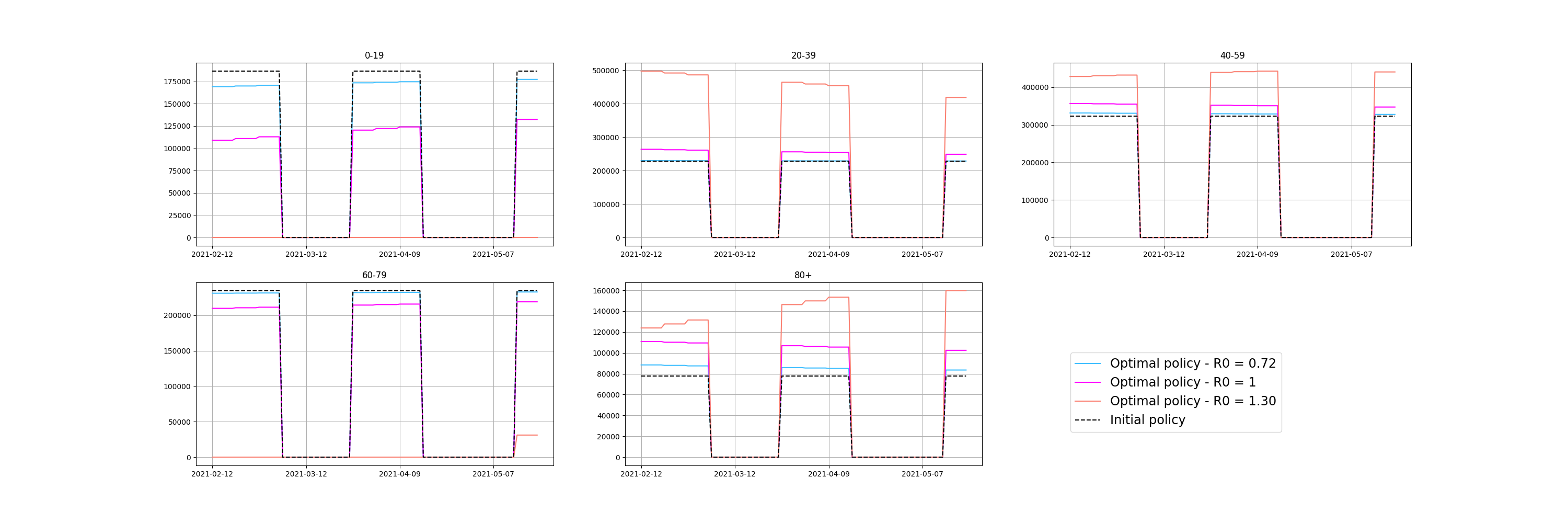}
    \caption{Weekly amount of doses delivered for each age-stratification in the solutions minimizing hospitalized starting from IG3 and considering the three different outbreaks.}
    \label{poliHosDifferentR0square}
\end{figure}

In Table \ref{tableSpared} we report, for different values of $\mathcal{R}_0$ and for different choice of the initial guess, the reductions of infected, deceased and hospitalized individuals with respect to the corresponding value obtained from the initial guess. We notice that as far as the reproduction number increases the projected gradient method is able to retrieve optimal policies corresponding to relevant improvement in the compartments to be minimized.
\begin{table}[H]
\centering
\begin{tabular}{c|c|c|c|c}
\multicolumn{1}{l|}{} &
  \multicolumn{1}{l|}{} &
  \multicolumn{1}{l|}{\textbf{IG1}} &
  \multicolumn{1}{l|}{\textbf{IG2}} &
  \multicolumn{1}{l}{\textbf{IG3}} \\ \hline
 & $\mathcal{R}_0 = 0.72$ & 5623                         & \color{ForestGreen}8711\color{black}  & 2003  \\
 \cline{2-5} 
 & $\mathcal{R}_0 = 1.01$    & 29902                        & \color{ForestGreen}53140\color{black} & 20080 \\ \cline{2-5} 
\multirow{-3}{*}{\textbf{\begin{tabular}[c]{@{}c@{}}Minimizing $\mathcal{J}_I$\\ Average infected saved\end{tabular}}} &
  $\mathcal{R}_0 = 1.30$ &
  158184 &
  \color{ForestGreen}388583\color{black} &
  147602 \\ \hline
 & $\mathcal{R}_0 = 0.72$ & \color{ForestGreen}2175\color{black} & 30                            & 250   \\ \cline{2-5} 
 & $\mathcal{R}_0 = 1.01$    & \color{ForestGreen}6051\color{black} & 469                           & 1732  \\ \cline{2-5} 
\multirow{-3}{*}{\textbf{\begin{tabular}[c]{@{}c@{}}Minimizing $\mathcal{J}_D$\\ Total Deceased saved\end{tabular}}} &
  $\mathcal{R}_0 = 1.30$ &
  \color{ForestGreen}22671\color{black} &
  4330 &
  21156 \\ \hline
 & $\mathcal{R}_0 = 0.72$ & \color{ForestGreen}1767\color{black} & 325                           & 57    \\ \cline{2-5} 
 & $\mathcal{R}_0 = 1.01$    & 7602                         & \color{ForestGreen}11242\color{black} & 1423  \\ \cline{2-5} 
\multirow{-3}{*}{\textbf{\begin{tabular}[c]{@{}c@{}}Minimizing $\mathcal{J}_H$\\ Average hospitalized saved\end{tabular}}} &
  $\mathcal{R}_0 = 1.30$ &
  42444 &
  \color{ForestGreen}116299\color{black} &
  39598 \\ 
\end{tabular}
\caption{Amount of saved infected, deceased and hospitalized individuals with respect to the corresponding value obtained from the initial guess (for different values of the 
initial reproduction number and for different choice of the initial guess).}
\label{tableSpared}
\end{table}

\section{Conclusions}
\label{Conclusions}
In this work, we presented an epidemic age-stratified compartmental model and formulated an optimal control problem for obtaining the vaccine distributions across ages minimizing specific goals. The model governing the optimal control problems is an age-stratified model, where the six compartments are split  into five age-classes. It incorporates  all the necessary features to effectively deal with the COVID19-vaccination campaign, and it can consistently incorporate the available data from DPC \cite{dpcdata, dpcdatavax}. Following the Pontryagin-KKT approach, we derived the optimal control formulation in order to minimize infectious, deceased or hospitalized caused by the disease. Finally, we detailed the iterative numerical strategy (projected gradient method) employed for the solution of the minimization problem.
We applied this architecture at the italian scenario during the period February 12th, 2021- June 1st, 2021. 

A first set of optimal vaccination problems has been solved by employing as initial guess for the iterative scheme the vaccination campaign actually implemented in Italy during the first six month of 2021. The resulting optimal vaccination policies for reducing the amount of infected individuals suggest to administer more doses to people ageing (20$\div$59), reducing the other age-classes. This result can be explained in light of the contact matrix weighted by age-dependent susceptibility, and assessing that individuals in these age-classes are the ones more active from a social perspective. However, this solution does not coincide with the one minimizing deceased individuals, which indeed requires to increase the amount of doses for the over eighties, who are the ones more prone to contract the disease in a more severe way. Finally, the solution minimizing the hospitalized individuals tends to increase the amount of doses to the age-class (20$\div$59) as in the minimizing-infected case. 

The second set of simulations compared the vaccination policies from the previous set of experiments with the optimal policies obtained by  employing an initial guess that allocates to each age class an amount of doses proportionally to the numerousness of each age-compartment. In this case, we qualitatively obtained the same guidelines for the vaccination policy that we obtained in the first set of experiments,  though the dependence of the solution on the initial guess (typical of gradient-based methods) clearly stands.

The last set of simulations explored different \textit{what-if} scenarios. In particular, we assumed  to deal with: (a) constant total amount of delivered doses in each week; (b) constant transmission rates corresponding to different level of severity of the outbreak. 
The obtained results showed the robustness of the proposed methodology to deal with the planning of a vaccination campaign in a variety of different epidemiological situations.

\newpage
\beginsupplement
\begin{center}
{\LARGE \textbf{Supporting information}}
\end{center}
\section{Parameters Setting for the Italian scenario}
In Section 3, the parameters involved in the model have been set as follows:
\begin{itemize}
    { \item $I = (0, T_f], \, T_f=151$ (days): the epidemiological model is set in Italy on the period ranging from January 1\textsuperscript{st}, 2021 to June 1\textsuperscript{st}, 2021 (\textit{i.e.} the six-month period that starts with the beginning of the vaccination campaign and concurrent with the third pandemic wave).}  
    \item $\beta \in (0,1)$: transmission rate, depending on the implemented Non-Pharmaceutical Interventions (NPIs) and virus transmissibility. It is assumed to be constant across all ages as in \cite{marziano2021retrospective};
    \item $\sigma_V, \sigma_W \in (0,1)$: vaccine effectiveness on transmissibility after administration of first dose (the former) or completing the cycle (the latter). It can be interpreted as the ratio of transmissibility between vaccinated individuals and unvaccinated ones. The value 0 means that the vaccine is fully effective, 1 totally ineffective;
    \item $\theta_V, \theta_W \in (0,1)$: vaccine effect on mortality after administration of first dose (the former) or completing the cycle (the latter). It can be interpreted as the ratio of probability of getting severe symptoms between vaccinated individuals and unvaccinated ones. The value 0 means that the vaccine is fully effective, 1 totally ineffective;
    \item $IFR_i$: age-dependant Infection Fatality Rate estimated starting from available data from Dipartimento di Protezione Civile Italiana as in \cite{parolini2021suihter};
    \item In the fatality function $f_i(S_i,V_i,W_i)$ the amount of days after the inoculation that we consider for reaching the complete vaccine effectiveness $t_a$ is fixex at 15 days;
    \item $C_{ik}$: $i,k$-th entry of the contact matrix, tracing back contacts  between ages starting from the POLYMOD surveys \cite{10.1371/journal.pmed.0050074} (see Figure \ref{polymodMatrix});
    \item $r_i$: susceptibilities to infection depending on age, drawn from \cite{zhang2020changes}, previously adopted by \cite{marziano2021retrospective};
    \item $\gamma$: recovery rate from the disease infection, which is maintained constant across ages. Since infectious individuals are supposed to exit from the correspondent compartment  with  flux $\gamma I_i$, the parameter $\gamma$ is interpreted as the inverse of the average time of recovery $t_R$. The distribution of recovery times is treated as a Gaussian distribution with mean $\bar{t}_R = 14.20$ (days) and variance $\sigma^2 = 5.94$ ($\textrm{days}^2$) as in \cite{voinsky2020effects}. We actually draw the posterior distribution of {the parameter $\gamma$ from} the calibration stage (see \ref{Calibration} for more details);
    \item $N_i$: {Number} of individuals in the i-th age stratification;
    \item $\mu_R$: natural waning immunity rate, taking into account plausible reinfections coming from previously-recovered individuals. This parameter has not shown peculiar dependance on age and has been retrieved from \cite{hansen2021assessment, lumley2021antibody};
    \item $I_{u,i}(t) = (1 - \delta(t)) I_i(t)\; t \, \in \mathrm{I}$: approximate number of undetected individuals whose age falls in the $i$-th stratification at time $t$. As in \cite{parolini2022modelling}, we compute the detection rate $\delta$ from the ratio of the computed Case Fatality Ratio ($CFR$) and the theoretical Infectious Fatality Ratio ($IFR$) for COVID19 independently on ages;
    \item $\delta_w$: elapsing time among subsequent administrations to be chosen in the admissible set $\{21, 28, 35, 42\}$ days. In the Italian scenarios, it has been set at 21 days;
    \item $U_{1,i},\, U_{2,i},\, U_{R,i}$: daily amount of administered first doses, second doses and doses administered to the $i$-th age-class, respectively. The choice of these variables is coherent with actual implementation of the italian vaccination campaign: two consecutive doses to be administered for completing the vaccination cycle to Susceptible individuals, one single administration to Recovered ones. 
    {To reduce the computational complexity of the model we assume that the functions $U_{1,i}(t),\, U_{2,i}(t),\, U_{R,i}(t)$ are piecewise constant (constant on each week) and the weekly value of administrations  is supposed to be equally distributed among each day of the week. In Figure \ref{vax_calib} we present the history of the daily administrations in Italy during the period of interest, while in Figure \ref{vax_calib_av} the corresponding averages over the weeks are reported.} 
\end{itemize}
\begin{figure}
    \centering
    \includegraphics[scale=0.3]{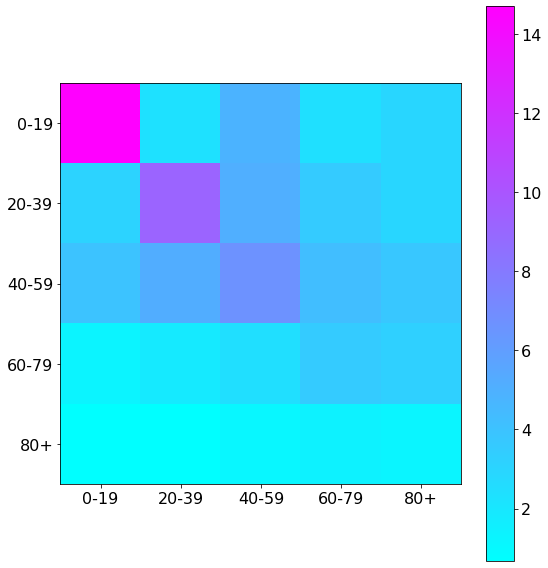}
    \caption{Contact matrix retrieved by \cite{10.1371/journal.pmed.0050074} for the specific age-classes contemplated in the model.}
    \label{polymodMatrix}
\end{figure}
\begin{remark}
During the Italian vaccination campaign, the mRNA-based vaccines Pfizer/BioNTech Comirnaty (BNT162b2) and SpikeVax (previously COVID-19 Moderna mRNA--1273) have been administered more extensively. The values of the parameters $\sigma$ and $\theta$ reported in Table \ref{ParamTab} have been deduced assuming to deal with a Pfizer/BioNTech kind of vaccine. The effectiveness of the vaccine on reducing both transmissibility and severity of symptoms has been previously assessed through medical trials (see \textit{e.g.} \cite{bernal2021effectiveness}) and then posteriorly confirmed through the available epidemic data \cite{tartof2021effectiveness}. However, given the potential relevant impact of $\sigma$ and $\theta$ on the optimality of the vaccination campaign, in \ref{Sensitivity} we carry out a sensitivity analysis of these parameters.  
\end{remark}

\begin{figure}
\minipage{0.45\textwidth}
    \centering
    \includegraphics[scale = 0.28]{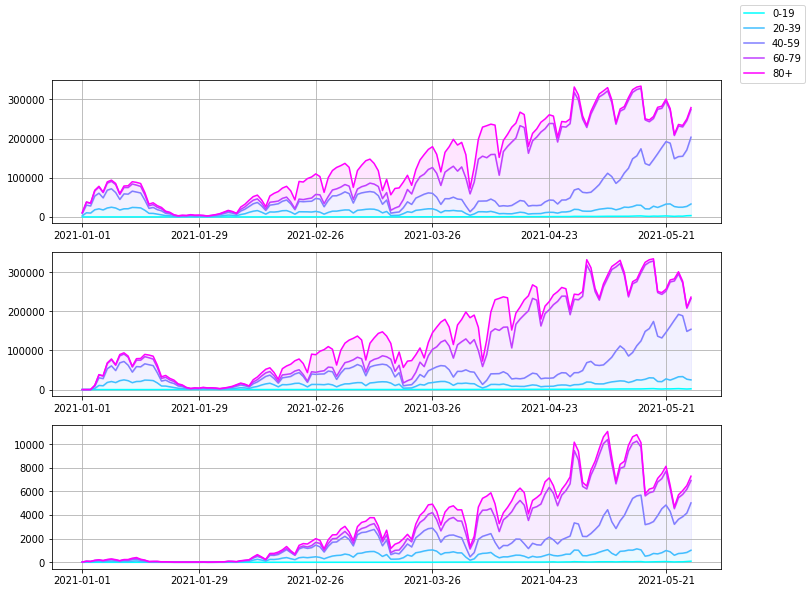}
    \caption{Amount of first doses, second doses and doses administered to recovered during the time frame of interest retrieved from DPC data \cite{dpcdatavax}.}
    \label{vax_calib}
\endminipage
\hspace{5mm}
\minipage{0.45\textwidth}
    \centering
    \includegraphics[scale = 0.28]{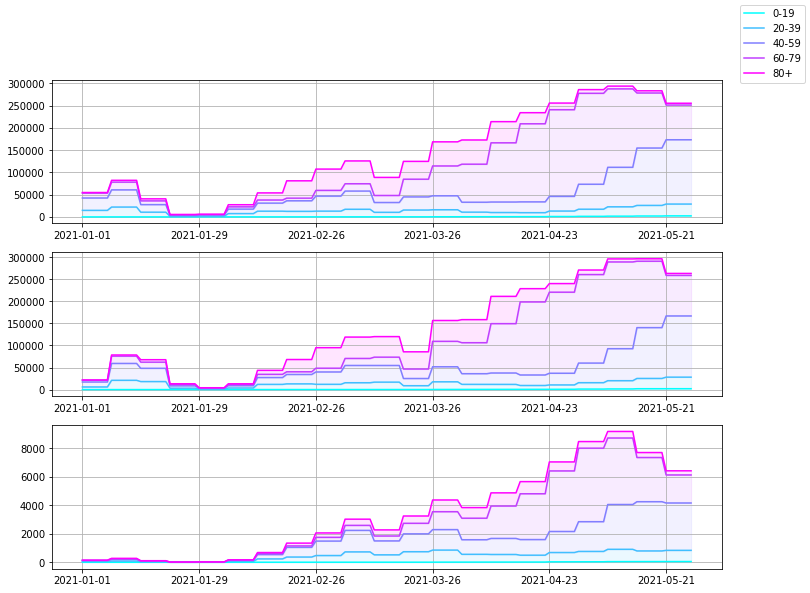}
    \caption{Amount of first doses, second doses and doses administered to recovered during the time frame of interest retrieved from DPC data \cite{dpcdatavax}, and then averaged per weeks.}
    \label{vax_calib_av}
\endminipage
\end{figure}

\begin{table}[]
\centering
\begin{tabular}{c|ccccc|c}
\textbf{Parameter} &
  \multicolumn{1}{c|}{\textbf{{(}0,19{)}}} &
  \multicolumn{1}{c|}{\textbf{{(}20,39{)}}} &
  \multicolumn{1}{c|}{\textbf{{(}40,59{)}}} &
  \multicolumn{1}{c|}{\textbf{{(}60,79{)}}} &
  \textbf{(80+)} & \textbf{Reference} \\ \hline
$\gamma$      & \multicolumn{5}{c|}{0.07}                                                                                                & \cite{voinsky2020effects}\\ \hline
$r_i$         & \multicolumn{1}{c|}{0.33} & \multicolumn{1}{c|}{1}    & \multicolumn{1}{c|}{1}    & \multicolumn{1}{c|}{1}      & 1.47  & \cite{zhang2020changes} \\ \hline
$IFR_i$ & \multicolumn{1}{c|}{1e-4} & \multicolumn{1}{c|}{6e-4} & \multicolumn{1}{c|}{4.5e-3} & \multicolumn{1}{c|}{2.3e-2} & 7.2e-2 & \cite{parolini2021suihter} \\ \hline
$\mu_R$       & \multicolumn{5}{c|}{0.006} & \cite{hansen2021assessment, lumley2021antibody} \\ \hline
$\sigma_V$    & \multicolumn{5}{c|}{0.21} & \cite{marziano2021retrospective}                                                                                               \\ \hline
$\sigma_W$    & \multicolumn{5}{c|}{0.21} & \cite{marziano2021retrospective}                                                                                               \\ \hline
$\theta_V$    & \multicolumn{5}{c|}{0.20} & \cite{marziano2021retrospective}                                                                                               \\ \hline
$\theta_W$    & \multicolumn{5}{c|}{0.037} & \cite{marziano2021retrospective}
             \\ 
\end{tabular}
\caption{Table of the parameters adopted during the calibration stage.}
\label{ParamTab}
\end{table}
{
The differential model has been endowed with proper initial conditions for each of the considered age-state compartment. The amount of initial deceased is deterministically fixed and it is retrieved by available data in \cite{dpcdata}. On the other hand, we introduce uncertainty on the initial conditions for the Susceptible, Infectious and Recovered classes. This choice is dictated by a different level of accuracy of the available data for the different compartments. 
Finally, vaccinated individuals belonging to $V$ and $W$ compartments have been set to zero since our period of interest 
ranges from January 1\textsuperscript{st},  2021 to June 1\textsuperscript{st}, 2021  (\textit{i.e.} the six-month period coinciding with the beginning of the vaccination campaign and therefore with the third pandemic wave).} Hence, we assume
\begin{equation}
    \begin{split}
        I_i(0) &= \zeta_{I,i} I_{D,i,0}\\
        R_i(0) &= \zeta_{R,i} R_{D,i,0}\\
        S_i(0) &= (N_i - (I_i(0) + R_i(0) + D_i(0))) \,\zeta_{S,i}
    \end{split}
    \label{UncInState}
\end{equation}
where $X_{D,i}$ is the value of the $X$-th compartment retrieved from the Italian DPC data on January 1\textsuperscript{st}, 2021 and $\zeta_{X,i}$ are independent and identical uniform distributions in [0.7, 1.3]. {The choice in \eqref{UncInState} allows to introduce prior uncertainty to be posteriorly reconstructed after the Monte Carlo Markov Chain stage of the calibration process (see \ref{Calibration})}.
\section{Calibration process}
\label{Calibration}
{Following \cite{parolini2022modelling}, we undergo a double-stage calibration for the transmission rate $\beta$.
The transmission rate embodies the effects of non-pharmaceutical interventions due to the Italian government's political decisions, updated on a weekly basis during the period of concern. Hence, the parameter is assumed to be piecewise-constant over one-week-long intervals ({\em time-phases} in the sequel).
Similarly, e.g., to \cite{marziano2021retrospective}), we assume that the transmission rate is independent of the age.
In view of the above assumptions, we need to calibrate 21 parameters representing the values of the transmission rate $\beta$ in each time-phase. 

Generally speaking, calibration problems can be interpreted as deterministic inverse optimal control problems \cite{comunian2020inversion, dantas2018calibration, marinov2020dynamics, marinov2014inverse}, or can be tackled through stochastic techniques \cite{daza2022bayesian, xu2021bayesian, korostil2013adaptive}. Alternatively, calibration problems can be solved by resorting to Machine-Learning based schemes ( see \textit{e.g.} \cite{wieczorek2020neural, zisad2021integrated} and  \cite{li2020predicting, zisad2021integrated, nabi2021forecasting, mohimont2021convolutional}). 
In this paper the calibration step is implemented by employing a classical two-stage process: (1) a Least Square (LS) phase is used for retrieving acceptable estimators of the parameters; (2) every output of the LS phase is then adopted as mean value for the corresponding prior distribution imposed during the classical Monte Carlo Markov Chain (MCMC) approach (we refer to \cite{rasmussen2011inference} for more details). For what concerns the calibration of the initial conditions \eqref{UncInState}, we impose uniform priors on the initial values of $S,I,R$ variables, while a Gaussian prior is employed for the calibration of the recovery time variable $t_R = \frac{1}{\gamma}$. During both LS and MCMC stages of the calibration, we minimize the adherence, for each age class $i$, of the deceased compartment $D_i$ to the public available data $D_{D,i}$}:
\begin{equation}
    \mathcal{E} = \sum_{i = 0}^{N_{ages}} \int_0^{T_f} (D_i(t) - D_{D,i}(t))^2 \, dt .
    \label{error}
\end{equation}
In Figure \ref{calibDec} the results of the calibration are shown. In particular, Figure \ref{calibDec}(a) represents the total amount of deceased during the first six months of 2021 (dashed line) and the behavior of the median value after the MCMC calibration, together with the credible interval of order 95\% (shaded area). In Figure \ref{calibDec}(b) we present the age repartition of deceased among age-classes. More precisely, the orange curve stands for the amount of deceased after the LS-calibration, the red one for the median of the posterior, the black dashed line represents the DPC data of deceased, while the shaded areas represents the 95\% credible intervals for each age-class.
Notice that, as the cumulative number of recorded deceased cannot decrease, the decreasing behaviour in Figure \ref{calibDec}(b) of the dashed line is probably due to reporting errors. {\color{black} The median trend of the total amount of deceased is in agreement with the actual trend of the detected deaths recorded by the DPC; the same variable lies in the 95\% credible interval posteriorly recovered, as showed in Figure \ref{calibDec}(b). 
Moreover, notice that the percentage repartition of deaths attained by the model with {\color{black} the calibrated transmission rate} agrees qualitatively with the age-repartition of deceased actually observed in Italy during the same period (see Figure \ref{calibDec}(c)). To further assess the reliability of the output of the calibration process, the calibrated model has also been validated by comparing the simulated amount of infected individuals with the corresponding data value. In particular, in Figure \ref{calibInf} the simulated detected infected individuals, computed as $I_{d,i}(t) = \delta I_i(t)$, $\delta$ being the detection rate, are compared with the actual positive individuals recorded by \cite{dpcdata}. As the the number of infected individuals are not included in \eqref{error}, the results of Figure \ref{calibInf} show that the performance of the implemented model with the calibrated parameter can be considered fully satisfactory. In this respect, we also note that Figure \ref{calibInf}(b) shows an excellent accordance, in terms of age-percentage repartition of infected, between the public available data and the output of our calibrated model.}

\minipage{0.65\textwidth}
\begin{figure}[H]
    \centering
    \subfloat[$\beta$ at week 0.]{\includegraphics[scale=0.2]{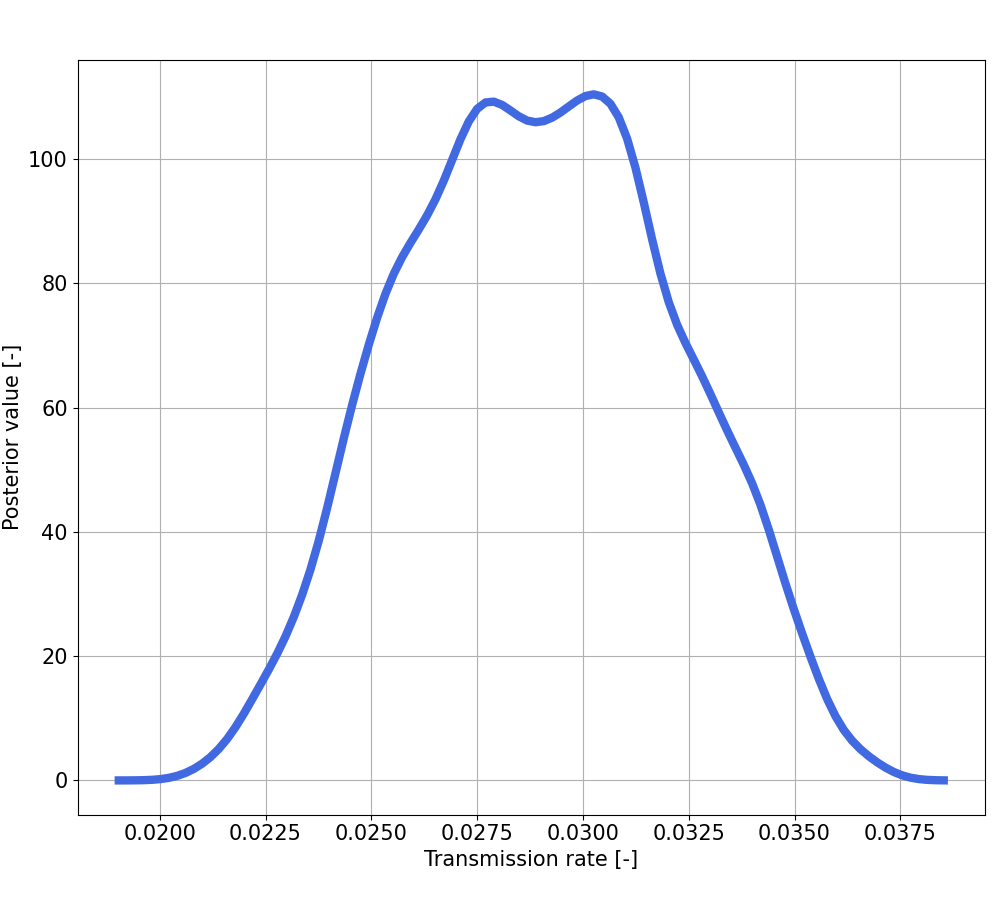}}
    \subfloat[$\beta$ at week 9.]{\includegraphics[scale=0.2]{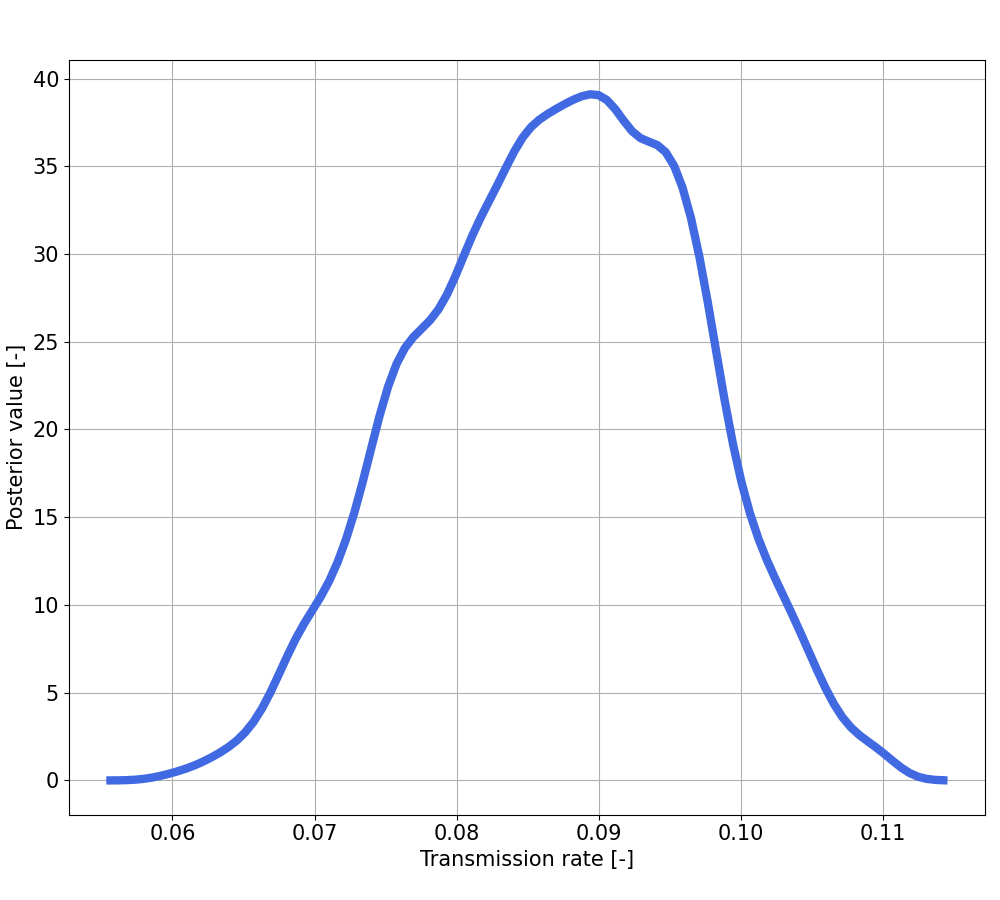}}\\
    \subfloat[$\beta$ at week 19.]{\includegraphics[scale=0.2]{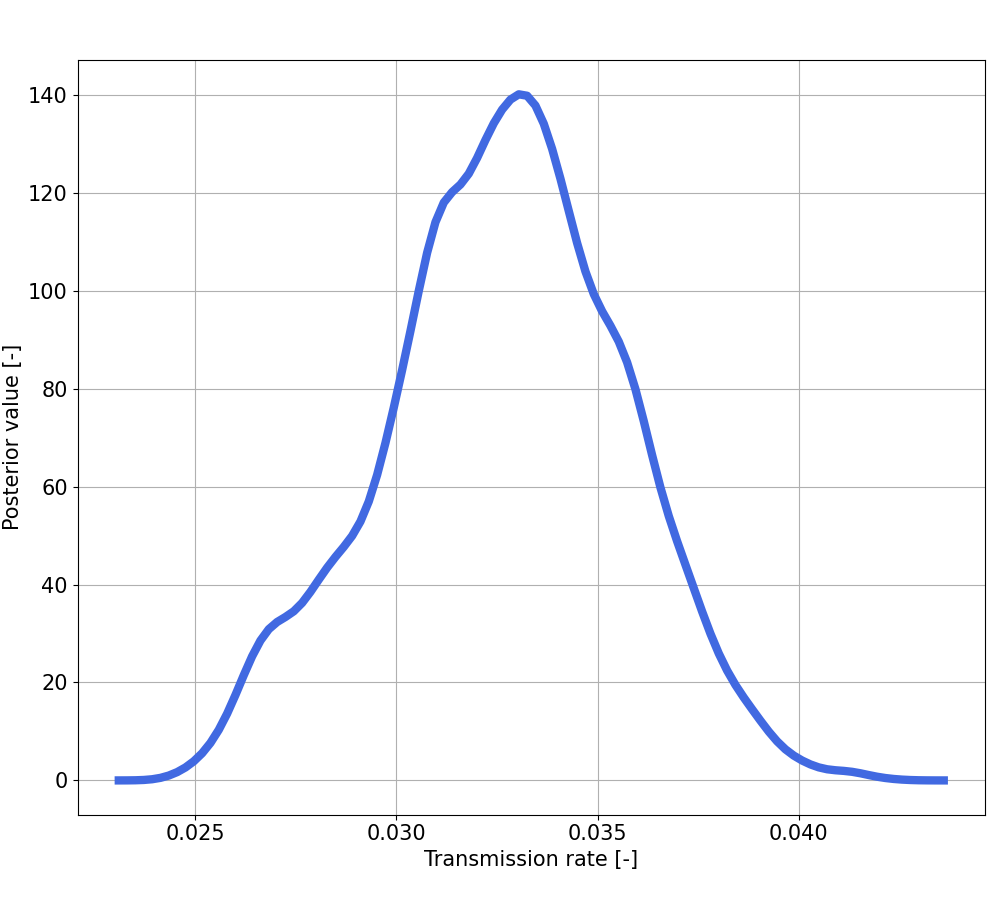}}
    \subfloat[Recovery time.]{\includegraphics[scale=0.2]{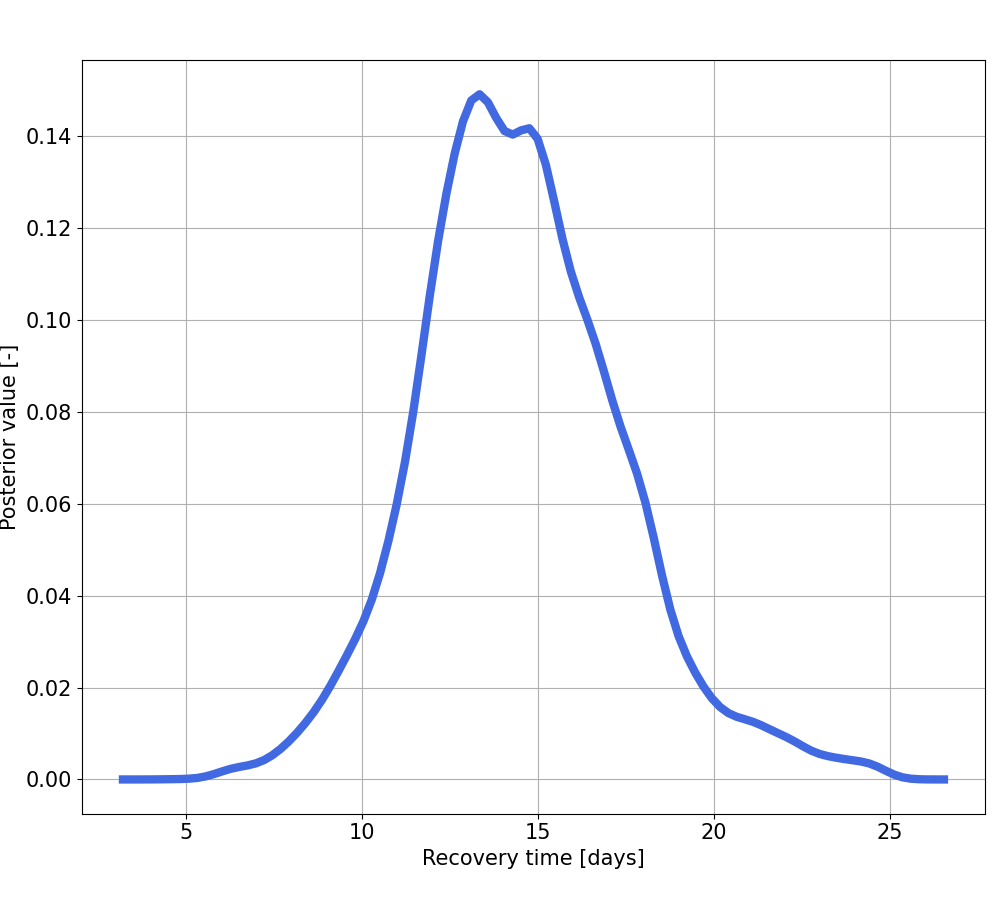}}
    \caption{Posterior distribution of three different transmission rates at the first (a), 9th (b) and 19th (c) week from the beginning of the simulation. (d) Posterior distribution of the recovery time variable.}
    \label{fig:my_label}
\end{figure}
\endminipage
\hspace{2mm}
\minipage{0.25\textwidth}
\begin{table}[H]
    \centering
    \begin{tabular}{c|c}
     \textbf{Week} & \textbf{$\beta$ median value}\\\hline
    1 & 0.02866168\\ \hline
    2 & 0.08487043\\ \hline
    3 & 0.05793916\\ \hline
    4 & 0.03441935\\ \hline
    5 & 0.04757512\\ \hline
    6 & 0.07276475\\ \hline
    7 & 0.08298755\\ \hline
    8 & 0.07913101\\ \hline
    9 & 0.07514123\\ \hline
    10 & 0.07197278\\ \hline
    11 & 0.07040295\\ \hline
    12 & 0.06946485\\ \hline
    13 & 0.06860226\\ \hline
    14 & 0.06789474\\ \hline
    15 & 0.06733712\\ \hline
    16 & 0.06685438\\ \hline
    17 & 0.06643925\\ \hline
    18 & 0.0661041\\ \hline
    19 & 0.03511503\\ \hline
    20 & 0.02103584\\ \hline
    21 & 0.01821423\\ 
    \end{tabular}
    \caption{Median values retrieved by the posterior distributions during each week of the calibration interval.}
    \label{transRateTab}
\end{table}
\endminipage


\begin{figure}[H]
    \centering
    \subfloat[Total deceased after the calibration stages compared with the DPC data function.]{\includegraphics[trim={0 0 0 2cm},clip,scale = 0.4]{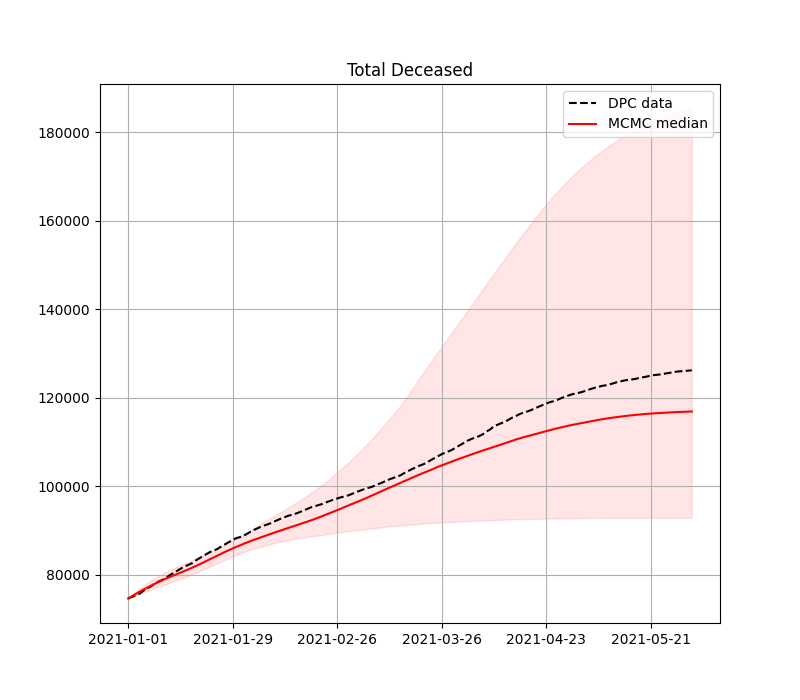}}\\
    \subfloat[Age-dependent evolution of deceased after the calibration stages compared with the DPC data \cite{dpcdata} function.]{\includegraphics[scale = 0.5]{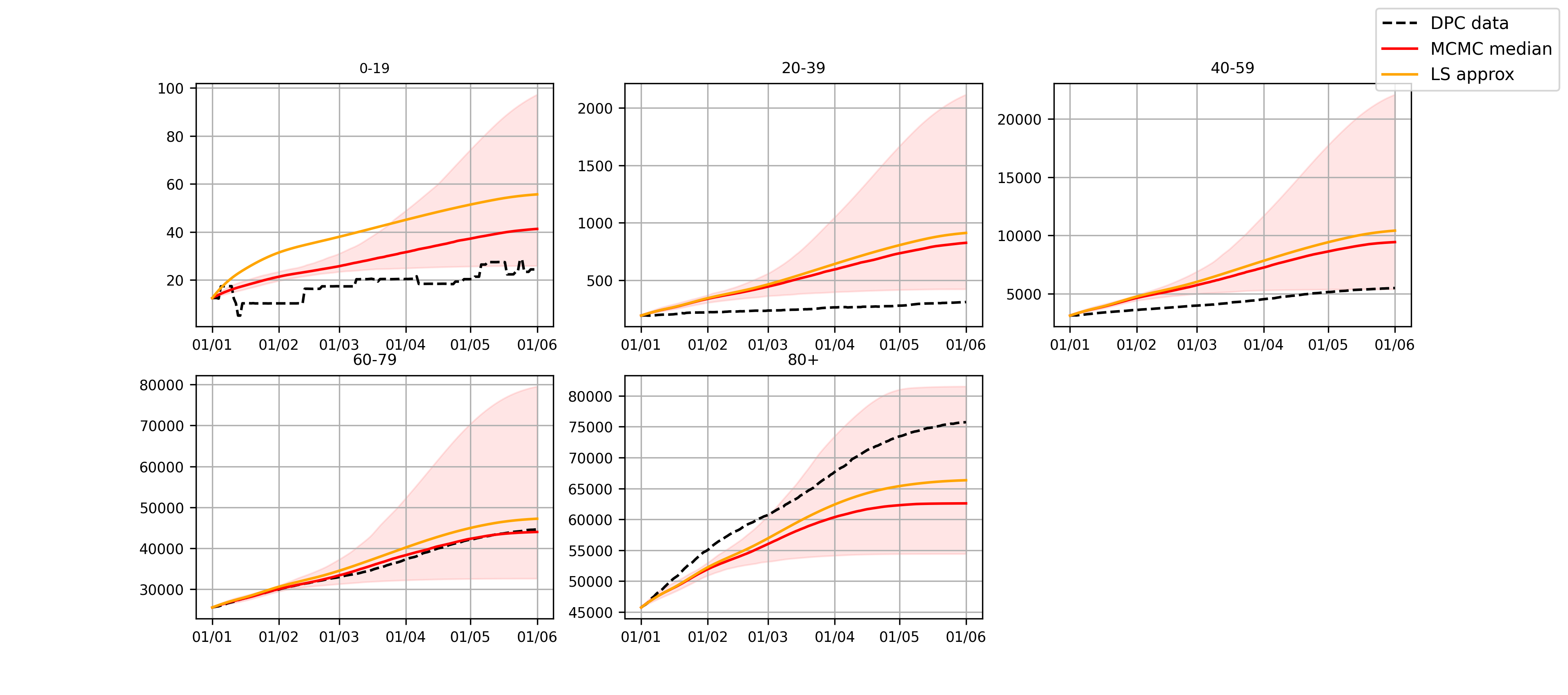}}\\
    \subfloat[Age-percentage repartition of deceased after the calibration stages compared with the DPC data \cite{dpcdata} function.]{\includegraphics[scale = 0.6]{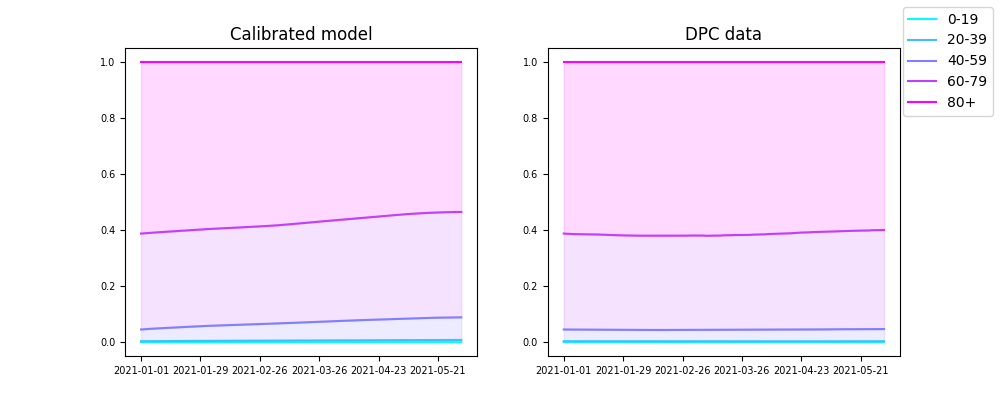}}
    \caption{Deceased evolutions after the calibration.}
    \label{calibDec}
\end{figure}

\begin{figure}[H]
    \centering
    \subfloat[Age-dependent evolution of detected infected after the calibration stages compared with the DPC data \cite{dpcdata} function.]{\includegraphics[scale = 0.5]{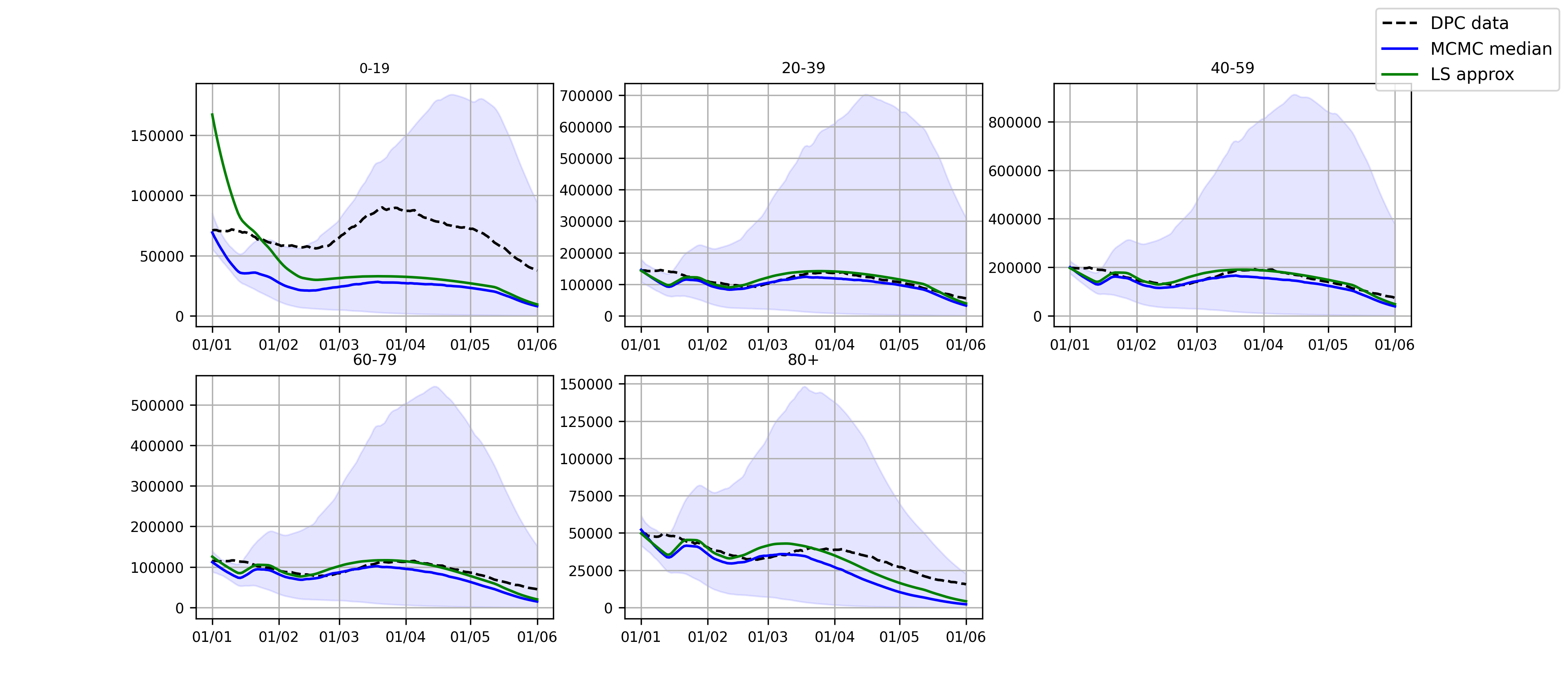}}\\
    \subfloat[Age-percentage repartition of deceased after the calibration stages compared with the DPC data \cite{dpcdata} function.]{\includegraphics[scale = 0.6]{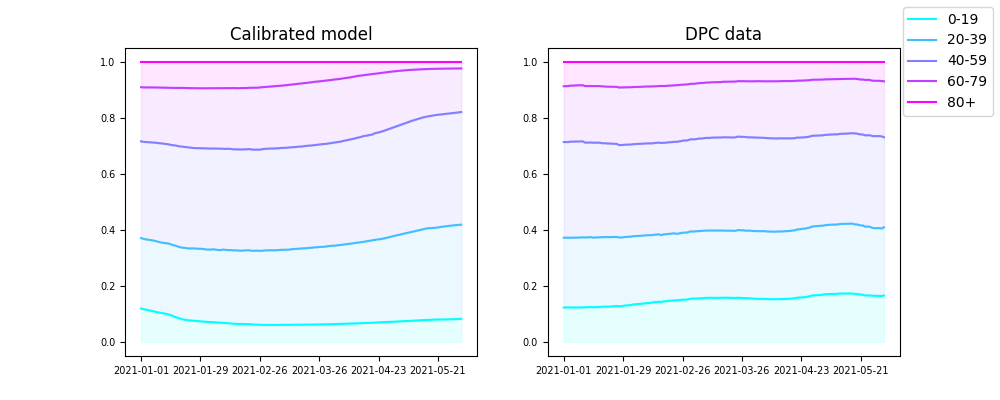}}
    \caption{Infected evolutions after the calibration.}
    \label{calibInf}
\end{figure}
\section{Sensitivity analyses}
\label{Sensitivity}
In this section we show the results of a sensitivity analysis of the reproduction number $\mathcal{R}_t$ with respect to the vaccination parameters $\theta_V,\theta_W,\sigma_V,\sigma_W$. {The reproduction number of the model, measuring the amount of secondary cases infected by a primal one in a fully susceptible population, has been recovered following the Next Generation Matrix approach \cite{diekmann2010construction}, \textit{i.e.} by computing the spectral radius of the so-called next generation matrix associated to the model}. In the following, we consider two different scenarios:
\begin{itemize}
    \item \textbf{Scenario 1 (Figure \ref{sensSigma})}: Consider the dependency of $\mathcal{R}_t$ on the effectiveness in reducing transmissibility, covering the space of admissibile parameters  $\sigma_{V}, \, \sigma_W \, \in [0, 1]$.
    \item \textbf{Scenario 2 (Figure \ref{sensTheta})}: Consider the dependency of $\mathcal{R}_t$ on the effectiveness in reducing severe infections, covering the space of admissibile parameters $\theta_{V}, \, \theta_W \, \in [0, 1]$.
\end{itemize}
It is well-known that the value of the reproduction number plays a key-role in the dynamics of the epidemic. Indeed, it can be interpreted as a bifurcation parameter assessing whether an outbreak is starting or if the epidemic is slowing down (threshold value for $\mathcal{R}_t$ is 1). {In all the simulations, the value of the transmission rate in each time phase is set equal to the median value obtained by the calibration process}, while the recovery time has been chosen equal  to the median value of the posterior distribution of the MCMC calibration, \textit{i.e.} 14.20 days. From Figures \ref{sensSigma}(a) and \ref{sensTheta}(a) we deduce that there is no evidence of the impact of the vaccination parameters on the initial Reproduction number (also known as $\mathcal{R}_0$).  {Moreover, from Figure \ref{sensTheta} we observe that variations of the severity reduction parameters $\theta_V,\theta_W$ have very little impact  (order 1e-9/1e-10)} on the reproduction number at different weeks. 
{This is not unexpected since these parameters enter in the model only through the fatality-function, thus impacting only the compartments $R$ and $D$. More precisely, as the reproduction number takes into account the amount of new infections and the infections cannot come from $D$ and $R$ compartments ($\mu_R$ is a very slow rate), $\theta_V,\theta_W$ are not relevant in the computation of $\mathcal{R}_t$.}
On the other hand, the results reported in Figure \ref{sensSigma} show that the variations in the effectiveness of the vaccine (parameters $\sigma_V$ and $\sigma_W$) have relevant impact (order 1e-2) on reducing the transmissibility after the 8th week of simulation, as a result of the increase in the amount of administrations.
{Looking at the slopes of the isolines in Figure \ref{sensSigma} we observe as the influence of $\sigma_V$ becomes more prominent over the one of $\sigma_W$ starting from week 10.} 
{Finally, a close inspection of the isolines in the neighbourhoods of the values retrieved by medical analyses (red dots in Figures \ref{sensSigma} and \ref{sensTheta}) reveals that fixing deterministically the values of the vaccination parameters, rather than dealing with its uncertainty, does not have a significant influence on the reproduction number.}
\begin{figure}[H]
    \centering
    \subfloat[Week 0.]{\includegraphics[scale = 0.32]{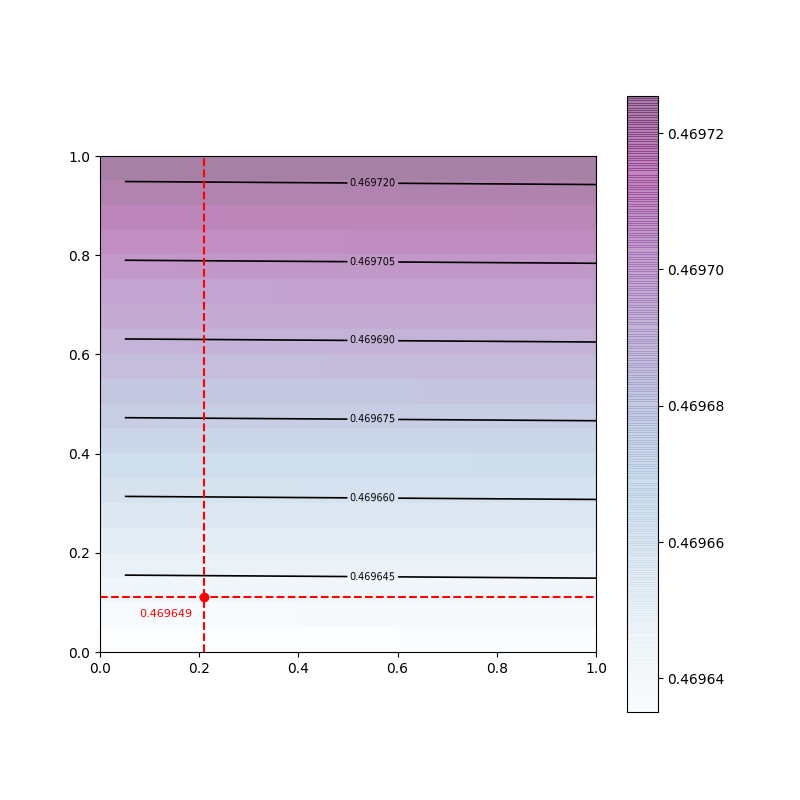}}
    \subfloat[Week 4.]{\includegraphics[scale = 0.32]{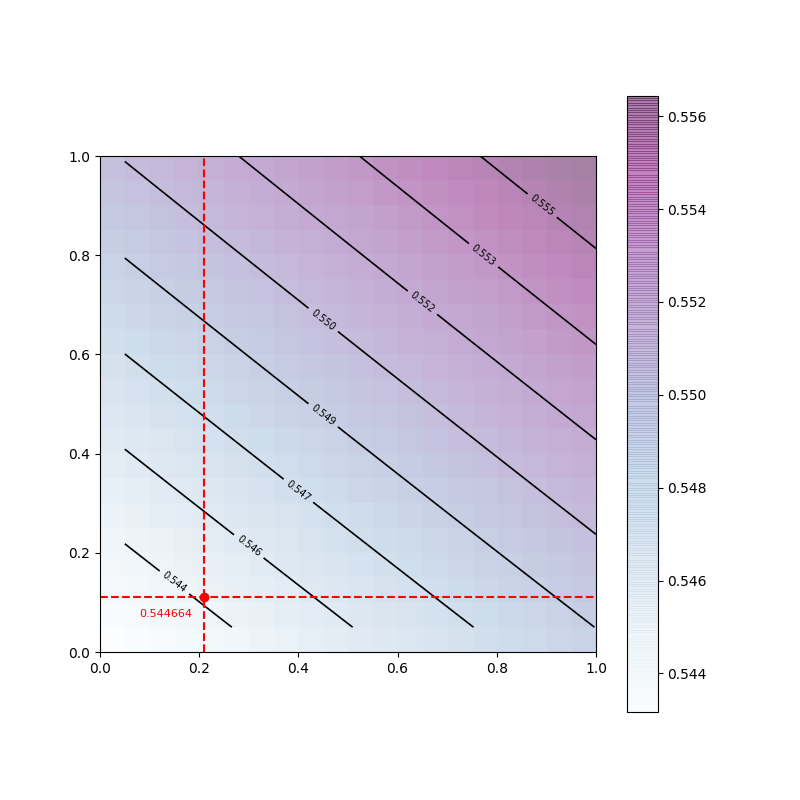}}\\
    \subfloat[Week 8.]{\includegraphics[scale = 0.32]{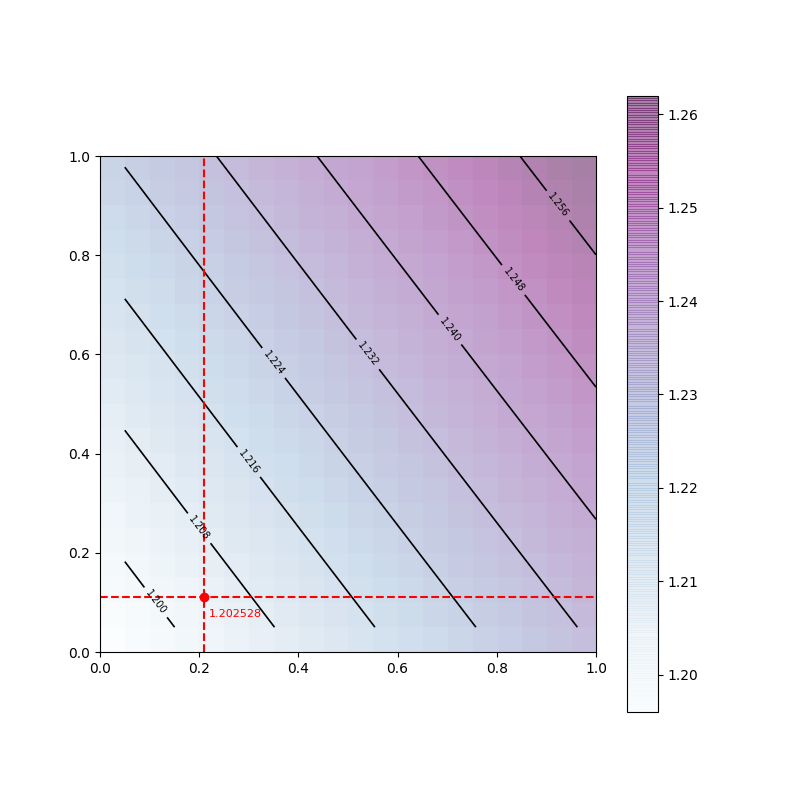}}
    \subfloat[Week 12.]{\includegraphics[scale = 0.32]{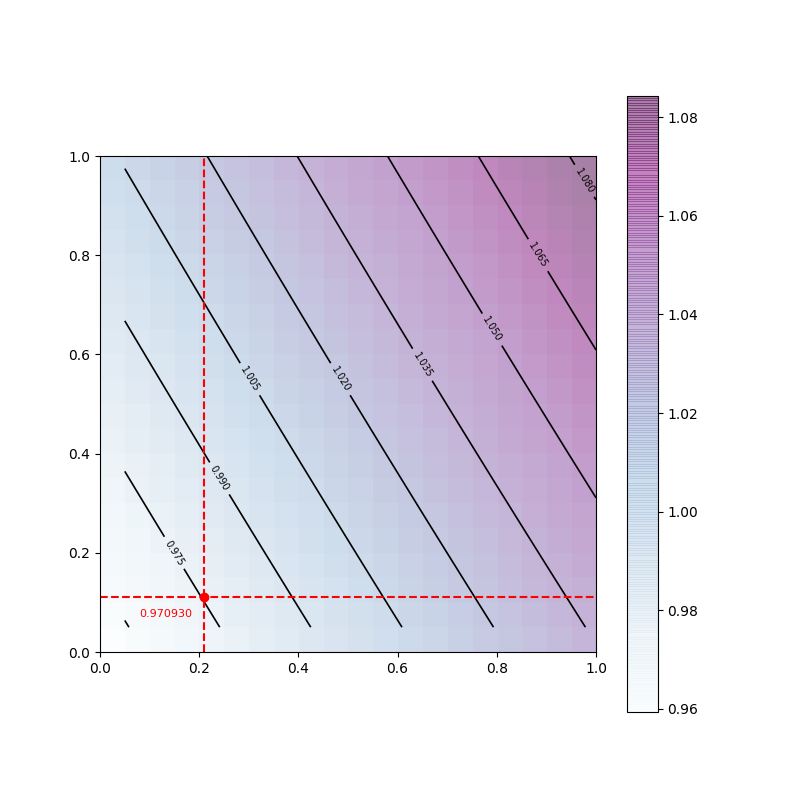}}\\
    \subfloat[Week 16.]{\includegraphics[scale = 0.32]{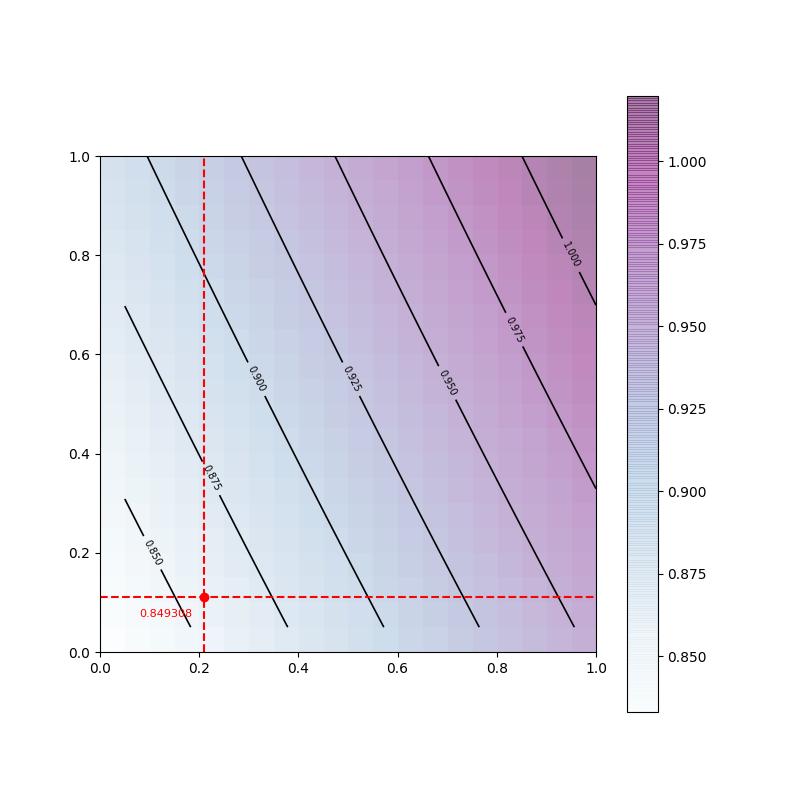}}
    \subfloat[Week 20.]{\includegraphics[scale = 0.32]{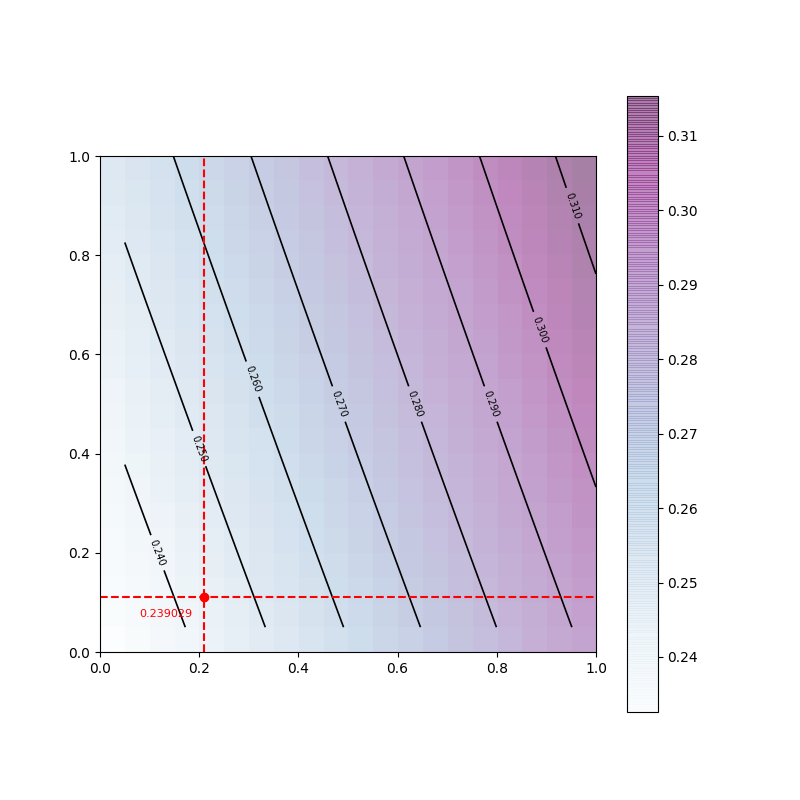}}\\
    \caption{Reproduction number at different weeks depending on the value assigned to the effectivenesses of doses on transmissibility. X-axis represents $\sigma_V$ and Y-axis represents $\sigma_W$. The red point corresponds to the reference value, \textit{i.e.} $\mathcal{R}_t$ value with $\sigma_V,\, \sigma_W$ as in Table \ref{ParamTab}.}
    \label{sensSigma}
\end{figure}

\begin{figure}[H]
    \centering
    \subfloat[Week 0.]{\includegraphics[scale = 0.32]{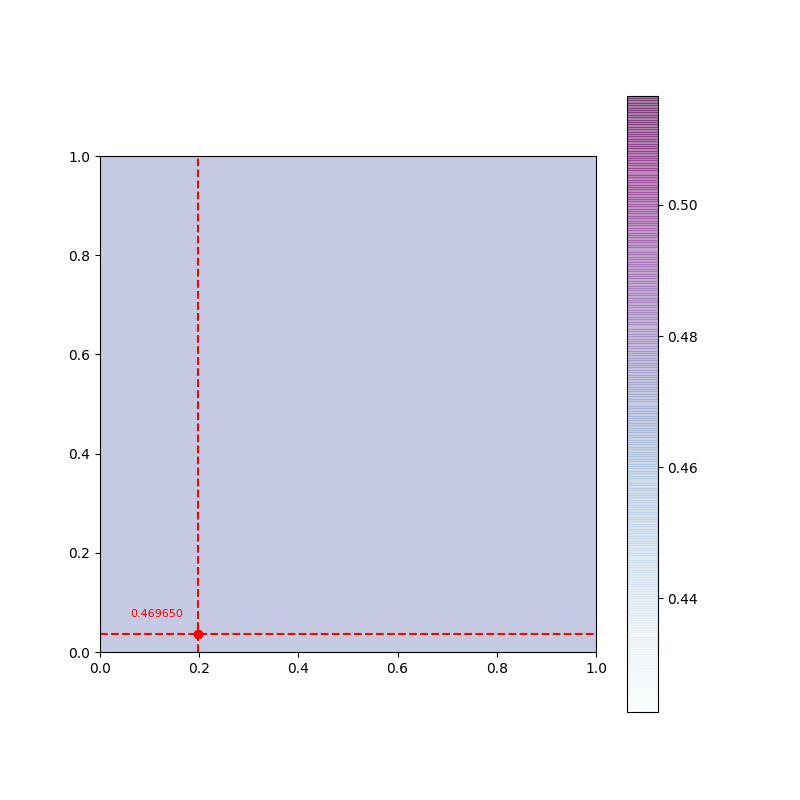}}
    \subfloat[Week 4. Reference $\mathcal{R}_t$ = 0.54480551917.]{\includegraphics[scale = 0.32]{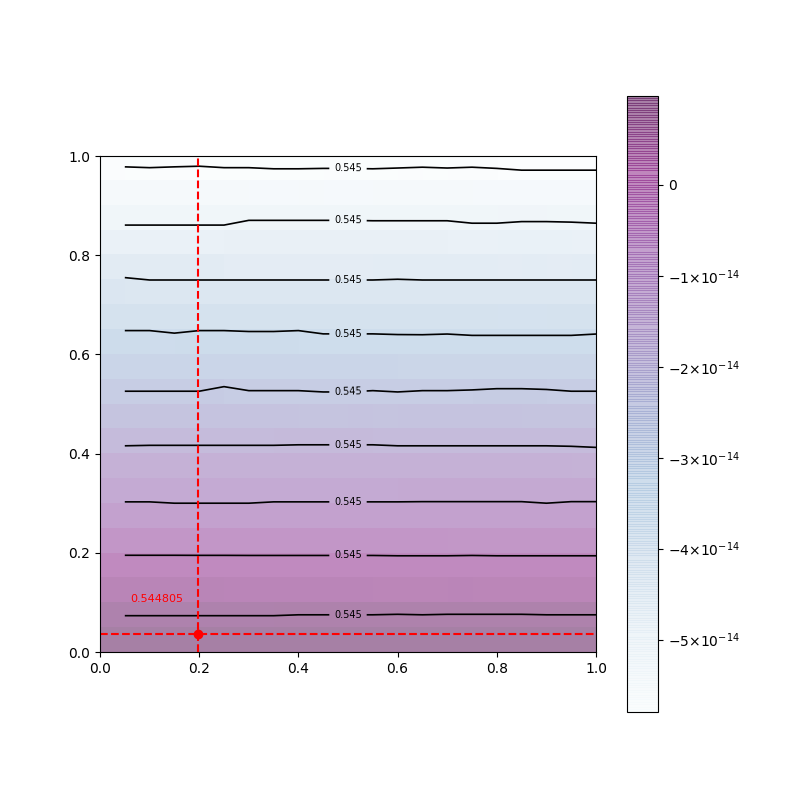}}\\
    \subfloat[Week 8. Reference $\mathcal{R}_t$ = 1.20322965963.]{\includegraphics[scale = 0.32]{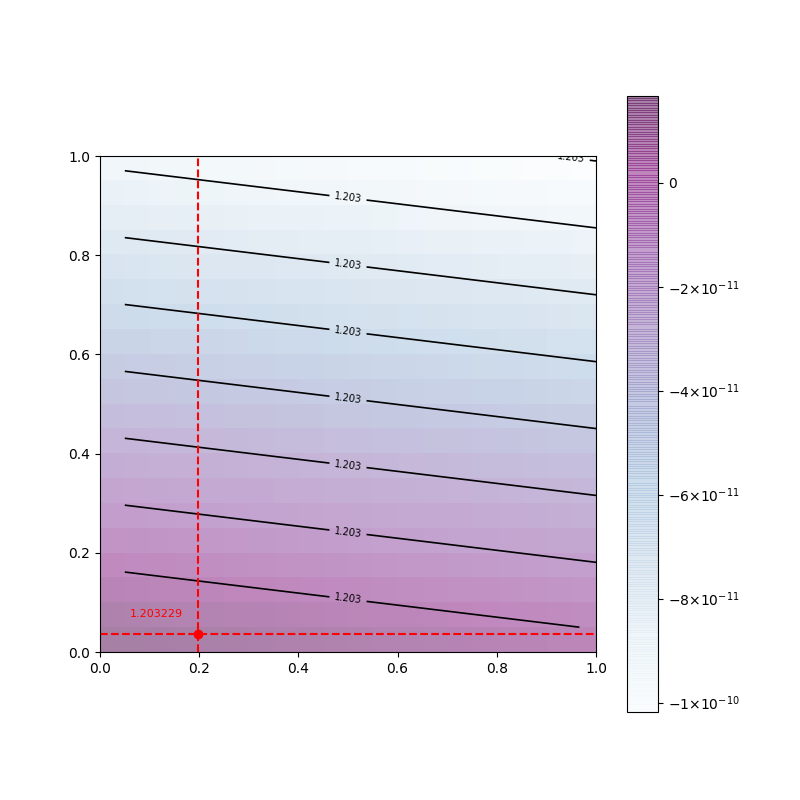}}
    \subfloat[Week 12. Reference $\mathcal{R}_t$ = 0.97225228190.]{\includegraphics[scale = 0.32]{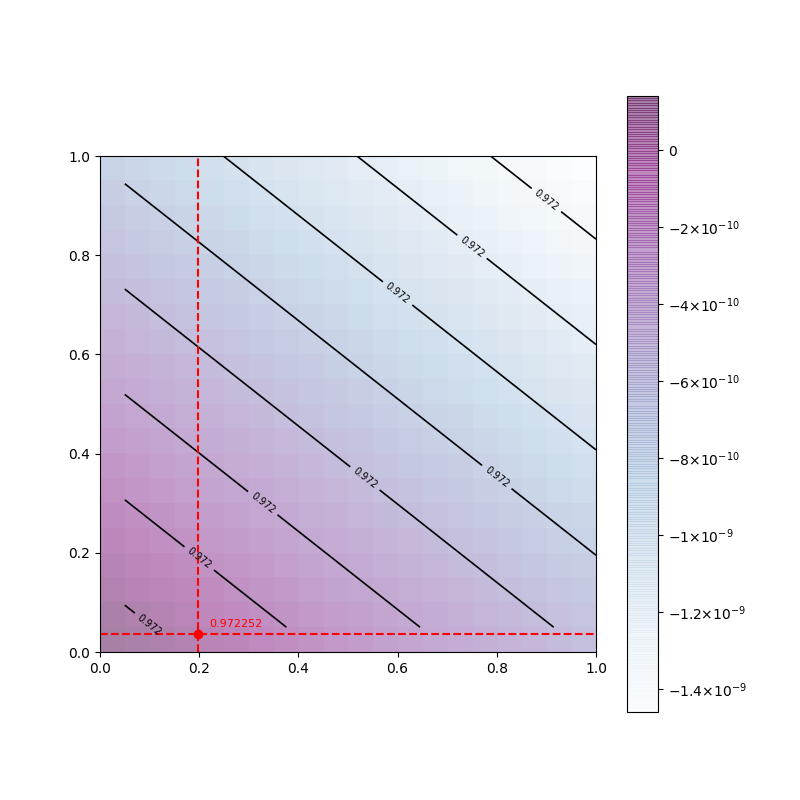}}\\
    \subfloat[Week 16. Reference $\mathcal{R}_t$ = 0.85123543297.]{\includegraphics[scale = 0.32]{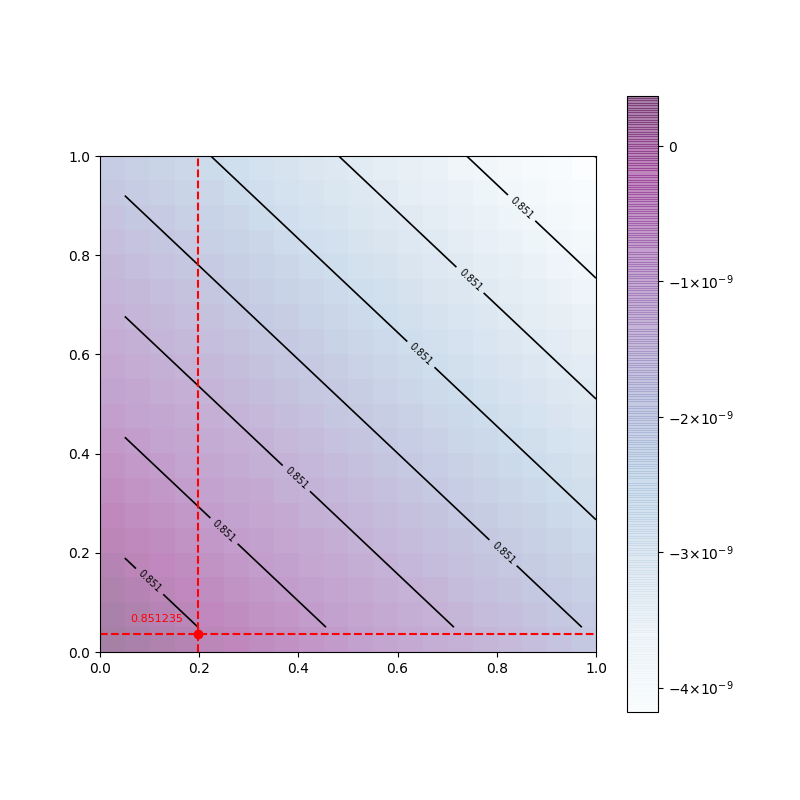}}
    \subfloat[Week 20. Reference $\mathcal{R}_t$ = 0.23986708323.]{\includegraphics[scale = 0.32]{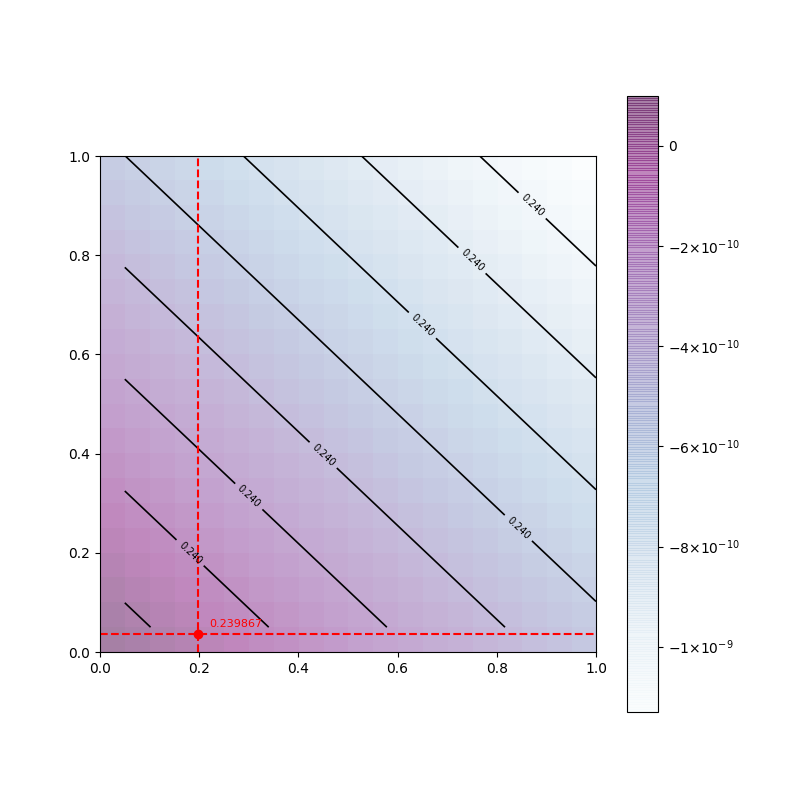}}\\
    \caption{Reproduction number at different weeks depending on the value assigned to the effectivenesses of doses on severity. X-axis represents $\theta_V$ and Y-axis represents $\theta_W$. The red point corresponds to the reference value, \textit{i.e.} $\mathcal{R}_t$ value with $\theta_V,\, \theta_W$ as in Table \ref{ParamTab}. \textit{Remark}: In Figures (b)-(f) the colorbar value stands for the increment or decrement with respect to the reference value.}
    \label{sensTheta}
\end{figure}

\printbibliography
\end{document}